\documentclass[gmd]{copernicus}

\usepackage{booktabs}
\newcommand{\epswitch}[2]{#1}
\begin{document}
\nolinenumbers

\title{asQ: parallel-in-time finite element simulations using ParaDiag for geoscientific models and beyond}

% \Author[affil]{given_name}{surname}

\Author[1]{Joshua}{Hope-Collins}
\Author[1,2]{Abdalaziz}{Hamdan}
\Author[3]{Werner}{Bauer}
\Author[4]{Lawrence}{Mitchell}
\Author[1]{Colin}{Cotter}

\affil[1]{Department of Mathematics, Imperial College London, London SW7 2AZ, UK}
\affil[2]{Institute for Mathematical Innovation, University of Bath, Bath, BA2 7AY, UK}
\affil[3]{School of Mathematics and Physics, University of Surrey, Guildford, GU2 7XH, UK}
\affil[4]{Independent researcher, Edinburgh, UK}
%% The [] brackets identify the author with the corresponding affiliation. 1, 2, 3, etc. should be inserted.

%% If an author is deceased, please mark the respective author name(s) with a dagger, e.g. "\Author[2,$\dag$]{Anton}{Smith}", and add a further "\affil[$\dag$]{deceased, 1 July 2019}".

%% If authors contributed equally, please mark the respective author names with an asterisk, e.g. "\Author[2,*]{Anton}{Smith}" and "\Author[3,*]{Bradley}{Miller}" and add a further affiliation: "\affil[*]{These authors contributed equally to this work.}".

\correspondence{Joshua Hope-Collins (joshua.hope-collins13@imperial.ac.uk)}

\runningtitle{asQ: parallel-in-time finite element simulations}

\runningauthor{Hope-Collins, Hamdan, Bauer, Mitchell, and Cotter}

%% These dates will be inserted by Copernicus Publications during the typesetting process.
\received{}
\pubdiscuss{} %% only important for two-stage journals
\revised{}
\accepted{}
\published{}

\firstpage{1}

\maketitle
\begin{abstract}
Modern high performance computers are massively parallel; for many PDE applications spatial parallelism saturates long before the computer's capability is reached.
Parallel-in-time methods enable further speedup beyond spatial saturation by solving multiple timesteps simultaneously to expose additional parallelism. ParaDiag is a 
particular approach to parallel-in-time based on preconditioning the simultaneous timestep system with a perturbation that allows block diagonalisation via a Fourier transform in time.
In this article, we introduce asQ, a new library for implementing ParaDiag parallel-in-time methods, with a focus on applications in the geosciences, especially weather and climate.
asQ is built on Firedrake, a library for the automated solution of finite element models, and the PETSc library of scalable linear and nonlinear solvers.
This enables asQ to build ParaDiag solvers for general finite element models and provide a range of solution strategies, making testing a wide array of problems straightforward.
We use a quasi-Newton formulation that encompasses a range of ParaDiag methods, and expose building blocks for constructing more complex methods.
The performance and flexibility of asQ is demonstrated on a hierarchy of linear and nonlinear atmospheric flow models.
We show that ParaDiag can offer promising speedups and that asQ is a productive testbed for further developing these methods.
\end{abstract}

%\copyrightstatement{TEXT} %% This section is optional and can be used for copyright transfers.

\tableofcontents

\introduction\label{sec:intro}  %% \introduction[modified heading if necessary]
In this article, we present asQ, a software framework for investigating the performance of ParaDiag parallel-in-time methods.
We focus our attention on geophysical fluid models relevant to the development of simulation systems for oceans, weather and climate.
This library allows researchers to rapidly prototype implementations of ParaDiag for time dependent partial differential equations (PDEs) discretised using Firedrake \citep{FiredrakeUserManual}, selecting options for the solution strategy facilitated by the composable design of PETSc, the Portable, Extensible Toolkit for Scientific Computation \citep{petsc-user-ref}.
\emph{The goal of the paper is not to advocate for ParaDiag as superior to other methods, but to demonstrate that asQ serves this purpose.}

Parallel-in-time (PinT) is the name for the general class of methods that introduce parallelism in the time direction as well as in space.
The motivation for PinT methods is that eventually it is not possible to achieve acceptable time to solution only through spatial domain decomposition as one moves to higher and higher fidelity solutions, and so one would need to look to the time dimension for further speedups.
\citet{friedhoff2012multigrid} provided a quantitative argument that any sufficiently high resolution time dependent simulation will eventually require the use of time parallel methods (the question is just when, and how). 

PinT methods have a long history, surveyed by \citet{gander201550}, but the topic has really exploded since the late 1990s, with many algorithms being proposed including space-time-concurrent multigrid waveform relaxation (WRMG) \citep{vandewalle1994space}, space time multigrid \citep{horton1995space}, Parareal \citep{maday2002parareal}, revisionist integral deferred correction (RIDC) \citep{christlieb2010parallel}, multigrid reduction in time (MGRIT) \citep{friedhoff2012multigrid}, parallel full approximation scheme in space and time (PFASST) \citep{emmett2012toward},   ParaEXP \citep{gander2013paraexp}, and, the subject of this work, ParaDiag.
As we shall elaborate later with references, ParaDiag is a computational linear algebra approach to PinT, solving an implicit system for several timesteps at once, using Fourier transforms in time to obtain a block diagonal system whose components can be solved independently and in parallel.
A review of the most common forms of ParaDiag and important analysis results can be found in \citet{gander_paradiag_2021}.

Although there are a small number of software libraries for PinT methods, most PinT software implementations are written from scratch as small standalone codes or individual scripts; increased availability and use of PinT libraries could increase research productivity \citep{pintsoftware}.
At the time of writing, three PinT libraries are mature and available open-source.
XBraid \citep{xbraid-package} and pySDC \citep{Speck2019} are reference implementations of MGRIT and spectral deferred corrections, written in C and Python, respectively.
These frameworks are designed to be non-intrusive so users can plug in existing serial-in-time code to quickly experiment with these PinT methods.
SWEET \citep{sweet-package} is a testbed for time-integration of the shallow water equations, a model of geophysical flow.
As opposed to implementing a single family of methods for any problem, SWEET implements many methods for a specific class of problem.
In comparison to these libraries, asQ implements a particular class of method (ParaDiag), for a particular class of discretisation (finite elements) but for general PDEs.

For the purposes of later discussion, we briefly define two scaling paradigms when exploiting parallelism.
Strong scaling is where a larger number of processors is used to obtain the \textit{same} solution in a \textit{shorter} wallclock time.
Weak scaling is where a larger number of processors is used to obtain a \textit{higher resolution} solution in the \textit{same} wallclock time.

The rest of this article is structured as follows. In Sect. \ref{subsec:pint_hyperbolic} we survey the varying PinT approaches to geophysical fluids models in particular, which incorporate transport and wave propagation processes that exhibit the general challenge of PinT methods for hyperbolic problems (and hence we also discuss aspects of hyperbolic problems more broadly). In Sect. \ref{subsec:prior_paradiag}, we complete this introduction with a survey previous research on the ParaDiag approach to PinT. Then, in Sect. \ref{sec:methods} we review the basic ParaDiag idea and discuss the extension to nonlinear problems 
which are more relevant to weather and climate, highlighting the wide range of choices that need to be made when using ParaDiag for these problems. In Sect. \ref{sec:asq}, we describe the asQ library, and explain how it addresses the need to rapidly explore different options in a high performance computing environment. In Sect. \ref{sec:examples} we present some numerical examples to demonstrate that we have achieved this goal. Finally, in Sect. \ref{sec:conclusions} we provide a summary and outlook.

\subsection{Parallel-in-time methods for hyperbolic and geophysical models}
\label{subsec:pint_hyperbolic}

The potential for PinT methods in oceans, weather and climate simulation is attractive because of the drive to higher resolution.
For example, the Met Office Science Strategy \citep{met2022}, and \citet{Bauer_Thorpe_Brunet_2015} highlight the need for global convection permitting atmosphere models, eddy resolving ocean models, eddy permitting local area atmosphere models, and estuary resolving shelf-seas models, to better predict hazards and extremes, which will require sub 10km global resolution.
In operational forecasting, the model needs to run to a particular end time (e.g.~ten simulation days) within a particular wallclock time in order to complete the forecasting procedure in time for the next cycle (e.g.~three wallclock hours).
Similarly, climate scientists have a requirement for simulations to complete within a feasible time for a model to be scientifically useful.
To try to maintain the operational wallclock limit when resolution is increased, weak scaling is used so that each timestep at the higher spatial resolution can be completed in the same wallclock time as the timestep at the lower spatial resolution.
However, even when this weak scaling is achievable, high spatial resolution yields dynamics (e.g.~transport and waves) with higher temporal frequencies that should be resolved in the timestep (and sometimes we are forced to resolve them due to stability restrictions in timestepping methods). This means that more timesteps are required at the higher resolution than the lower resolution.
To satisfy the operational wallclock limit we now need to be able to strong scale the model to reduce the wallclock time for each timestep, so that the same simulation end time can be reached without breaking the wallclock time limit at the higher resolution.
Achieving this scaling with purely spatial parallelism is very challenging because these models are already run close to the scaling limits. This motivates us to consider other approaches to timestepping such as PinT methods.

The challenge to designing effective PinT methods for these geophysical fluid dynamics models is that their equations support high frequency wave components coupled to slow balanced motion that governs the large scale flow, such as the fronts, cyclones, jets and Rossby waves that are the  familiar features of midlatitude atmospheric weather.
The hyperbolic nature of these waves makes them difficult to treat efficiently using the classical Parareal algorithm, as discussed in \citet{gander201550}.
The difficulty is that the errors are dominated by dispersion error, and there is a mismatch in the dispersion relation between coarse and fine model \citep{ruprecht2018wave}.
Similarly, \citet{De_Sterck_Friedhoff_Krzysik_MacLachlan_2023} showed that standard MGRIT has deteriorated convergence for hyperbolic problems due to the removal of some characteristic components on the coarse grid if simple rediscretisation is used.
However, \citet{de2023fast} showed that a carefully modified semi-Lagrangian method can overcome this deficiency by ensuring that the coarse grid operator approximates the fine grid operator to a higher order of accuracy than it approximates the PDE, analogously to previous findings for MGRIT applied to chaotic systems \citep{Vargas_Falgout_Gunther_Schroder_2023}.
Using this approach, \citet{De_Sterck_Falgout_Krzysik_Schroder_2023,de2023fast} demonstrated real speedups for the variable coefficient scalar advection.
Scalable iteration counts were obtained for nonlinear PDEs using a preconditioned quasi-Newton iteration \citep{De_Sterck_Falgout_Krzysik_Schroder_2024a} and systems of linear and nonlinear PDEs \citep{De_Sterck_Falgout_Krzysik_Schroder_2024b}.
\citet{hamon2020parallel} and \citet{Caldas_Steinstraesser_Peixoto_Schreiber_2024} have demonstrated parallel speedups for the nonlinear rotating shallow water equations (a prototypical highly oscillatory PDE for geophysical fluid dynamics), using multilevel methods. 

For linear systems with pure imaginary eigenvalues (e.g. discretisations of wave equations), a parallel technique based on sums of rational approximations of the exponential function (referred to in later literature as rational exponential integrators (REXI)) restricted to the imaginary axis was proposed in \citet{haut2016high}.
The terms in the sum are mutually independent so can be evaluated in parallel.
Each of these terms requires the solution of a problem that resembles a backward Euler integrator with a complex valued timestep.
For long time intervals in the wave equation case, some of these problems resemble shifted Helmholtz problems with coefficients close to the negative real axis and far from the origin.
These problems are known to be very unsuited to be solved by multigrid methods which are otherwise a scalable approach for linear timestepping problems \citep{gander2015applying}.
\cite{haut2016high} avoided this by using static condensation techniques for higher order finite element methods.
This approach was further investigated and developed in the geophysical fluid dynamics setting using pseudospectral methods (including on the sphere) in \citet{schreiber2018beyond, schreiber2019parallel}; \citet{schreiber2019exponential} used a related approach with coefficients derived from Cauchy integral methods.
In an alternative direction, Paraexp \citep{gander2013paraexp} provides a PinT mechanism for dealing with nonzero source terms for the linear wave equation, if a fast method for applying the exponential is available.

Multiscale methods are a different strategy to tackle the highly oscillatory components, whose phases are not tremendously important, but their bulk coupling to the large scale can be.
\citet{legoll2013micro} proposed a micro-macro Parareal approach where the coarse propagator is obtained by averaging the vector field over numerical solutions with frozen macroscopic dynamics, demonstrating parallel speedups for test problems using highly oscillatory ordinary differential equations.
\citet{haut2014asymptotic} proposed a different approach, based on previous analytical work \citep{schochet1994fast,embid1998low,majda1998averaging}, in which the highly oscillatory PDE is transformed using operator exponentials to a nonautonomous PDE with rapidly fluctuating explicit time dependence.
After averaging over the phase of this explicit time dependence, a slow PDE is obtained that approximates the transformed system under suitable assumptions.
To obtain a numerical algorithm, the "averaging window" (range of phase values to average over) is kept finite, and the average is replaced by a sum whose terms can be evaluated independently in parallel (providing parallelisation across the method for the averaged PDE).
\citet{haut2014asymptotic} used this approach to build a coarse propagator for the one dimensional rotating shallow water equations.
In experiments with a standard geophysical fluid dynamics test case on the sphere, \citet{yamazaki2023time} showed that the error due to averaging can actually be less than the time discretisation error in a standard method with the same timestep size, suggesting that phase averaging might be used as a PinT method in its own right without needing a Parareal iteration to correct it.
\citet{bauer2022higher} showed an alternative route to correcting the phase averaging error using a series of higher order terms, which may expose additional parallelism. 

\subsection{Prior ParaDiag research}

\label{subsec:prior_paradiag}

The software we present here is focused on $\alpha$-circulant diagonalisation techniques for all at once systems, which have come to be known in the literature as ``ParaDiag''.
In this class of methods, a linear constant coefficient ODE (e.g.~a discretisation in space of a linear constant coefficient PDE) is discretised in time, resulting in a system of equations with a tensor product structure in space and time.
This system is then diagonalised in time, leading to a block diagonal system with one block per timestep, which can each be solved in parallel before transforming back to obtain the solution at each timestep.
The first mathematical challenge is that the  block structure in time is not actually diagonalisable for constant timestep $\Delta t$, because of the nontrivial Jordan normal form.
In the original paper on diagonalisation in time, \citet{maday2008parallelization} tackled this problem by using timesteps forming a geometric progression, which then allows for a direct diagonalisation.
The main drawback is that the diagonalisation is badly conditioned for small geometric growth rate, which might otherwise be required for accurate solutions.
\citet{mcdonald2018preconditioning} proposed to use a time periodic (and thus diagonalisable) system to precondition the initial value system.
This works well for parabolic systems but is not robust to mesh size for wave equations (such as those arising in geophysical fluid applications).
\citet{gander2019convergence} (following \citet{wu_toward_2018} in the parareal context) proposed a modification, in which the periodicity condition $u(0)=u(T)$ is replaced by a periodicity condition $u^{k}(0)=u_{0} + \alpha(u^{k}(T)-u^{k-1}(T))$ in the preconditioner, with $u_{0}$ the real initial condition, $k$ the iteration index, and $0<\alpha<1$, with the resulting time block structure being called $\alpha$-circulant.
This system can be diagonalised by FFT in time after appropriate scaling by a geometric series.
When $\alpha<\frac{1}{2}$, an upper bound can be shown for the preconditioner which is independent of the mesh, the linear operator, and any parameters of the problem.
In particular, it produces mesh independent convergence for the wave equation. However, the diagonalisation is badly conditioned in the limit $\alpha\to 0$ \citep{gander2019convergence}.
\citet{vcaklovic2023parallel} addressed this by adapting \(\alpha\) at each iteration.

The technique has also been extended to other timestepping methods.
\citet{vcaklovic2023parallel} provided a general framework for higher order implicit collocation methods, using the diagonalisation of the polynomial integral matrix.
The method can also be extended to higher order multistep methods, such as BDF methods \citep{danieli2021all,gander2024new} or Runge-Kutta methods \citep{wu2021parallel,Kressner_Massei_Zhu_2023}.

Moving to nonlinear PDEs, the all at once system for multiple timesteps must now be solved using a Newton or quasi-Newton method.
The Jacobian system is not constant coefficient in general, and so must be approximated by some form of average, as first proposed by \citet{gander2017time}.
There are a few analytical results about this approach.
\citet{gander2019convergence} proved convergence for a fixed point iteration for the nonlinear problem when an $\alpha$-periodic time boundary condition is used.
\citet{caklovic_paradiag_2023} developed a theory for convergence of quasi-Newton methods, presented later in \eqref{eq:nonlinear-convergence-bound}.

Performance measurements for time-parallel ParaDiag implementations are still relatively sparse in the literature, and have been predominantly for linear problems with small scale parallelism.
\citet{goddard_note_2019}, \citet{gander_paradiag_2021} and \citet{Liu_Wang_Wu_Zhou_2022} found very good strong scaling up to 32, 128 and 256 processors respectively for the heat, advection, and wave equations, with \citet{Liu_Wang_Wu_Zhou_2022} noting the importance of a fast network for multi-node performance due to the collective communications required in ParaDiag.
Actual speedups vs a serial-in-time method of \(15\times\) were obtained by \citet{liu_fast_2020} for the wave equation on 128 processors.
A doubly time-parallel ParaDiag implementation was presented in \citet{vcaklovic2023parallel} for the collocation method with parallelism across both the collocation nodes and the timesteps.
For the heat and advection equations, they achieved speedups of 15-20\(\times\) over the serial-in-time/serial-in-space method on 192 processors, and speedups of 10\(\times\) over the serial-in-time/parallel-in-space method on 2304 processors (with a speedup of 85\(\times\) over the serial-in-time/serial-in-space method).
As far as the authors are aware, \citet{caklovic_paradiag_2023} is the only study showing speedups vs a serial-in-time method for nonlinear problems, achieving speedups of 25\(\times\) for the Allen-Cahn equations and speedups of 5.4\(\times\) for the hyperbolic Boltzmann equation.

\section{ParaDiag methods}\label{sec:methods}
In this section we review the ParaDiag method. First, we discuss application to linear problems for which it was originally developed. Second, we discuss the adaptation to nonlinear problems which are of interest in many practical applications.
We then use a simple performance model to highlight the requirements for ParaDiag to be an effective method.

ParaDiag is a method to accelerate the solution of an existing time-integrator, so we will define the serial-in-time method first.
In our exposition, we are interested in time-dependent (and initially linear with constant coefficient) PDEs combined with a finite element semi-discretisation in space,
\begin{equation}\label{eq:linear-pde}
    M\partial_{t}u + Ku = b(t),
\end{equation}
where \(u\in\mathbb{R}^{N_{x}}\) is the unknown, \(M\in\mathbb{R}^{N_{x}\times N_{x}}\) is the mass matrix, \(K\in\mathbb{R}^{N_{x}\times N_{x}}\) is the stiffness matrix arising from the discretisation of the spatial terms, \(b(t)\in\mathbb{R}^{N_{x}}\) is some forcing term not dependent on the solution, and \(N_{x}\) is the number of spatial degrees of freedom (DoFs).
Throughout Sects. \ref{sec:methods:linear} and \ref{sec:methods:nonlinear}: lower case Roman letters denote vectors (except \(t\), which is reserved for time, and \(n\) and \(j\) which are reserved for integers); upper case Roman letters denote matrices (except \(N\), which is reserved for integers); vectors and matrices defined over both space and time are in boldface; and Greek letters denote scalars.
(\ref{eq:linear-pde}) is solved on a spatial domain \(x\in\Omega\) with appropriate boundary conditions on the boundary \(\partial\Omega\).
We want to find the solution in the time range \(t\in\left[t_{0}, t_{0}+N_{t}\Delta t\right]\) from an initial condition \(u_{0} = u(t_{0})\) at time \(t_{0}\).
To achieve this, (\ref{eq:linear-pde}) is discretised in time using the implicit \(\theta\)-method,
\begin{equation}\label{eq:theta-method}
    \frac{1}{\Delta t}M\left(u^{n+1} - u^{n}\right)
    + \theta K u^{n+1} + \left(1-\theta\right)K u^{n}
    = \tilde{b}^{n+1},
\end{equation}
where the right hand side is
\begin{equation}\label{eq:theta-rhs}
    \tilde{b}^{n+1} = \theta b^{n+1} + \left(1-\theta\right) b^{n},
\end{equation}
\(\Delta t\) is the timestep size, \(n\in[0,N_{t}]\) is the timestep index, \(u^{n}\) is the discrete solution at time \(t=t_{0}+n\Delta t\), and \(\theta\in[0,1]\) is a parameter.
ParaDiag is not restricted to the \(\theta\)-method, but this is currently the only method implemented in asQ (other methods may be added in the future).

As written, (\ref{eq:theta-method}) is an inherently serial method, requiring the solution of an implicit system for each timestep \(u^{n+1}\) given \(u^{n}\).

\subsection{Linear problems}\label{sec:methods:linear}
The ParaDiag method begins by constructing a single monolithic system for multiple timesteps, called the ``all-at-once system''.
We illustrate an all-at-once system below for four timesteps (\(N_{t}=4\)), created by combining equation (\ref{eq:theta-method}) for \(n = 0,1,2,3\), 
to obtain
\begin{equation}\label{eq:linear-aaos}
    \mathbf{A}\vec{u}
    = \left(B_{1}\otimes M
    + B_{2}\otimes K\right)\vec{u}
    = \vec{\tilde{b}},
\end{equation}
where \(\otimes\) is the Kronecker product.
The all-at-once matrix \(\mathbf{A}\in\mathbb{R}^{N_{t}N_{x}\times N_{t}N_{x}}\) (which is the Jacobian of the full all-at-once residual \(\mathbf{A}\vec{u}-\vec{\tilde{b}}\)) is written using Kronecker products of the mass and stiffness matrices with the two matrices \(B_{1,2}\in\mathbb{R}^{N_{t}\times N_{t}}\).
These are Toeplitz matrices which define the time-integrator,
\begin{align}\label{eq:toeplitz-matrices}
    B_{1} = &
    \frac{1}{\Delta t}
    \begin{pmatrix}
        1  &  0 &  0 &  0 \\
        -1 &  1 &  0 &  0 \\
        0  & -1 &  1 &  0 \\
        0  &  0 & -1 &  1 \\
    \end{pmatrix},\\
    \label{eq:toeplitz-B2}
    B_{2} = &
    \begin{pmatrix}
        \theta  &  0 &  0 &  0 \\
        (1-\theta) &  \theta &  0 &  0 \\
        0  & (1-\theta) &  \theta &  0 \\
        0  &  0 & (1-\theta) &  \theta \\
    \end{pmatrix},
\end{align}
and the vector \(\vec{u}\in\mathbb{R}^{N_{t}N_{x}}\) of unknowns for the whole time-series is
\begin{equation}\label{eq:aaofunction}
    \vec{u} =
    \begin{pmatrix}
        u^{1} \\ u^{2} \\ u^{3} \\ u^{4}
    \end{pmatrix}.
\end{equation}
The right hand side \(\vec{\tilde{b}}\in\mathbb{R}^{N_{t}N_{x}}\) includes the initial conditions,
\begin{equation}\label{eq:aaorhs}
    \vec{\tilde{b}} =
    \begin{pmatrix}
        \tilde{b}^{1} \\
        \tilde{b}^{2} \\
        \tilde{b}^{3} \\
        \tilde{b}^{4} \\
    \end{pmatrix}
    +
    \begin{pmatrix}
        \left(\frac{1}{\Delta t}M - (1-\theta)K\right)u^{0} \\
        0 \\
        0 \\
        0 \\
    \end{pmatrix}.
\end{equation}
Due to the properties of the Kronecker product, if \({B}_{1}\) and \({B}_{2}\) are simultaneously diagonalisable then the Jacobian \(\mathbf{A}\) is block diagonalisable.
A block-diagonal matrix can be efficiently inverted by solving each block separately in parallel.
Note that when forming an \(\mathbb{R}^{N_{t}N_{x}\times N_{t}N_{x}}\) space-time matrix using a Kronecker product, the \(\mathbb{R}^{N_{t}\times N_{t}}\) matrix is always on the left of the Kronecker product, and the \(\mathbb{R}^{N_{x}\times N_{x}}\) matrix is always on the right.

The original ParaDiag-I method introduced by \cite{maday2008parallelization} premultiplied (\ref{eq:linear-aaos}) by \({B}_{2}^{-1}\otimes I_{x}\) where \(I_{x}\) is the spatial identity matrix, and chose the time discretisation such that \(B_{2}^{-1}B_{1}\) is diagonalisable.
This gives a direct solution to (\ref{eq:linear-aaos}) but diagonalisation is only possible if the timesteps are distinct; \cite{maday2008parallelization} and \cite{gander2017time} used geometrically increasing timesteps which leads to large values of \(\Delta t\) and poor numerical conditioning for large \(N_{t}\).

Here, we focus on an alternative approach, named ParaDiag-II.
The Toeplitz matrices \eqref{eq:toeplitz-matrices} and \eqref{eq:toeplitz-B2} are approximated by two \(\alpha\)-circulant matrices,
\begin{align}\label{eq:circulant-matrices}
    C_{1} = &
    \frac{1}{\Delta t}
    \begin{pmatrix}
        1  &  0 &  0 &  -\alpha \\
        -1 &  1 &  0 &  0 \\
        0  & -1 &  1 &  0 \\
        0  &  0 & -1 &  1 \\
    \end{pmatrix},\\
    C_{2} = &
    \begin{pmatrix}
        \theta  &  0 &  0 &  \alpha(1-\theta) \\
        (1-\theta) &  \theta &  0 &  0 \\
        0  & (1-\theta) &  \theta &  0 \\
        0  &  0 & (1-\theta) &  \theta \\
    \end{pmatrix},
\end{align}
so that the preconditioning operator is
\begin{equation}\label{eq:linear-preconditioner}
    \mathbf{P} = C_{1}\otimes M + C_{2}\otimes K,
\end{equation}
where \(\alpha\in(0,1]\).
The approximation properties of \(\alpha\)-circulant matrices to Toeplitz matrices are well established \citep{gray2006toeplitz}, and are especially favourable for triangular or low bandwidth Toeplitz matrices as is the case here.
The advantange of using (\ref{eq:linear-preconditioner}) is that all \(\alpha\)-circulant matrices are simultaneously diagonalisable with the weighted Fourier transform \citep{benzi2007numerical},
\begin{equation}\label{eq:circulant-decomposition}
    C_{j} = VD_{j}V^{-1},
    \quad
    V = \Gamma^{-1}\mathbb{F}^{*},
    \quad
    D_{j} = \text{diag}(\Gamma\mathbb{F}c_{j}),
\end{equation}
where \(j\in\{1,2\}\), \(\mathbb{F}\) is the discrete Fourier transform matrix, and the eigenvalues in \(D_{j}\) are the weighted Fourier transforms of \(c_{j}\), the first column of \(C_{j}\).
The weighting matrix is \(\Gamma = \text{diag}(\alpha^{(n-1)/N_{t}})\), \(n\in[1,N_{t}]\).
Using the mixed product property of the Kronecker product, \((AB)\otimes(CD)=(A\otimes C)(B\otimes D)\), the eigendecomposition (\ref{eq:circulant-decomposition}) leads to the following factorisation of (\ref{eq:linear-preconditioner}),
\begin{equation}\label{eq:preconditioner-decomposition}
    \mathbf{P} =
    \left(V\otimes I_{x}\right)
    \left(D_{1}\otimes M + D_{2}\otimes K\right)
    \left(V^{-1}\otimes I_{x}\right).
\end{equation}
The inverse \(\mathbf{P}^{-1}\mathbf{x}=\mathbf{b}\) can then be applied efficiently in parallel in three steps, as follows.
\begin{enumerate}\setlength{\itemindent}{+.3in}
    \item[Step 1:]\(\mathbf{z} = (V^{-1}\otimes I_{x})\mathbf{b}\),
    \item[Step 2:] \(\mathbf{y} = \left(D_{1}\otimes M + D_{2}\otimes K\right)^{-1}\mathbf{z} \),
    \item[Step 3:] \(\mathbf{x} = (V\otimes I_{x})\mathbf{y}\).
\end{enumerate}
Steps 1 and 3 correspond to a (weighted) FFT/IFFT in time at each spatial degree of freedom, which is ``embarrassingly'' parallel in space.
Step 2 corresponds to solving a block-diagonal matrix, which is achieved by solving \(N_{t}\) complex valued blocks with structure similar to the implicit problem required for (\ref{eq:theta-method}).
The blocks are independent, so Step 2 is ``embarrassingly'' parallel in time.
Steps 1 and 3 require data aligned in the time direction, and Step 2 requires data aligned in the space direction. 
Switching between these two layouts corresponds to transposing a distributed array (similar transposes are required in parallel multidimensional FFTs). This requires all-to-all collective communication.
The implications of these communications on efficient implementation will be discussed later.

The matrix \(\mathbf{P}\) can then be used as a preconditioner for an iterative solution method for the all-at-once system (\ref{eq:linear-aaos}).
The effectiveness of this preconditioner relies on how well \(\mathbf{P}\) approximates \(\mathbf{A}\).
\cite{mcdonald2018preconditioning} showed that the matrix \(\mathbf{P}^{-1}\mathbf{A}\) has at least \((N_{t}-1)N_{x}\) unit eigenvalues independent of \(N_{t}\) or discretisation and problem parameters.
Although this is a favourable result, the values of the \(N_{x}\) non-unit eigenvalues may depend on the problem parameters so a good convergence rate is not guaranteed.
However, \cite{gander2019convergence} proved that the convergence rate \(\eta\) of Richardson iterations for (\ref{eq:linear-aaos}) preconditioned with (\ref{eq:linear-preconditioner}) scales according to
\begin{equation}\label{eq:linear-convergence-bound}
	\eta \sim \frac{\alpha}{1-\alpha},
\end{equation}
when $\alpha<1/2$.
\cite{caklovic_parallel_2021} proved the same bound for collocation time integration methods, and \cite{wu_stability_2021} proved that (\ref{eq:linear-convergence-bound}) holds for any stable one step time integrator.

This bound on the convergence rate (\ref{eq:linear-convergence-bound}) implies that a very small \(\alpha\) should be used.
However, the roundoff error of the three step algorithm above is \(\mathcal{O}(N_{t}\epsilon\alpha^{-2})\) where \(\epsilon\) is the machine precision \citep{gander2019convergence}, so for very small values of \(\alpha\) the round off errors will become large.
A value of around \(10^{-4}\) is often recommended to provide a good convergence rate without suffering significant round off error \citep{gander2019convergence}.
For improved performance, \cite{caklovic_parallel_2021} devised a method for adapting \(\alpha\) through the iteration to achieve excellent convergence without loss of accuracy.

\begin{figure}
    \centering
    \epswitch{
    \includegraphics[width=\linewidth]{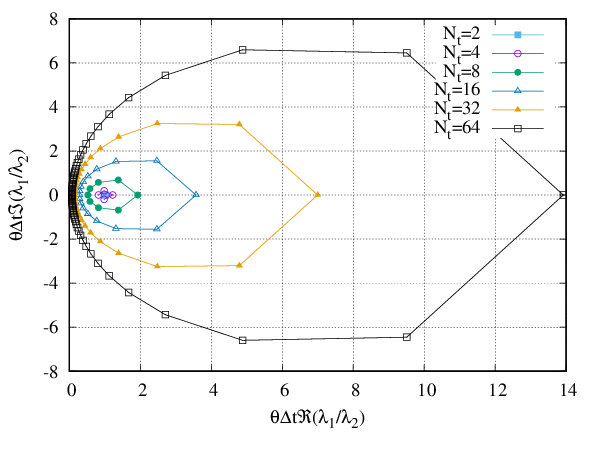}
    }{
    \includegraphics[width=\linewidth]{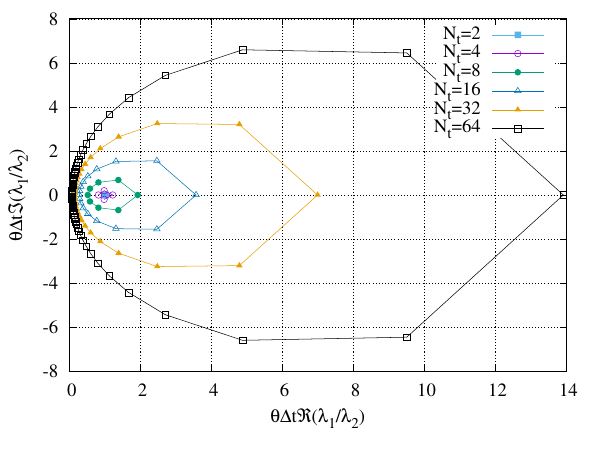}
    }
    \caption{\(\psi_{j}=(\lambda_{1,j}/\lambda_{2,j})/(1/(\Delta t\theta))\) in the complex plane for varying \(N_{t}\) with \(\theta=0.5\) and \(\alpha=10^{-4}\). As \(N_{t}\) increases the \(\psi_{j}\) for low frequencies cluster towards the imaginary axis.}
    \label{fig:circulant-eigenvalues}
\end{figure}

Before moving on to consider nonlinear problems, we will briefly discuss the complex-valued block systems in Step 2 of the algorithm above.
We compare these blocks to the implicit linear system \((M/\Delta t + \theta K)\) in the serial-in-time method (\ref{eq:theta-method}).
For consistent naming convention, we refer to the linear system from the serial-in-time method as the ``real-valued block'' contrasted with the ``complex-valued blocks'' in the circulant preconditioner.
In each case we need to solve blocks of the form
\begin{equation}\label{eq:blocks}
    \left(\frac{\beta_{1}}{\beta_{2}}M + K\right)x = \frac{1}{\beta_{2}}b
\end{equation}
where \((\beta_{1},\beta_{2})\) is \((1/\Delta t, \theta)\) in the serial-in-time method and \((\lambda_{1,j},\lambda_{2,j})\) in the parallel-in-time method where \(\lambda_{1,j}\) and \(\lambda_{2,j}\) are the \(j\)th eigenvalue of \(C_{1}\) and \(C_{2}\) respectively.
We have divided through by \(\beta_{2}\) to make comparison easier.
The ratio of the coefficient on the mass matrix \(M\) in the parallel-in-time method to the coefficient in the serial-in-time method is
\begin{equation}\label{eq:block-coeff-ratio}
    \psi_{j} = \frac{\lambda_{1,j}/\lambda_{2,j}}{1/(\Delta t \theta)}.
\end{equation}
Figure \ref{fig:circulant-eigenvalues} shows \(\psi\) plotted in the complex plane for increasing \(N_{t}\).
When \(N_{t}\) is small, \(\psi_{j}\) is clustered around unity and the blocks in the parallel-in-time method are almost identical to those in the serial-in-time method.
However as \(N_{t}\) increases \(\psi_{j}\) spreads further from unity in all directions.
The magnitude of \(\psi_{j}\) for high frequency modes increases, with a real part \(\geq1\), which is comparable to a small (albeit complex) timestep and is usually a favourable regime for iterative solvers.
On the other hand, the \(\psi_{j}\) for low frequency modes cluster closer and closer to the imaginary axis as \(N_{t}\) increases.
This resembles the case of a very large and mostly imaginary timestep (and hence a large, imaginary, Courant number), which is challenging for many iterative methods.

\subsection{Nonlinear problems}\label{sec:methods:nonlinear}
The method as presented so far is designed for linear problems with constant coefficients, however many problems of interest are nonlinear in nature or have non-constant coefficients.
In this section we will show the all-at-once system for nonlinear problems, and show how the ParaDiag method can be applied to such problems.
We consider PDEs of the form
\begin{equation}\label{eq:nonlinear-pde}
	M\partial_{t}u + f(u, t) = 0,
\end{equation}
where the function \(f(u, t): \mathbb{R}^{N_{x}}\times \mathbb{R} \to \mathbb{R}^{N_{x}}\) may be nonlinear.
The discretisation of (\ref{eq:nonlinear-pde}) with the implicit \(\theta\) method is analogous to (\ref{eq:theta-method}) with \(Ku^{n}\) replaced with \(f(u^{n}, t^{n})\).
The all-at-once system for the nonlinear PDE, analogous to (\ref{eq:linear-aaos}), is
\begin{equation}\label{eq:nonlinear-aaos}
    \left({B}_{1}\otimes M\right)\vec{u}
    + \left({B}_{2}\otimes I_{x}\right)\mathbf{f}(\vec{u, \vec{t}})
    = \vec{\tilde{b}},
\end{equation}
where \(\mathbf{f}(\vec{u}, \vec{t}): \mathbb{R}^{N_{t}N_{x}}\times \mathbb{R}^{N_{t}} \to \mathbb{R}^{N_{t}N_{x}}\) is the concatenation of the function evaluations for the whole timeseries, which we show again for four timesteps:
\begin{equation}\label{eq:aaoform}
    \mathbf{f}(\vec{u}, \vec{t}) =
    \begin{pmatrix}
        f(u^{1}, t^{1}) \\
        f(u^{2}, t^{2}) \\
        f(u^{3}, t^{3}) \\
        f(u^{4}, t^{4})
    \end{pmatrix}\\
\end{equation}
where the vector of time values \(\vec{t}\) is
\begin{equation}\label{eq:aaotime}
    \vec{t} =
    \begin{pmatrix}
        t^{1} \\ t^{2} \\ t^{3} \\ t^{4}
    \end{pmatrix}
\end{equation}
The right hand side \(\vec{\tilde{b}}\) of \eqref{eq:nonlinear-aaos} resembles \eqref{eq:aaorhs} with \(K u^{0}\) replaced with \(f(u^{0})\).
This nonlinear system can be solved with a quasi-Newton method once a suitable Jacobian has been chosen.
The exact Jacobian of (\ref{eq:nonlinear-aaos}) is
\begin{equation}\label{eq:nonlinear-jacobian}
    \mathbf{J}(\vec{u}, \vec{t}) =
    \left(B_{1}\otimes M\right)
    + \left(B_{2}\otimes I_{x}\right)\mathbf{\nabla}_{\vec{u}}\mathbf{f}(\vec{u}, \vec{t}),
\end{equation}
where \(\nabla_{\vec{u}}\vec{f}(\vec{u},\vec{t})\in\mathbb{R}^{N_{t}N_{x}\times N_{t}N_{x}}\) is the derivative of \(\vec{f}\) with respect to \(\vec{u}\), and is a block diagonal matrix with the \(n\)th diagonal block corresponding to the spatial Jacobian of \(f\) with respect to \(u^{n}\), i.e. \(\nabla_{u}f(u^{n},t^{n})\in\mathbb{R}^{N_{x}\times N_{x}}\).
Unlike \(K\) in the constant coefficient linear system (\ref{eq:linear-aaos}), the spatial Jacobian \(\nabla_{u}f\) varies at each timestep either through nonlinearity, time-dependent coefficients, or both.
As such, \eqref{eq:nonlinear-jacobian} cannot be written solely as a sum of Kronecker products like the Jacobian in equation (\ref{eq:linear-aaos}) and we cannot immediately apply the same \(\alpha\)-circulant trick as before to construct a preconditioner.
First, we must choose some constant-in-time value for the spatial Jacobian. \citet{gander2017time} proposed to average the spatial Jacobians over all timesteps \(\overline{\nabla_{u}f}=\sum^{N_{t}}_{n=1}\nabla_{u}f(u^{n},t^{n})/N_{t}\).
In our implementation, we use a constant-in-time reference state \(\hat{u}\) and reference time \(\hat{t}\) and evaluate the spatial Jacobian using these reference values (this choice is discussed further in Sect. \ref{sec:asq:parallelism}).
We also allow the preconditioner to be constructed from a different nonlinear function \(\hat{f}\); this allows for a variety of quasi-Newton methods.
Using \(\nabla_{u}\hat{f}(\hat{u}, \hat{t})\) we can construct an \(\alpha\)-circulant preconditioner according to
\begin{equation}\label{eq:nonlinear-preconditioner}
    \mathbf{P}(\hat{u}, \hat{t}) = C_{1}\otimes M + C_{2}\otimes\nabla_{u}\hat{f}(\hat{u}, \hat{t}).
\end{equation}
The preconditioner (\ref{eq:nonlinear-preconditioner}) can be used with a quasi-Newton-Krylov method for (\ref{eq:nonlinear-aaos}).

The preconditioner (\ref{eq:nonlinear-preconditioner}) has two sources of error: the $\alpha$-circulant block as in the linear case, and now also the difference between \(\nabla_{u}\hat{f}(\hat{u}, \hat{t})\) and \(\nabla_{u}f(u^{n}, t^{n})\) at each step.
By analogy to Picard iterations for the Dahlquist equation with an approximate preconditioner, \cite{caklovic_paradiag_2023} estimates the convergence of the inexact Newton method for nonlinear problems with a collocation method time integrator as 
\begin{equation}\label{eq:nonlinear-convergence-bound}
    \eta \sim \frac{\vartheta\kappa N_{t}\Delta t + \alpha}{1-\alpha},
\end{equation}
where \(N_{t}\Delta t\) is the duration of the time-window, \(\vartheta\) is a constant that depends on the time integrator, and \(\kappa\) is the Lipschitz constant of the ``error'' in the reference spatial Jacobian,
\begin{equation}\label{eq:average-lipschitz}
    \kappa = \text{Lip}\left(f(u) - \left(\nabla_{u}\hat{f}(\hat{u})\right)u \right).
\end{equation}
This result can be understood as stating that the convergence rate deteriorates as the nonlinearity of the problem gets stronger (\(\kappa\) increases) or as the all-at-once system encompasses a longer time window (\(N_{t}\Delta t\) increases).
The estimate implies that if \(\alpha\) is small enough then the predominant error source will be the choice of reference state for any moderate nonlinearity.

We now highlight the flexibility that expressing (\ref{eq:nonlinear-aaos}) as a preconditioned quasi-Newton-Krylov method provides. There are four main aspects, as follows.
\begin{enumerate}
	\item At each Newton iteration, the Jacobian (\ref{eq:nonlinear-jacobian}) can be solved exactly (leading to a “true” Newton method), or inexactly (quasi-Newton).
    In the extreme case of the least exact Jacobian, the Jacobian is simply replaced by the preconditioner (\ref{eq:nonlinear-preconditioner}), as in previous studies \citep{gander2017time,gander2019convergence,caklovic_paradiag_2023, wu2021parallel}.
   	\item The nonlinear function \(\hat{f}\) used in the preconditioner (\ref{eq:nonlinear-preconditioner}) does not need to be identical to the nonlinear function \(f\) in (\ref{eq:nonlinear-aaos}). The Jacobian (\ref{eq:nonlinear-jacobian}) has been written using \(f\), but in general could also be constructed from yet another function (not necessarily equal to either \(f\) or \(\hat{f}\)).
    For example if (\ref{eq:nonlinear-aaos}) contains both linear and nonlinear terms, the Jacobian and preconditioner may be constructed solely from the linear terms, as in \cite{caklovic_paradiag_2023}.
	\item The Jacobian \(\mathbf{J}(\vec{u}, \vec{t})\) (\ref{eq:nonlinear-jacobian}) need not be linearised around the current Newton iterate \(\vec{u}\), but instead could be linearised around some other time-varying state \(\hat{\vec{u}}\), e.g. some reduced order reconstruction of \(\vec{u}\).
	\item Usually the reference state \(\hat{u}\) for the preconditioner \(\mathbf{P}(\hat{u}, \hat{t})\) (\ref{eq:nonlinear-preconditioner}) is chosen to be the time-average state, and is updated at every Newton iteration, however any constant-in-time state may be used, e.g. the initial conditions.
	This may be favourable if \(\hat{u}\) does not change between Newton iterations and a factorisation of \(\nabla_{u}\hat{f}(\hat{u}, \hat{t})\) may be reused across multiple iterations.
\end{enumerate}

\subsection{Performance model}\label{sec:methods:performance-model}

In this section we present a simple performance model for ParaDiag. The purpose of the model is firstly to identify the factors determining the effectiveness of the method, and secondly to help us later demonstrate that we have produced a reasonably performant implementation in asQ.
The model extends those presented in \cite{maday2008parallelization} and \cite{caklovic_parallel_2021} by a more quantitative treatment of the block solve cost.
In \cite{maday2008parallelization} it is assumed that the cost of the block solves in the all-at-once preconditioner is identical to the block solves in the serial-in-time method.
In \cite{caklovic_parallel_2021} the costs are allowed to be different, but the difference is not quantified.
Here we assume that an iterative Krylov method is used to solve both the real- and complex-valued blocks and that, because of the variations in \(\psi_{j}\) \eqref{eq:block-coeff-ratio}, the number of iterations required is different for the blocks in the circulant preconditioner and in the serial-in-time method.
In the numerical examples we will see that accounting for the difference in the number of block iterations is essential for accurately predicting performance.
We assume perfect weak scaling in space for the block iterations i.e.~the time taken per Krylov iteration for the complex-valued blocks in Step 2 of applying the preconditioner is the same as the time taken per Krylov iteration for the real-valued blocks in the serial-in-time method so long as twice the number of processors are used for spatial parallelism.
%\footnote{Finite element computation is almost always memory bound, and the increase in arithmetic intensity from real to complex arithmetic will not change this in the vast majority of cases, therefore the weak scaling assumption is reasonable assuming memory bandwidth per core is also maintained. \jhc{Is this point even necessary or is it reviewer bait? If so, word it better. Assumes bandwidth is linear in core count. Recap discussion with Jack.}}

We start with a simple performance model for the serial-in-time method for \(N_{t}\) timesteps of the system \eqref{eq:nonlinear-pde} with \(N_{x}\) degrees of freedom (DoFs) in space.
Each timestep is solved sequentially using Newton's method, requiring a certain number of (quasi-)Newton iterations per timestep where, at each Newton iteration, the real-valued block \((M/\Delta t + \theta \nabla_{u}f)\) is solved (possibly inexactly) using a Krylov method.
We assume that the block solves constitute the vast majority of the computational work in both the serial- and parallel-in-time methods, and will revisit this assumption later.
The cost of each timestep is then proportional to the number of Krylov iterations \(k_{s}\) per solve of the real-valued block, and to the number of times \(m_{s}\) that the real-valued block must be solved per timestep.
For linear problems \(m_{s}=1\) and for nonlinear problems \(m_{s}\) is the number of Newton iterations.
We assume that each Krylov iteration on the block requires work proportional to \(N_{x}^{q}\), where the exponent \(q\) determines the scalability in space, for example a textbook multigrid method has \(q=1\), and a sparse direct solve of a 2D finite element matrix has \(k_{s}\)=1 and \(q=3/2\) (up to some upper limit on \(N_{x}\)).
The cost of solving \(N_{t}\) timesteps serial-in-time is therefore proportional to
\begin{equation}\label{eq:work-serial}
    W_{s} \sim k_{s}m_{s}N_{x}^{q}N_{t}
\end{equation}
The \(s\) subscript refers to ``serial''(-in-time).
If we parallelise only in space and assume weak scaling, then the number of processors is proportional to the number of spatial degrees of freedom, i.e. \(P_{s} \sim N_{x}\) (this relation also holds once we have reached saturation when strong scaling).
If the spatial parallelism weak scales perfectly with \(N_{x}\), then the time taken to calculate the entire timeseries using the serial-in-time method is
\begin{equation}\label{eq:time-taken-serial}
    T_{s} \approx \frac{W_{s}}{P_{s}} \sim k_{s}m_{s}N_{x}^{q-1}N_{t}.
\end{equation}
from which see that \(T_{s}\) is linear in \(N_{t}\).

The nonlinear all-at-once system (\ref{eq:nonlinear-aaos}) is solved using a (quasi-)Newton method where, at each Newton iteration, a linear system with the Jacobian matrix \(\mathbf{J}\) \eqref{eq:nonlinear-jacobian} is solved (possibly inexactly) using a Krylov method preconditioned with the circulant matrix \(\mathbf{P}\) \eqref{eq:nonlinear-preconditioner} (referred to as the outer Krylov method).
At each outer Krylov iteration, the inverse of the preconditioner \(\mathbf{P}\) is applied by the three step algorithm above, which requires a solve of the complex-valued blocks in Step 2.
The complex-valued blocks are solved (possibly inexactly) with a separate Krylov method (referred to as the inner Krylov method).\footnote{
The outer Krylov method must be flexible if a nonlinear Krylov method is used for the inner iterations.}

The total work is proportional to the number of blocks \(N_{t}\), to the number of inner Krylov iterations \(k_{p}\) per solve of the complex-valued blocks, and to the total number of times \(m_{p}\) that the complex-valued blocks must be solved.
\(m_{p}\) is the number of times the preconditioner \(\mathbf{P}\) is applied throughout the entire Newton solve, i.e. the total number of outer Krylov iterations across all Newton iterations of (\ref{eq:nonlinear-aaos}).
(Note that for linear systems, where \(\mathbf{J}\) reduces to \(\mathbf{A}\), \(m_{p}\neq1\) but is equal to the number of outer Krylov iterations to solve \(\mathbf{A}\)).
If the same Krylov method that was used for the real-valued blocks in the serial-in-time method is also used for the inner Krylov method for the complex-valued blocks, then the cost of each inner Krylov iteration also scales with \(N_{x}^{q}\).
Finally, this gives the following estimate for the total computational work to solve (\ref{eq:nonlinear-aaos}):
\begin{equation}\label{eq:work-parallel}
    W_{p} \sim 2k_{p}m_{p}N_{x}^{q}N_{t}
\end{equation}
where the \(p\) subscript refers to ``parallel''(-in-time), and the factor of 2 accounts for the blocks being complex-valued.

We revisit the assumption that \(W_{p}\) is dominated by the complex-valued block solves.
The main computational components of the ParaDiag method are: the evaluation of the residual of (\ref{eq:nonlinear-aaos}); the action of the all-at-once Jacobian (\ref{eq:nonlinear-jacobian}); solving the complex-valued blocks in the circulant preconditioner; the collective communications for the space-time data transposes; and the (I)FFTs.
In the numerical examples we will see that, in practice, the contributions to the runtime profile of the (I)FFTs and evaluating the residual are negligible, and the contribution of the Jacobian action is only a fraction of that of the block solves, so we do not include these operations in the performance model.
This leaves just the computational work of the block solves, and the time taken for the space-time transpose communications in Steps 1 \& 3.

The parallel-in-time method is parallelised across both time and space, giving
\begin{equation}\label{eq:proc-parallel}
    P_{p} \sim 2N_{x}N_{t}.
\end{equation}
The factor of two ensures that that the number of floating point numbers per processor in the complex-valued block solves is the same as in the real-valued block solves.
The time taken to complete the calculation is then
\begin{equation}\label{eq:time-taken-parallel}
    T_{p} \approx \frac{W_{p}}{P_{p}} + T_{c} \sim k_{p}m_{p}N_{x}^{q-1} + T_{c},
\end{equation}
where \(T_{c}\) is the time taken for the collective communications in the space-time data transposes.
This leads to an estimate \(S\) of the speedup over the serial-in-time (but possibly parallel-in-space) method of
\begin{equation}\label{eq:speedup}
\begin{aligned}
    \frac{T_{s}}{T_{p}} \approx S & = N_{t}\frac{k_{s}m_{s}N_{x}^{q-1}}{k_{p}m_{p}N_{x}^{q-1} + T_{c}} \\
     & = \left(\frac{N_{t}}{\gamma\omega}\right)\frac{1}{1 + T_{c}/T_{b}},
\end{aligned}
\end{equation}
where
\begin{equation}\label{eq:performance-parameters}
    \gamma = \frac{k_{p}}{k_{s}},
    \quad
    \omega = \frac{m_{p}}{m_{s}},
    \quad
    \mbox{ and } \quad
    T_{b} = k_{p}m_{p}N_{x}^{q-1}.
\end{equation}
Here, \(\gamma\) quantifies how much more `difficult' the blocks in the circulant preconditioner are to solve than the blocks in the serial method, \(\omega\) quantifies how much more `difficult' the all-at-once system is to solve than one timestep of the serial method, and \(T_{b}\) is the total time spent solving the blocks of the circulant preconditioner.
If \(T_{c}<<T_{b}\) then the speedup estimate simplifies to \(N_{t}/(\gamma\omega)\) and depends solely on the `algorithmic scaling' i.e. the dependence on \(N_{t}\) of the block and all-at-once system iteration counts.
However the speedup will suffer once \(T_{c}\) becomes non-negligible compared to \(T_{b}\), which we would expect only for larger \(N_{t}\).

From the speedup estimate (\ref{eq:speedup}) we can estimate the parallel efficiency as
\begin{equation}\label{eq:efficiency}
    E = \frac{S}{P_{p}/P_{s}} = \frac{S}{2N_{t}}
    = \left(\frac{1}{2\gamma\omega}\right)\frac{1}{1 + T_{c}/T_{b}}.
\end{equation}
Note that using \(P_{p}/P_{s}\) as the processor count means that (\ref{eq:efficiency}) estimates the efficiency of the time-parallelism independently of the efficiency of the space-parallelism, just as (\ref{eq:speedup}) estimates the speedup over the equivalent serial-in-time method independently of the space-parallelism.

\section{asQ library}\label{sec:asq}
In this section we introduce and describe asQ.
First we state the aims of the library.
Secondly, we take the reader through a simple example of solving the heat equation with asQ.
Once the reader has a picture of the basic usage of asQ, we discuss how the library is structured to meet the stated aims.
Next, we briefly describe the space-time MPI parallelism.
Lastly, we describe the main classes in the library with reference to the mathematical objects from section \ref{sec:methods} that they represent.

asQ is open source under the MIT license and is available at \url{https://github.com/firedrakeproject/asQ}.
It can be installed by either: passing the \verb~--install asQ~ flag when installing or updating Firedrake;\footnote{Instructions on installing Firedrake can be found at \url{https://www.firedrakeproject.org/download.html}} cloning/forking directly from the repository and pip installing in a Firedrake virtual environment; or downloading the Firedrake Docker container.

\subsection{Library overview}\label{sec:asq:overview}

A major difficulty in the development and adoption of parallel-in-time methods is their difficulty of implementation.
Adapting an existing serial-in-time code to be parallel-in-time often requires major overhauls from the top level, and many parallel-in-time codes are written from the ground up as small codes or individual scripts for developing new methods.
The simplest implementation is often to hard code a specific problem and solution method, which then necessitates significant additional effort to test new cases in the future.
Once a promising method is found, it must be tested at scale on a parallel machine to confirm whether it actually achieves real speedups.
However, efficiently parallelising a code can be time consuming in itself, and is a related but distinct skillset to developing a good numerical algorithm.

In light of these observations, we state the three aims of the asQ library for being a productive tool for developing ParaDiag methods.
\begin{enumerate}
    \item It must be straightforward for a user to test out different problems - both equations sets and test cases.
    \item It must be straightforward for a user to try different solution methods. We distinguish between two aspects:
    \begin{enumerate}
        \item The construction of the all-at-once systems, e.g.~specifying the form of the all-at-once Jacobian or the state to linearise the circulant preconditioner around.
        \item The linear/nonlinear solvers used, e.g.~the Newton method used for the all-at-once system, or the preconditioning used for the linear block systems.
    \end{enumerate}
    \item The implementation must be efficient enough that a user can scale up to large-scale parallelism and get a realistic indication of the performance of the ParaDiag method.
\end{enumerate}

Broadly, we attempt to fulfil these aims by building asQ on top of the UFL, Firedrake, and PETSc libraries, and by following the design ethos of these libraries.
The Unified Form Language (UFL) \citep{Alnaes_Logg_Olgaard_Rognes_Wells_2014} is a purely symbolic language for describing finite element forms, with no specification of how those forms are implemented or solved.
UFL expressions can be symbolically differentiated, which allows for the automatic generation of Jacobians and many matrix-free methods.
Firedrake is a Python library for the solution of finite element methods that takes UFL expressions and uses automatic code generation to compile high performance C kernels for the forms.
Firedrake integrates tightly with PETSc via \texttt{petsc4py} \citep{dalcinpazklercosimo2011} for solving the resulting linear and nonlinear systems.
PETSc provides a wide range of composable linear and nonlinear solvers that scale to massive parallelism \citep{Lange_Gorman_Weiland_Mitchell_Southern_2013,May_Sanan_Rupp_Knepley_Smith_2016}, as well as mechanisms for creating new user defined solver components e.g.~preconditioners.
After the example script below we will discuss further how the three aims above are met.

We reinforce that asQ is not intended to be a generic ParaDiag library into which users can port their existing serial-in-time codes, as XBraid is for MGRIT or pySDC is for SDC and collocation methods.
asQ requires a Firedrake installation and for the user to have some familiarity with basic Firedrake usage, and to be considering a modelling approach
that is within the Firedrake paradigm (essentially a discretisation expressible in the domain specific language UFL).
In return for this restriction we gain all the previously mentioned benefits of Firedrake and PETSc for increasing developer productivity and computational performance.
We believe that for ParaDiag in particular it is important to have easy access to a wide range of linear solvers because the complex valued coefficients on the block systems (\ref{eq:blocks}) can impair the performance of iterative methods, and so a range of strategies may have to be explored before finding a sufficiently efficient scheme.
Contrast this with MGRIT or Parareal, for example, where the inner solves are exactly the serial-in-time method so existing solvers are often still optimal.

\subsection{A heat equation example}\label{sec:asq:example}

The main components of asQ will be demonstrated through an example.
%%% TODO: add this back in %%% The full script is provided as supplementary material.
We solve the heat equation over the domain \(\Omega\) with boundary \(\Gamma\),
\begin{align}
    \partial_{t}u - \nabla^{2}u = 0 & \text{ in }\Omega, \label{eq:heat} \\
    n\cdot\nabla u = 0 & \text{ on } \Gamma_{N}, \label{eq:heat-neumann} \\
    u = 0 & \text{ on } \Gamma_{D}, \label{eq:heat-dirichlet}
\end{align}
where \(n\) is the normal to the boundary, \(\Gamma_{N}\) is the section of the boundary with zero Neumann boundary conditions, \(\Gamma_{D}\) is the section of the boundary with zero Dirichlet boundary conditions, and \(\Gamma_{N}\cup\Gamma_{D} = \Gamma\).
Equation (\ref{eq:heat}) will be discretised with a standard continuous Galerkin method in space.
If \(V\) is the space of piecewise linear functions over the mesh, and \(V_{0}\) is the space of functions in \(V\) which are 0 on \(\Gamma_{D}\), the solution \(u\in V_{0}\) is the solution of the weak form,
\begin{equation}
    \int_{\Omega}\left(v\partial_{t}u + \nabla u \cdot\nabla v\right)\text{d}\vec{x} = 0\quad \forall v \in V_{0}.
\end{equation}
The Neumann boundary conditions are enforced weakly by removing the corresponding surface integral, and the Dirichlet conditions are enforced strongly by restricting the solution to \(V_{0}\).

We will set up and solve an all-at-once system for the backwards Euler method for eight timesteps distributed over four spatial communicators.

The first part of the script will import Firedrake and asQ and set up the space-time parallel partition.
asQ distributes the timeseries over multiple processors in time, and each timestep may also be distributed in space.
Firedrake's \texttt{Ensemble} class manages this distribution by setting up a tensor product \(P_{x}\times P_{t}\) of MPI communicators for space and time parallelism (described in more detail later and illustrated in Fig. \ref{fig:ensemble-diagram}).
The time parallelism is determined by the list \texttt{time\_partition}; here we request four processors in time, each holding two timesteps for a total timeseries of eight timesteps.
We refer to all timesteps on a single partition in time as a ``slice'' of the timeseries; here, we have four slices of two timesteps each.
\begin{verbatim}
from firedrake import *
import asQ

time_partition = [2, 2, 2, 2]
ensemble = asQ.create_ensemble(
    time_partition, comm=COMM_WORLD)
\end{verbatim}
We do not need to provide the spatial partition because the total MPI ranks will automatically be distributed evenly in time. If the script is run with 4 MPI ranks, each spatial communicator will have a single rank, but if the script is run with 8 MPI ranks, each spatial communicator will have 2 MPI ranks.
Next we define a mesh and a finite element function space using Firedrake.
The domain is the square \((x, y) \in \Omega = (0, 1)\times(0, 1)\) defined on the local spatial communicator \texttt{ensemble.comm}, and we define a linear continuous Galerkin ("CG") function space \(V\) for a single timestep.
A Firedrake \texttt{Function} is created to hold the initial conditions \(u_{0}=\text{sin}(\pi x/4)\text{cos}(2\pi y)\).
\begin{verbatim}
mesh = SquareMesh(nx=32, ny=32, L=1,
                  comm=ensemble.comm)
x, y = SpatialCoordinate(mesh)

V = FunctionSpace(mesh, "CG", 1)
u0 = Function(V)
u0.interpolate(sin(0.25*pi*x)*cos(2*pi*y))
\end{verbatim}
Now we build the all-at-once system using asQ, starting with an \texttt{AllAtOnceFunction}.
This class represents \eqref{eq:aaofunction}, a timeseries of finite element functions in \(V\) distributed in time over an \texttt{Ensemble}.
\begin{verbatim}
aaofunc = asQ.AllAtOnceFunction(
    ensemble, time_partition, V)
aaofunc.initial_condition.assign(u0)
\end{verbatim}
In an \texttt{AllAtOnceFunction} each spatial communicator holds a slice of one or more timesteps of the timeseries, in this example two timesteps per communicator.

Next we need to define the all-at-once system itself.
To represent the finite element form (\ref{eq:linear-aaos}) or (\ref{eq:nonlinear-aaos}) over the timeseries \texttt{aaofunc} we use an \texttt{AllAtOnceForm}.
Building this requires the timestep \(\Delta t\), the implicit parameter \(\theta\), the boundary conditions, and a way of describing the mass matrix \(M\) and the function \(f(u,t)\) (which for linear equations is just \(Ku\)).
\begin{verbatim}
dt = 0.05
theta = 1

bcs = [DirichletBC(V, 0, sub_domain=1)]

def form_mass(u, v):
    return u*v*dx

def form_function(u, v, t):
    return inner(grad(u), grad(v))*dx

aaoform = asQ.AllAtOnceForm(
    aaofunc, dt, theta, form_mass,
    form_function, bcs=bcs)
\end{verbatim}
The timestep \(\Delta t=0.05\) gives a Courant number of around 13, and \(\theta=1\) gives the implicit Euler time-integration method.
The Python functions \texttt{form\_mass} and \texttt{form\_function} each take in a function \texttt{u} and a test function \texttt{v} in \(V\), and return the UFL form for the mass matrix \(M\) and \(f(u, t)\) respectively. 
asQ uses these two Python functions to generate all of the necessary finite element forms for the all-at-once system (\ref{eq:nonlinear-aaos}), the Jacobian (\ref{eq:nonlinear-jacobian}), and the circulant preconditioner (\ref{eq:nonlinear-preconditioner}), while the user need only define them for the semi-discrete form of a single timestep.
The boundary conditions for a single timestep are passed to the \texttt{AllAtOnceForm} in \texttt{bcs}, and are applied to all timesteps.
For \eqref{eq:heat-dirichlet} we set \(\Gamma_{D}\) to be the left boundary \(x=0\) by passing the \texttt{subdomain=1} argument to the \texttt{DirichletBC}, and for \eqref{eq:heat-neumann} Firedrake defaults to zero Neumann boundary conditions on all other boundaries.
Given \(V\) and a set of Dirichlet boundary conditions, Firedrake will automatically create the restricted space \(V_{0}\).
The initial conditions defined earlier conform to these boundary conditions.

We next need to specify how we solve for the solution of \texttt{aaoform}.
In PETSc, linear and nonlinear solvers are specified using solver parameters dictionaries.
Nonlinear problems are solved using a SNES (Scalable Nonlinear Equations Solver), and linear problems are solved using a KSP (Krylov SPace method) which is PreConditioned by a PC.
Our solver parameters are as follows.
\begin{verbatim}
solver_parameters = {
    'snes_type': 'ksponly',
    'mat_type': 'matfree',
    'ksp_type': 'richardson',
    'ksp_rtol': 1e-12,
    'ksp_monitor': None,
    'ksp_converged_rate': None,
    'pc_type': 'python',
    'pc_python_type': 'asQ.CirculantPC',
    'circulant_block': {'pc_type': 'lu'},
    'circulant_alpha': 1e-4}
\end{verbatim}
Starting with \texttt{solver\_parameters}, our problem is linear so we use \texttt{'snes\_type': 'ksponly'} to perform one Newton iteration and then return the result.
Assembling the entire space-time matrix would be very expensive and is unnecessary for ParaDiag, so asQ implements the Jacobian matrix-free, specified by the \texttt{'mat\_type': 'matfree'} option.
We select the Richardson iterative method, and require a drop in the residual of twelve orders of magnitude.
The \texttt{ksp\_monitor} and \texttt{ksp\_converged\_rate} options tell PETSc to print the residual at each Krylov iteration, and the contraction rate upon convergence respectively.
The Richardson iteration is preconditioned with the corresponding block circulant ParaDiag matrix, which is provided by asQ as a Python type preconditioner \texttt{CirculantPC} with the circulant parameter \(\alpha=10^{-4}\).
Lastly we need to specify how to solve the blocks in the preconditioner after the diagonalisation.
The composability of PETSc solvers means that to specify the block solver we simply provide another parameters dictionary with the \texttt{'circulant\_block'} key.
The block solver here is just a direct LU factorisation but could be any other Firedrake or PETSc solver configuration.

The tight integration of Firedrake and asQ with PETSc means that a wide range of solution strategies are available through the options dictionary.
For example, we can change the Krylov method for the all-at-once system simply by changing the \texttt{'ksp\_type'} option.
Rather than a direct method, an iterative method could be used for the block solves by changing the options in the \texttt{circulant\_block} dictionary.
Options can also be passed from the command line, in which case zero code changes are required to experiment with different solution methods.

The last all-at-once object we need to create is an \texttt{AllAtOnceSolver}, which solves \texttt{aaoform} for \texttt{aaofunc} using the solver parameters specified above.
\begin{verbatim}
aaosolver = asQ.AllAtOnceSolver(
    aaoform, aaofunc, solver_parameters)
\end{verbatim}

If we want to solve for say a total of 48 timesteps then, for an all-at-once system of size \(N_{t}=8\), we need to solve six windows of 8 timesteps each, where the final timestep of each window is used as the initial condition for the next window.
\begin{verbatim}
aaofunc.assign(u0)
for i in range(6):
    aaosolver.solve()
    aaofunc.bcast_field(
        -1, aaofunc.initial_condition)
    aaofunc.assign(
        aaofunc.initial_condition)
\end{verbatim}
After each window solve, we use the final timestep as the initial condition for the next window.
To achieve this, \texttt{aaofunc.bcast\_field(i, uout)} wraps \texttt{MPI\_Bcast} to broadcast timestep \texttt{i} across the temporal communicators, so that all spatial communicators now hold a copy of timestep \texttt{i} in \texttt{uout} (in the code above \texttt{i} is the Pythonic -1 for the last element, and \texttt{uout} is the \texttt{initial\_condition}).
After broadcasting the new initial conditions, every timestep in \texttt{aaofunc} is then assigned the value of the new initial condition as the initial guess of the next solve. We have been able to set up and solve a problem parallel in time in under 50 lines of code.

We now return to the three requirements stated in the library overview and discuss how asQ attempts to meet each one.

\begin{enumerate}
    \item \textit{Straightforward to test different problems}: This requirement is met by only requiring the user to provide UFL expressions for the mass matrix \(M\) and the function \(f(u, t)\); a Firedrake \texttt{FunctionSpace} for a single timestep; and a Firedrake \texttt{Function} for the initial conditions.
    These are all components that a user would already have, or would need anyway, to implement the corresponding serial-in-time method.
    From these components asQ can generate the UFL for the different all-at-once system components which it then hands to Firedrake to evaluate numerically.
    Changing to a different problem simply requires changing one or more of the UFL expressions, the function space, or the initial conditions.
    \item \textit{Straightforward to test different solution methods}:
    \begin{enumerate}
        \item Requirement 2a is met by allowing users to optionally pass additional UFL expressions (\texttt{form\_function}) for constructing the different all-at-once system components (see the \texttt{AllAtOnceSolver}, \texttt{CirculantPC} and \texttt{AuxiliaryBlockPC} descriptions below).
        \item Requirement 2b is met through the use of PETSc's solver options interface.
        Changing between many different methods is as simple as changing some options strings. For more advanced users, novel methods can be written as bespoke petsc4py Python PCs.
        The use of UFL means that symbolic information is retained all the way down to the block systems, so solution methods that rely on certain structure can be applied without issue, for example Schur factorisations or additive Schwarz methods defined on topological entities.
    \end{enumerate}
    \item \textit{Efficient parallel implementation}: This requirement involves both spatial and temporal parallelism.
    The spatial parallelism is entirely provided by Firedrake and PETSc.
    The temporal parallelism is implemented using a mixture of PETSc objects defined over the global communicator, and mpi4py calls via the Firedrake \texttt{Ensemble} or via the mpi4py-fftw library \citep{dalcin2021mpi4py}.
    The temporal parallelism is discussed in more depth in Sect. \ref{sec:asq:parallelism} and is profiled in the examples in Sect. \ref{sec:examples}.
\end{enumerate}

\subsection{Space-time parallelism}\label{sec:asq:parallelism}

As stated previously, time parallelism is only used once space parallelism is saturated due to the typically lower parallel efficiency.
This means that any practical implementation of parallel-in-time methods must be space-time parallel.
In terms of evaluating performance, full space-time parallelism is especially important for ParaDiag methods due to the need for all-to-all communication patterns which can be significantly affected by network congestion.

The three steps in applying the block circulant preconditioner require two distinct data access patterns.
The (I)FFTs in Steps 1 and 3 require values at a particular spatial degree of freedom (DoF) at all \(N_{t}\) timesteps/frequencies, whereas the block solves in Step 2 require values for all spatial DoFs at a particular frequency.
The data layout in asQ places a single slice of the timeseries on each spatial communicator. The slice length is assumed to be small and is usually just one i.e. a single timestep per spatial communicator, to maximise time parallelism during Step 2.
This layout minimises the number of ranks per spatial communicator, hence reducing the overhead of spatial halo swaps during the block solve.
However, parallel FFTs are not efficient with so few values per processor.
Instead of using a parallel FFT, the space-time data is transposed so that each rank holds the entire timeseries for a smaller number of spatial DoFs.
Each rank can then apply a serial FFT at each spatial DoF.
After the transform the data is transposed back to its original layout.
These transposes require all-to-all communication rounds, which are carried out over each partition of the mesh separately so we have \(P_{x}\sim\mathcal{O}(N_{x})\) all-to-all communications involving \(P_{t}\sim\mathcal{O}(N_{t})\) processors each, instead of a single all-to-all communication involving all \(P_{t}P_{x}\) processors.
asQ currently uses the \texttt{mpi4py-fftw} library for these communications, which implements the transposes using \texttt{MPI\_Alltoallw} and MPI derived datatypes.

Collective communication patterns, particularly \textit{Alltoall} and its variants, do not scale well for large core counts compared to the point-to-point communications typical required for spatial parallelism.
There are two common implementations of Alltoall in MPI: pairwise and Bruck \citep{netterville2022visual}.
The pairwise algorithm minimises the total communication volume but has communication complexity \(\mathcal{O}(P)\).
The Bruck algorithm has communication complexity \(\mathcal{O}(\log P)\) at the expense of a higher total communication volume, so is only more efficient than the pairwise algorithm for latency bound communications (i.e., small message sizes).
MPI implementations usually select which algorithm to use at run-time based on a variety of factors, primarily the number of processors in the communicator and the message size.
%The typical message size above which the pairwise algorithm is used is usually around **bytes.
%As stated above, typically time-parallelism will only be used once spatial parallelism is saturated, so we expect the DoFs/core to typically be in the range \(10^{4}-10^{5}\).
%In the numerical examples later we will see that for nonlinear problems the highest speedup usually occurs for \(N_{t}<128\).
%Assuming double precision this means that the message size between each processor will be at least several kilobytes, so the pairwise algorithm will be selected.
The typical message size for ParaDiag is above the usual thresholds for the Bruck algorithm to be used, therefore we expect the communication time \(T_{c}\) in our performance model above to scale with \(\mathcal{O}(P_{t})\) and become a limiting factor on the parallel scaling for large \(N_{t}\).
It is worth noting that, although to a first approximation \(T_{c}\) should depend on \(P_{t}\) but not \(P_{x}\) for constant DoFs/core, in practice network congestion may cause \(T_{c}\) to increase with \(P_{x}\) even if \(P_{t}\) is fixed.

\begin{figure}
    \centering
    \includegraphics[width=8cm]{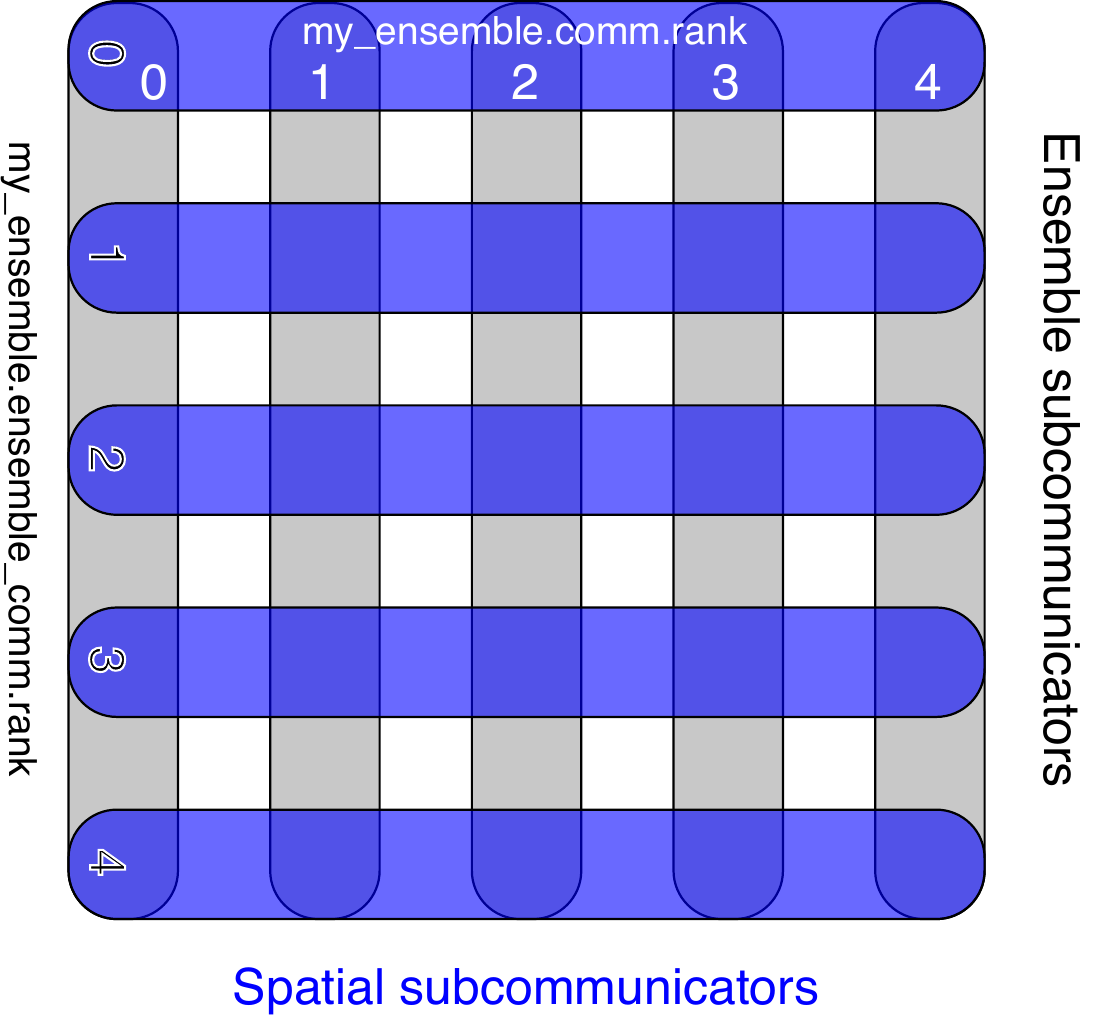}
    \caption{Space-time parallelism using Firedrake's \texttt{Ensemble} class. \texttt{ensemble.comm.size} is \(P_{x} = 5\) and \texttt{ensemble.ensemble\_comm.size} is \(P_{t} = 5\). The dark blue horizontal lines each represent a spatial communicator \texttt{Ensemble.comm} over which a mesh is partitioned.
    Every spatial communicator has the same number of ranks and the same mesh partitioning.
    The grey vertical lines represent communicators in time, \texttt{Ensemble.ensemble\_comm}.
    These communicators are responsible for connecting the ranks on each spatial communicator with the same section of the mesh partition.}
    \label{fig:ensemble-diagram}
\end{figure}

The space-time parallelism in asQ is implemented using Firedrake's \texttt{Ensemble} class, which splits a global communicator (usually \texttt{COMM\_WORLD}) into a Cartesian product of spatial and temporal communicators.
The resulting layout of communicators in the \texttt{Ensemble} is shown in Fig. \ref{fig:ensemble-diagram}.
Firedrake objects for each timestep are defined on each spatial communicator, and PETSc objects and asQ objects for the all-at-once system are defined over the global space-time communicator.
The spatial parallelism required is already fully abstracted away by Firedrake, so asQ need only implement the temporal parallelism on top of this.
The \texttt{Ensemble} provides wrappers around mpi4py calls for sending Firedrake \texttt{Function} objects between spatial communicators over the temporal communicators, and mpi4py-fftw is used for space-time transposes separately across each temporal communicator.

As stated earlier, in asQ we construct the circulant preconditioner from the spatial Jacobian of a constant reference state (usually the time averaged state), rather than the time average of the spatial Jacobians.
The time average of the spatial Jacobian is optimal in the L2 norm, but using the reference state has two advantages.
First, the communication volume of averaging the state is smaller than averaging the Jacobians. 
For scalar PDEs or low order methods the difference is minor, but for systems of PDEs or higher order methods the difference can be substantial.
Second, explicitly assembling the Jacobian prevents the use of matrix-free methods provided by Firedrake and PETSc for the block solves.
Matrix free methods have lower memory requirements than assembled methods, which is important for finite element methods because they are usually memory bound.
The matrix free implementation in Firedrake \citep{Kirby_Mitchell_2018} also enables the use of block preconditioning techniques such as Schur factorisation by retaining symbolic information that is lost when assembling a matrix numerically.
We also note that for quadratic nonlinearities (such as in the shallow water equations in primal form) the Jacobian of the average and the average of the Jacobian are identical.
However, for cubic or higher nonlinearities (such as the compressible Euler equations in conservative form) the two forms are not identical.

\subsection{asQ components}\label{sec:asq:components}

Now we have seen an example of the basic usage of asQ, we give a more detailed description of each component of the library, with reference to which mathematical objects described in Sect. \ref{sec:methods} they represent.
First the components of the all-at-once system are described, then a variety of Python type PETSc preconditioners for both the all-at-once system and the blocks, then the submodule responsible for the complex-valued block solves in the circulant preconditioner.
There are a number of other components implemented in asQ, but these are not required for the numerical examples later so are omited here (see Appendix \ref{app:other-asq} for a description of some of these components).

\subsubsection*{\texttt{AllAtOnceFunction}}

The \texttt{AllAtOnceFunction} is a timeseries of finite element functions distributed over an \texttt{Ensemble} representing the vector \(\vec{u}\) in (\ref{eq:aaofunction}) plus an initial condition.
Each ensemble member holds a slice of the timeseries of one or more timesteps, and  timestep \texttt{i} can be accessed on its local spatial communicator as a Firedrake \texttt{Function} using \texttt{aaofunc[i]}.
Evaluation of the \(\theta\)-method residual \eqref{eq:theta-method} at a timestep \(u^{n+1}\) requires the solution at the previous timestep \(u^{n}\), so the \texttt{AllAtOnceFunction} holds time-halos on each spatial communicator for the last timestep on the previous slice.
Some other useful operations are also defined, such as calculating the time-average, linear vector space operations (axpy etc.), and broadcasting particular timesteps to all slices (\texttt{bcast\_field}).
Internally, a PETSc \texttt{Vec} is created for the entire space-time solution to interface with PETSc solvers.

\subsubsection*{\texttt{AllAtOnceCofunction}}
In addition to functions in the primal space \(V\), it is also useful to have a representation of cofunctions in the dual function space \(V^{*}\).
asQ provides this dual object as an \texttt{AllAtOnceCofunction} which represents a timeseries of the dual space \(V^{*}\).
Cofunctions are used when assembling residuals of (\ref{eq:linear-aaos}) or (\ref{eq:nonlinear-aaos}) over the timeseries, and providing the right hand side (the constant part) of linear or nonlinear systems.
Both \texttt{AllAtOnceFunction} and \texttt{AllAtOnceCofunction} have a \texttt{riesz\_representation} method for converting between the two spaces.

\subsubsection*{\texttt{AllAtOnceForm}}
The \texttt{AllAtOnceForm} holds the finite element form for the all-at-once system (\ref{eq:nonlinear-aaos}), and is primarily responsible for assembling the residual of this form to be used by PETSc.
It holds the \(\Delta t\), \(\theta\), \texttt{form\_mass}, \texttt{form\_function}, and boundary conditions which are used for constructing the Jacobian and preconditioners.
The contribution of the initial conditions to the right hand side \(\vec{\tilde{b}}\) is automatically included in the system.
After updating the time-halos in the \texttt{AllAtOnceFunction}, assembly of the residual can be calculated in parallel across all slices.
The \texttt{AllAtOnceForm} will also accept an \(\alpha\) parameter to introduce the circulant approximation at the PDE level, such as required for the waveform relaxation iterative method of \cite{gander2019convergence} which uses the modified initial condition \(u^{k}(0) = u_{0} + \alpha(u^{k}(T)-u^{k-1}(T))\). 

\subsubsection*{\texttt{AllAtOnceJacobian}}
The \texttt{AllAtOnceJacobian} represents the action of the Jacobian (\ref{eq:nonlinear-jacobian}) of an \texttt{AllAtOnceForm} on an all-at-once function.
This class is rarely created explicitly by the user, but is instead created automatically by the \texttt{AllAtOnceSolver} or \texttt{LinearSolver} classes described below, and used to create a PETSc \texttt{Mat} (matrix) for the matrix-vector multiplications required for Krylov methods.
As for assembly of the \texttt{AllAtOnceForm}, calculating the action of the \texttt{AllAtOnceJacobian} is parallel across all slices after the time-halos are updated.
The automatic differentiation provided by UFL and Firedrake means that the action of the Jacobian is computed matrix free, removing any need to ever explicitly construct the entire space-time matrix.
For a quasi-Newton method the Jacobian need not be linearised around the current solution iterate.
The \texttt{AllAtOnceJacobian} has a single PETSc option for specifying alternate states to linearise around, including the time-average, the initial conditions, or a completely user defined state.

\subsubsection*{\texttt{AllAtOnceSolver}}
The \texttt{AllAtOnceSolver} is responsible for setting up and managing the PETSc SNES for the all-at-once system defined by an \texttt{AllAtOnceForm} using a given dictionary of solver options.
The nonlinear residual is evaluated using the \texttt{AllAtOnceForm}, and an \texttt{AllAtOnceJacobian} is automatically created.
By default this Jacobian is constructed from the \texttt{AllAtOnceForm} being solved, but \texttt{AllAtOnceSolver} also accepts an optional argument for a different \texttt{AllAtOnceForm} to construct the Jacobian from, for example one where some of the terms have been dropped or simplified.
Just like Firedrake solver objects, \texttt{AllAtOnceSolver} accepts an \texttt{appctx} dictionary containing any objects required by the solver which cannot go into an options dictionary (i.e. not strings or numbers).
For example for \texttt{CirculantPC} this could include an alternative \texttt{form\_function} for \(\hat{f}\) or a Firedrake \texttt{Function} for \(\hat{u}\) in (\ref{eq:nonlinear-preconditioner}).

Callbacks can also be provided to the \texttt{AllAtOnceSolver} to be called before and after assembling the residual and the Jacobian, similar to the callbacks in Firedrake's \texttt{NonlinearVariationalSolver}.
There are no restrictions on what these callbacks are allowed to do, but potential uses include updating solution dependent parameters before each residual evaluation, and manually setting the state around which to linearise the \texttt{AllAtOnceJacobian}.

\subsubsection*{\texttt{Paradiag}}
In most cases, users will want to follow a very similar workflow to that in the example above: setting up a Firedrake mesh and function space, building the various components of the all-at-once system, then solving one or more windows of the timeseries.
To simplify this, asQ provides a convenience class \texttt{Paradiag} to handle this case.
A \texttt{Paradiag} object is created with the following arguments,
\begin{verbatim}
paradiag = Paradiag(
    ensemble=ensemble,
    form_mass=form_mass,
    form_function=form_function,
    ics=u0, dt=dt, theta=theta,
    time_partition=time_partition,
    solver_parameters=solver_parameters)
\end{verbatim}
and creates all the necessary all-at-once objects from these arguments.
The windowing loop in the last snippet of the example above is also automated, and can be replaced with
\begin{verbatim}
paradiag.solve(nwindows=6)
\end{verbatim}
The solution is then available in \texttt{paradiag.aaofunc}, and the other all-at-once objects are similarly available.
Callbacks can be provided to be called before and after each window solve e.g. for writing output or collecting diagnostics.

\subsubsection*{\texttt{SerialMiniapp}}                                                                                                                                                       
The \texttt{SerialMiniApp} is a small convenience class for setting up a serial-in-time solver for the implicit \(\theta\)-method, and takes most of the same arguments as the \texttt{Paradiag} class.
\begin{verbatim}
serial_solver = SerialMiniApp(
   dt, theta, u0,
   form_mass, form_function,
   solver_parameters)
\end{verbatim}
The consistency in the interface means that, once the user is satisfied with the performance of the serial-in-time method, setting up the parallel-in-time method is straightforward, with the only new requirements being the solver options and the description of the time-parallelism (ensemble and time partition).
Importantly, this also makes ensuring fair performance comparisons easier.
In Sect. \ref{sec:examples}, all serial-in-time results are obtained using the \texttt{SerialMiniApp} class.

\subsubsection*{\texttt{CirculantPC}}
The block \(\alpha\)-circulant matrix (\ref{eq:nonlinear-preconditioner}) is implemented as a Python type petsc4py preconditioner, which allows it to be applied matrix free.
By default the preconditioner is constructed using the same \(\Delta t\), \(\theta\), \texttt{form\_mass}, and \texttt{form\_function} as the \texttt{AllAtOnceJacobian}, and linearising the blocks around the time average \(\hat{u}\).
However, alternatives for all of these can be set via the PETSc options and the \texttt{appctx}.
For example, an alternate \(\hat{f}\) can be passed in through the appctx, and \(\hat{u}\) could be chosen as the initial conditions via the PETSc options.
Additional solver options can also be set for the blocks, as shown in the example above.
The implementation of the complex-valued linear system used for the blocks can also be selected via the PETSc options: more detail on this below.

\subsubsection*{\texttt{AuxiliaryRealBlockPC} \& \texttt{AuxiliaryComplexBlockPC}}
These two classes are Python preconditioners for a single real- or complex-valued block respectively, constructed from an ``auxiliary operator'' i.e. from a different finite element form than the one used to construct the block.
These are implemented by subclassing Firedrake's \texttt{AuxiliaryOperatorPC}.
Examples of when preconditioning the blocks with a different operator may be of interest include preconditioning the nonlinear shallow water equations with the linear shallow water equations, or preconditioning a variable diffusion heat equation with the constant diffusion heat equation.
This preconditioner can be supplied with any combination of alternative values for \texttt{form\_mass}, \texttt{form\_function}, boundary conditions, \((\Delta t,\theta)\) or \((\lambda_{1},\lambda_{2})\) as appropriate, and \((\hat{u},\hat{t})\), as well as PETSc options for how to solve the matrix resulting from the auxiliary operator.

\subsubsection*{\texttt{complex\_proxy}}
When PETSc is compiled, a scalar type must be selected (e.g. double precision float, single precision complex etc) and then only this scalar type is available.
asQ uses PETSc compiled with a real-valued scalar because complex scalars would double the required memory access when carrying out real-valued parts of the computation, and not all of Firedrake's preconditioning methods are available in complex mode yet (notably geometric multigrid).
However, this means that the complex valued linear systems in the blocks of the \texttt{CirculantPC} must be manually constructed.
The \texttt{complex\_proxy} module does this by writing  the linear system \(\lambda Ax = b\) with complex number \(\lambda = \lambda_{r} + i\lambda_{i}\), complex vectors \(x = x_{r} + ix_{i}\) and \(b = b_{r} + ib_{i}\), and real matrix \(A\) as:
\begin{equation}\label{eq:complex-linear-system}
    \begin{pmatrix}
        \lambda_{r}A & -\lambda_{i}A \\
        \lambda_{i}A & \phantom{-}\lambda_{r}A
    \end{pmatrix}
    \begin{pmatrix}
        x_{r} \\ x_{i}
    \end{pmatrix}
    =
    \begin{pmatrix}
        b_{r} \\ b_{i}
    \end{pmatrix}
\end{equation}
Two versions of \texttt{complex\_proxy} are available, both having the same API but different implementations.
One implements the complex function space as a two component Firedrake \texttt{VectorFunctionSpace}, and one which implements the complex function space as a two component Firedrake \texttt{MixedFunctionSpace}.
The \texttt{VectorFunctionSpace} implementation is the default because for many block preconditioning methods it acts more like a true complex-valued implementation.
The \texttt{complex\_proxy} submodule API will not be detailed here, but example scripts are available in the asQ repository demonstrating its use for testing block preconditioning strategies without having to setup an entire all-at-once system.

\section{Numerical examples}\label{sec:examples}

We now present a set of numerical examples of increasing complexity in order to demonstrate two points, firstly that asQ is a correct and efficient implementation of the ParaDiag method, and secondly that the ParaDiag method is capable of delivering speedups on relevant linear and nonlinear test cases from the literature.

The correctness of asQ will be demonstrated by verifying the convergence rate (\ref{eq:linear-convergence-bound}) from Sect. \ref{sec:methods:linear}, and the efficiency of the implementation will be demonstrated by comparing actual wallclock speedups to the performance model of Sect. \ref{sec:methods:performance-model}.

In the current manuscript we are interested only in speedup over the equivalent serial-in-time method, i.e.~       can ParaDiag accelerate the calculation of the solution to the implicit-\(\theta\) method.
We are not considering here the question of whether ParaDiag can beat the best of some set of serial-in-time methods (e.g.~on error vs wallclock time).
This question is the topic of a later publication.

In the performance model (\ref{eq:speedup}) it is assumed that once spatial parallelism is saturated, the number of processors is proportional to the number of timesteps.
We have not found a situation in which there is a speedup advantage to having more timesteps than processors by having multiple timesteps per time slice.
Therefore in all examples we keep the DoFs/processor constant when varying \(N_{t}\) or \(N_{x}\) which, when combined with the framing of weak and strong scaling in terms of resolution and wallclock time, leads to the following interpretations of strong and weak scaling.
For strong scaling in time, the spatial (\(\Delta x\)) and temporal (\(\Delta t\)) resolution is fixed, and we attempt to decrease the wallclock time per timestep by increasing \(N_{t}\).
Although the number of DoFs being computed at any one time is increasing with DoFs/processor fixed (as it would be for traditional weak scaling in space), we call this strong scaling because the resolution is fixed. For weak scaling in time and space, the spatial and temporal resolutions are refined simultaneously, and we attempt to keep the wallclock time fixed by increasing \(N_{t}\) proportionally to keep \(T=N_{t}\Delta t\) constant.
Note that for the nonlinear test cases we split the full timeseries of \(N_{T}\) timesteps into \(N_{w}\) windows of \(N_{t}\) timesteps each, where \(N_{w}N_{t}=N_{T}\), and solve the all-at-once system for each window sequentially.
In these cases \(N_{t}\) always refers to the number of timesteps in the window, not the total number of timesteps in the timeseries.

All results presented here were obtained on the ARCHER2 HPE Cray supercomputer at the EPCC, using a Singularity container with asQ, Firedrake, and PETSc installed.
ARCHER2 consists of 5860 nodes, each with 2 AMD EPYC 7742 CPUs with 64 cores per CPU for a total of 128 cores per node.
The nodes are connected with an HPE Slingshot network.
The EPYC 7742 CPUs have a deep memory hierarchy, with the lowest level of shared memory being a 16MB L3 cache shared between 4 cores.
During preliminary testing we found that, due to the memory bound nature of finite element computation, the best performance was obtained by underfilling each node by allocating only two cores per L3 cache.
This strategy is used in all examples unless otherwise stated, giving a maximum of 64 cores per node.
For all cases we have only one timestep per \texttt{Ensemble} member to maximise the available time-parallelism.
Unless otherwise stated, we use twice as many cores per timestep in the parallel-in-time method as in the serial-in-time method to account for the complex-valued nature of the blocks in the circulant preconditioner.

The software used to generate the results in this section is available via Zenodo: \citet{hope_collins_2024_asq_scripts} for the Python scripts; \citet{hope_collins_2024_asq_container} for the Singularity image containing a build of PETSc, Firedrake, and asQ; and \citet{firedrake_zenodo_2024_asq} for the versioned Firedrake and PETSc libararies, and their dependencies.
The Python scripts may be run with the Singularity container without requiring a separate Firedrake installation.

\subsection{Advection equation}\label{sec:examples:advection}

The first example we show is the linear scalar advection equation, the prototypical hyperbolic PDE
\begin{equation}\label{eq:scalar-advection}
    \partial_{t}q + \vec{u}\cdot\nabla q = 0,
\end{equation}
where \(q\) is the advected quantity and \(\vec{u}\) is the advecting velocity.
This very simple equation demonstrates the fundamental wave propagation behaviour of hyperbolic problems, and has proved challenging for many parallel-in-time schemes.
This equation has a single characteristic wave travelling at speed \(|\vec{u}|\).
(\ref{eq:scalar-advection}) is discretised with a linear discontinuous Galerkin method in space with an upwind flux at the element interfaces, and the trapezium rule in time.
%%% TODO: add this back in %%% The weak form used can be found in \gmd{the supplementary materials}{Appendix \ref{app:scalar-advection}}.
This discretisation is second order in the $L_2$ norm in both space and time and is unconditionally stable.
The test case is a Gaussian bump advected in a periodic unit square, with a constant velocity \(\vec{u}\) at an angle \(30\degree\) to the mesh lines.

\begin{table*}[t]
  \caption{Iteration counts \(m_{p}\), contraction rates \(\eta\), and measured speedup \(S\) for the preconditioned Richardson iterations with varying \(\alpha\) and window lengths \(N_{t}\).
    The expected contraction rate is \(\eta_{e} = \alpha/(1-\alpha)\).}
  \label{tab:advection-contraction-rates}
  \centering
  \begin{tabular}{rccccccccccccccc}
    \toprule
  & \multicolumn{3}{c}{\(\alpha=10^{-1}\)}
            & \multicolumn{3}{c}{\(\alpha=10^{-2}\)}
             & \multicolumn{3}{c}{\(\alpha=10^{-3}\)}
              & \multicolumn{3}{c}{\(\alpha=10^{-4}\)}
    & \multicolumn{3}{c}{\(\alpha=10^{-6}\)} \\
    \cmidrule(rl){2-4} \cmidrule(rl){5-7} \cmidrule(rl){8-10}
    \cmidrule(rl){11-13} \cmidrule(rl){14-16}
 \(N_{t}\)   & {\(m_{p}\)} & {\(\eta/\eta_{e}\)} & \(S\)
     & {\(m_{p}\)} & {\(\eta/\eta_{e}\)} & \(S\)
      & {\(m_{p}\)} & {\(\eta/\eta_{e}\)} & \(S\)
       & {\(m_{p}\)} & {\(\eta/\eta_{e}\)} & \(S\)
        & {\(m_{p}\)} & {\(\eta/\eta_{e}\)} & \(S\) \\
    \midrule
    2 & {12} & {1.007} & 0.047 & {6} & {1.026} & 0.088 & {4} & {1.048} & 0.123 & {3} & {1.073} & 0.152 & {2} & {1.125} & 0.208 \\
    4 & {12} & {1.010} & 0.090 & {6} & {1.040} & 0.168 & {4} & {1.079} & 0.237 & {3} & {1.120} & 0.296 & {2} & {1.208} & 0.401 \\
    8 & {12} & {1.008} & 0.179 & {6} & {1.049} & 0.336 & {4} & {1.096} & 0.472 & {3} & {1.148} & 0.588 & {2} & {1.258} & 0.796 \\
   16 & {12} & {0.995} & 0.356 & {6} & {1.052} & 0.671 & {4} & {1.104} & 0.938 & {3} & {1.161} & 1.17  & {2} & {1.282} & 1.58  \\
   32 & {12} & {0.961} & 0.697 & {6} & {1.048} & 1.32  & {4} & {1.105} & 1.86  & {3} & {1.162} & 2.29  & {2} & {1.284} & 3.14  \\
   64 & {12} & {0.905} & 1.37  & {6} & {1.034} & 2.62  & {4} & {1.090} & 3.68  & {3} & {1.140} & 4.55  & {2} & {1.243} & 6.25  \\
  128 & {11} & {0.869} & 3.00  & {6} & {1.007} & 5.13  & {4} & {1.040} & 7.31  & {3} & {1.063} & 8.97  & {2} & {1.108} & 12.1  \\
  256 & {11} & {0.890} & 5.96  & {6} & {0.994} & 10.1  & {4} & {1.013} & 14.5  & {3} & {1.022} & 17.8  & {2} & {1.038} & 24.3  \\
  512 & {11} & {0.884} & 11.2  & {6} & {0.997} & 19.7  & {4} & {1.021} & 27.5  & {3} & {1.034} & 32.6  & {2} & {1.058} & 46.7  \\
 1024 & {11} & {0.887} & 21.9  & {6} & {1.000} & 36.3  & {4} & {1.029} & 51.9  & {3} & {1.045} & 64.5  & {2} & {1.077} & 86.2  \\
 2048 & {11} & {0.884} & 40.3  & {6} & {0.988} & 68.4  & {4} & {1.008} & 98.7  & {3} & {1.015} & 126   & {2} & {1.025} & 162   \\
 4096 & {12} & {0.921} & 62.0  & {6} & {0.969} & 119   & {4} & {0.957} & 163   & {3} & {0.938} & 201   & {2} & {0.899} & 267   \\
 8192 & {11} & {0.863} & 103   & {6} & {1.004} & 179   & {4} & {1.037} & 249   & {3} & {1.058} & 298   & {2} & {1.102} & 420   \\
16384 & {11} & {0.896} & 136   & {6} & {0.991} & 241   & {4} & {1.007} & 325   & {3} & {1.013} & 419   & {2} & {1.020} & 528   \\
\bottomrule
  \end{tabular}
\end{table*}

We use this case to demonstrate the linear convergence rate (\ref{eq:linear-convergence-bound}) to be independent of \(N_{t}\).
Measured speedups will be shown versus the performance model predictions (\ref{eq:speedup}) to verify the efficiency of the implementation in asQ.
All results will be presented as strong scaling in time.
A quadrilateral mesh with \(128^{2}\) elements is used resulting in \(\approx\)65kDoFs, which is small enough to fit on a single core for the serial-in-time method.
The Courant number is \(\sigma = |\vec{u}|\Delta t/\Delta x = 0.8\), with window sizes in powers of two from \(N_{t}=2\) up to \(N_{t}=2^{14}=16384\).
The strong scaling wallclock measurements will be reported as time taken per timestep calculated.

The linear system for the serial-in-time method is solved directly using LU factorisation with the MUMPS package \citep{MUMPS:1,MUMPS:2}.
The factorisation is precalculated and reused across all timesteps.
The all-at-once system is solved to a tolerance of \(10^{-11}\) using preconditioned Richardson iterations.
The direct solver is also used for the blocks in the circulant preconditioner, meaning that \(\gamma=1\) in the performance model (\ref{eq:performance-parameters}).
Direct solvers do not scale well at the low DoFs/core count of this case, so the number of cores is kept equal to the number of timesteps in the ParaDiag method instead of twice the number of timesteps.
The speedup predictions from the performance model are halved accordingly.
It was found that the best results were obtained by allocating only a single core per L3 cache up to a maximum of 32 cores per node. This essentially ensures that the memory resource per core remains the same between the serial-in-time method and the parallel method.

The number of Richardson iterations per window \(m_{p}\) and the measured convergence rates are shown in Table \ref{tab:advection-contraction-rates} for varying \(\alpha\) from \(10^{-1}\) to \(10^{-6}\).
The convergence rate for all cases is very close to the \(\alpha/(1-\alpha)\) rate predicted by (\ref{eq:linear-convergence-bound}), and varies very little with \(N_{t}\).
For larger \(\alpha\) the convergence rate is slightly lower than the theoretical prediction, while for very small \(\alpha\) the convergence rate is slightly higher than the theoretical prediction, likely as a result of the increased roundoff error of the diagonalisation.
Note that because the Richardson iteration requires one additional preconditioner application to calculate the initial preconditioned residual, the speedup prediction (\ref{eq:speedup}) is \(N_{t}/\left(2(m_{p}+1)\right)\) assuming \(T_{c} << T_{b}\).

\begin{figure}[t]
    \centering
    \epswitch{
    \includegraphics[width=0.99\linewidth]{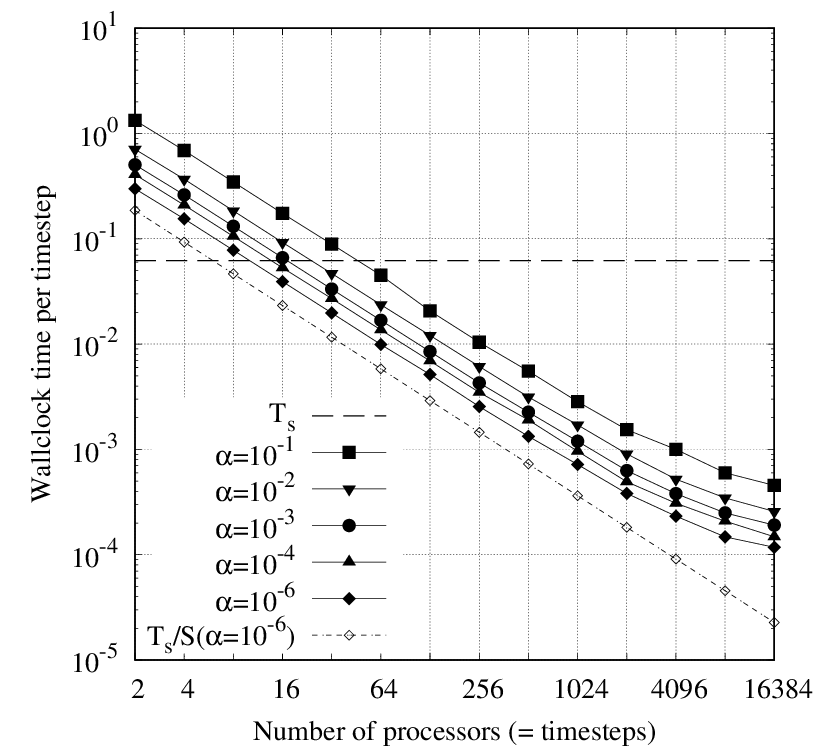}
    }{
    \includegraphics[width=0.99\linewidth]{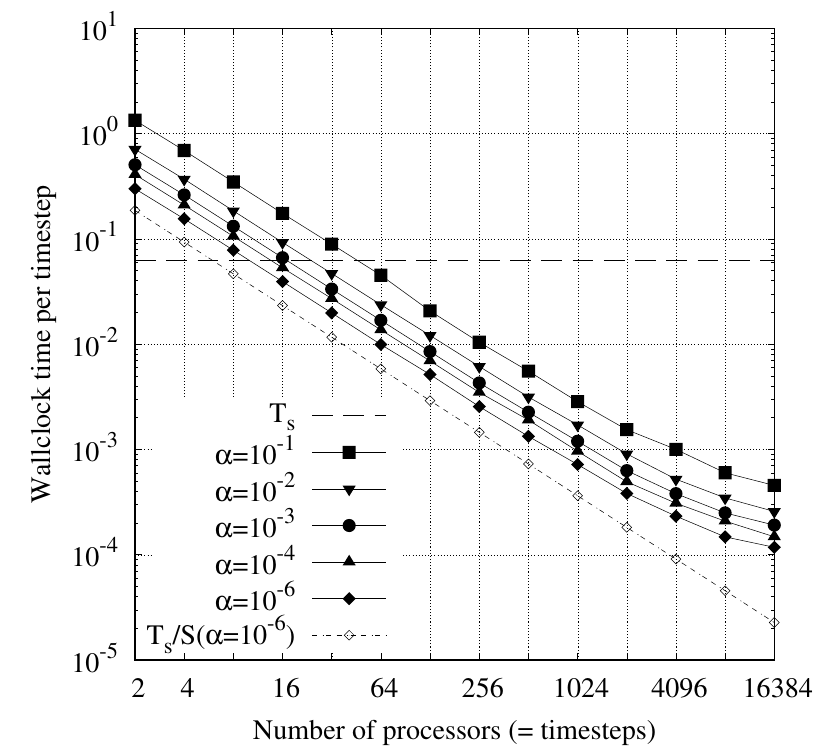}
    }
    \caption{Strong scaling in time of the linear scalar advection equation test case with \(N_{x}=128^{2}\), \(\sigma=0.8\)
    is demonstrated by measured wallclock times per timestep (\(T_{p}/N_{t}\)) for varying \(N_{t}\) and \(\alpha\).
    The wallclock time for the serial-in-time method is plotted as the horizontal line \(T_{s}\).
    The prediction of the performance model \(T_{s}/S\) (\ref{eq:speedup}) is shown for \(\alpha=10^{-6}\) using the measured \(T_{s}\) and \(S=N_{t}/\left(2(m_{p}+1)\right)\) (assuming \(T_{c}=0\)).
    Maximum speedup of 528 achieved with \(\alpha=10^{-6}\) and \(N_{t}=16384\).}
    \label{fig:advection-speedup}
\end{figure}

The measured wallclock times per timestep \(T_{p}/N_{t}\) for each value of \(\alpha\) are shown in Fig. \ref{fig:advection-speedup} compared to the serial-in-time method \(T_{s}\).
The corresponding speedups are shown in Table \ref{tab:advection-contraction-rates}.
For smaller \(N_{t}\) the available parallelism is not enough to outweigh the overhead of the parallel method, with the crossover to speedup occurring at \(N_{t}\approx64\) for \(\alpha=10^{-1}\) down to \(N_{t}\approx16\) for \(\alpha=10^{-6}\).
The scaling in time is almost perfect from \(N_{t}=2\) to \(N_{t}=1024\), achieving \(45-55\%\) of the ideal speedup \(S\) \eqref{eq:speedup} across this range compared to the measured \(T_{s}\).
For larger window sizes the speedup increases further but the scaling deteriorates, with maximum speedups at \(N_{t}=16384\) of 129 with \(\alpha=10^{-1}\) up to 517 with \(\alpha=10^{-6}\).
As \(\alpha\) decreases the speedup increases almost exactly proportionally to the reduction in the iteration count. {These speedup results are very competitive with previous state-of-the-art results for the constant coefficient scalar advection equation: 10-15 using ParaDiag with a collocation method \citep{vcaklovic2023parallel}; and up to 18 using MGRIT with optimized coarse grid operators \citep{De_Sterck_Falgout_Friedhoff_Krzysik_MacLachlan_2021}.
\citet{de2023fast} and \citet{De_Sterck_Falgout_Krzysik_Schroder_2023} achieved speedups of up to 12 on \textit{variable} coefficient advection using MGRIT with modified semi-Lagrangian coarse operators.
Although the variable coefficient case follows very naturally for the semi-Lagrangian method, for ParaDiag this case is more difficult because it requires a constant coefficient state in the circulant preconditioner.
In this article, problems with variable coefficients are demonstrated using the nonlinear examples below; variable coefficient linear equations are left for later work.}

\begin{figure}[t]
    \centering
    \epswitch{
    \includegraphics[width=0.99\linewidth]{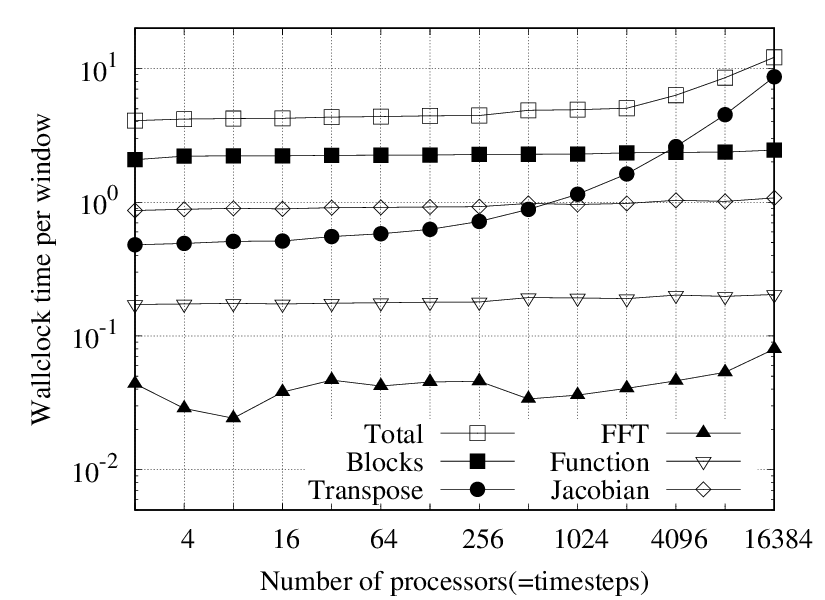}
    }{
    \includegraphics[width=0.99\linewidth]{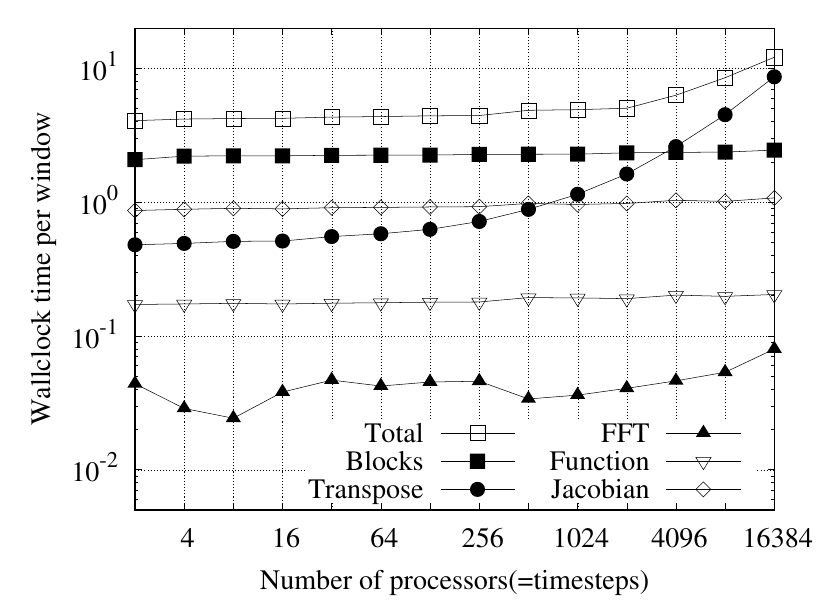}
    }
    \caption{Weak scaling in time profile of the linear scalar advection equations test case with \(N_{x}=128^{2}\), \(\sigma=0.8\), \(\alpha=10^{-4}\), and varying window length \(N_{t}\). Measured wallclock times per window for each component of the algorithm: total, block solves, space-time transpose, (I)FFTs, all-at-once function evaluation, and all-at-once Jacobian action.}
    \label{fig:advection-scaling}
\end{figure}

The ideal lower bound on the runtime \(T_{s}/S\) from \eqref{eq:speedup} using the measured
\(T_{s}\) is shown in Fig. \ref{fig:advection-speedup} for \(\alpha=10^{-6}\), labelled \(T_{s}/S(\alpha^{-6})\).
The measured speedup closely follows the prediction until the very longest window lengths.
To explore this scaling behaviour, the time taken by the most expensive sections of the computation are shown in Fig. \ref{fig:advection-scaling}.
Only the profile for \(\alpha=10^{-4}\) is shown as this is the recommended value from \citet{gander2019convergence}, and the profiles for the other values of \(\alpha\) are almost identical.
For clarity, the timings are shown \textit{per window solve} instead of \textit{per timestep} i.e. multiplied by \(N_{t}\) compared to the times in Fig. \ref{fig:advection-speedup}.
Shown are the total time per solve and the time taken for: the block solves; the space-time transpose; the (I)FFTs; evaluation of the all-at-once function (\ref{eq:nonlinear-aaos}); and the action of the all-at-once Jacobian (\ref{eq:nonlinear-jacobian}).
With perfect parallelism in time each of these components would take constant time independent of \(N_{t}\).

We first inspect the operations in the preconditioner.
The first point to notice is that the time spent in the block solves in Step 2 scales almost perfectly across the entire range of \(N_{t}\), as expected due to this step being embarrassingly parallel in time.
There is a very slight increase for the highest \(N_{t}\) due to a minor increase in the fill in of the LU factorisation but this effect is negligible.

The (I)FFTs are a minimal part of the profile, taking approximately 1\% of the total solution time for all window sizes.
The time taken actually decreases as \(N_{t}\) increases up to \(N_{t}=512\).
Each core calculates \(N_{x}/N_{t}\) transforms of length \(N_{t}\).
The total work required from each core is \(\mathcal{O}(N_{x}\text{log}N_{t})\), but it seems that the implementation used (\texttt{scipy.fft}) is more efficient for fewer longer transforms, leading to the initial decrease with \(N_{t}\), before the gradual increase after that.

Looking at the final component of the preconditioner, it is clear that the space-time transposes are the main culprit for the loss of scaling for large \(N_{t}\).
For \(N_{t}\leq32\) the entire problem fits on a single node and the communication time is almost constant.
For \(N_{t}\geq64\) internode communication is required and the communication time increases steadily as expected for an \texttt{Alltoall} communication.
At \(N_{t}=8192\) the transpose communications require more time than any other part of the computation, and at \(N_{t}=16384\) they take almost 60\% of the total solution time.

Outside of the preconditioner, the Jacobian action and the function evaluation have the same execution pattern: a one sided point-to-point communication round to update the time halos, followed by embarrassingly parallel-in-time computation of the result.
The Jacobian action takes the longer time but, just as for the block solves, both scale very well with \(N_{t}\) due to the majority of the time being spent on parallel computation and the communication complexity being independent of \(N_{t}\).

\subsection{Linear shallow water equations}\label{sec:examples:lswe}

The next example is the linearised shallow water equations on the rotating sphere.
The shallow water equations are obtained by assuming that the vertical (normal to the sphere) lengthscales of the flow are orders of magnitude smaller than the horizontal (tangential to the sphere) lengthscales and taking an average over the depth of the fluid.
The equations can be further simplified by linearising around a state of rest and a mean depth, which results in the following,
\begin{equation}\label{eq:lswe}
\begin{aligned}
    \partial_{t}\vec{u} + f\vec{u}^{\perp} + g\nabla h & = 0, \\
    \partial_{t}h + H \nabla\cdot\vec{u} & = 0,
\end{aligned}
\end{equation}
where \(\vec{u}\) is the velocity perturbation tangent to the sphere around the state of rest, \(f\) is the Coriolis parameter due to the rotating reference frame, \(\vec{u}^{\perp}=\vec{u}\times\hat{k}\) where \(\hat{k}\) is the unit normal vector to the sphere, \(g\) is gravity, and \(h\) is the depth perturbation around the mean depth \(H\).
This is an oscillatory PDE system with a hyperbolic part resembling a wave equation, plus the Coriolis term.

\begin{figure}
\begin{center}
\includegraphics[width=7cm]{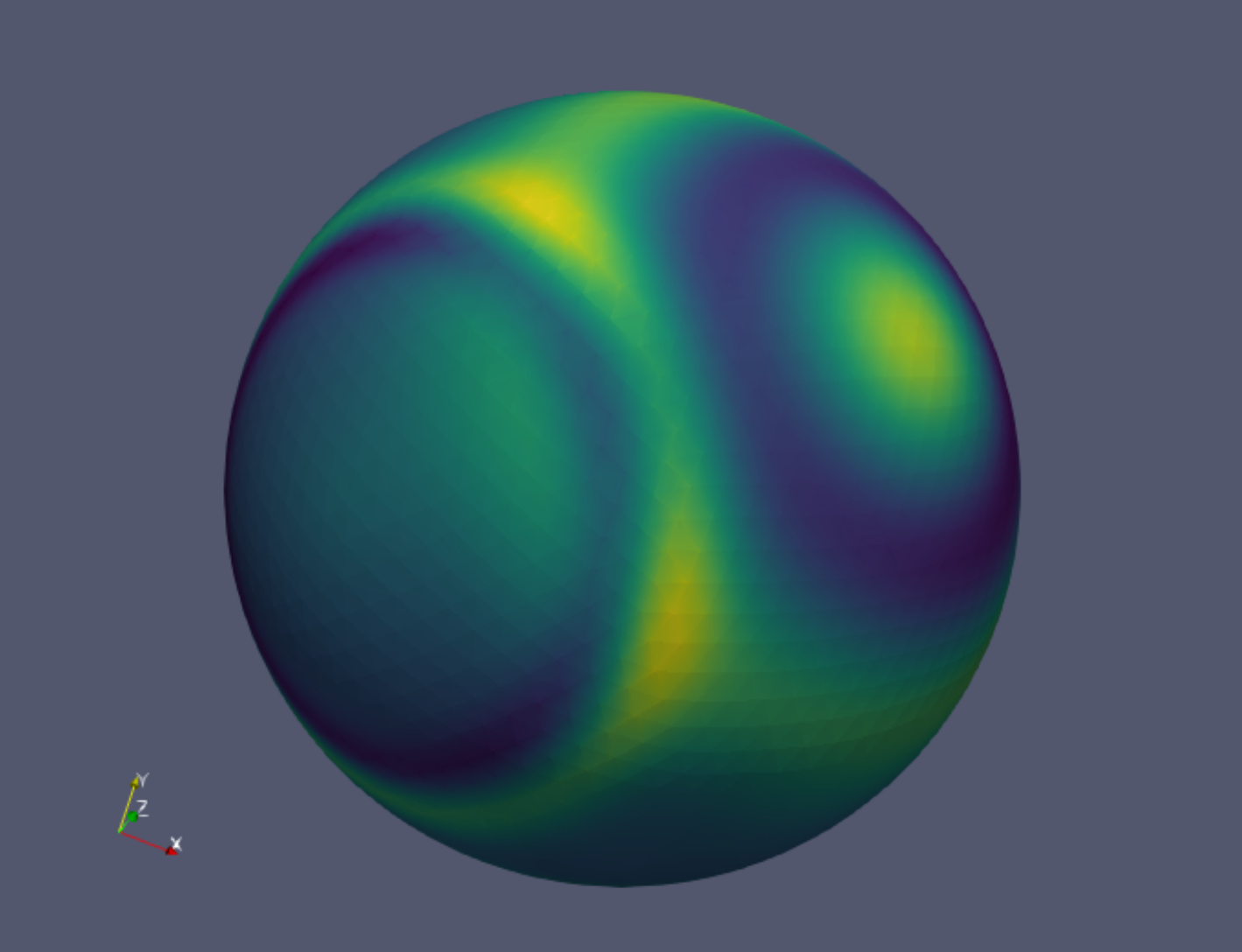}
\end{center}
\caption{\label{fig:gw}Colour plot of the free surface elevation \(h\) in the linear rotating shallow water equation testcase.}
\end{figure}

We use the test case of \cite{schreiber2019parallel} which simulates propagation of gravity waves around the Earth's surface and was used to test REXI parallel-in-time algorithms.
The initial conditions are zero velocity and three Gaussian perturbations to the depth.
As time advances, each Gaussian relaxes into a wave travelling around the globe in all directions, shown in Fig. \ref{fig:gw}.

\cite{schreiber2019parallel} show that the dispersion error of the trapezium rule integration method can lead to higher errors than the REXI schemes, however here we are primarily concerned with ParaDiag's ability to speed up the calculation of the base integrator, rather than with the ultimate accuracy of the base integrator.

An icosahedral sphere mesh with 5 refinement levels is used, resulting in \(\approx20k\) simplex elements.
The equations (\ref{eq:lswe}) are discretised using a compatible finite element method \citep{cotter2023compatible}.
We choose the quadratic Brezzi-Douglas-Marini elements \(V_{u}=\text{BDM}_{2}\) for the discrete velocity space \citep{Brezzi_Douglas_Marini_1985}, and the linear discontinuous Lagrange elements \(V_{h}=\text{DG}_{1}\) for the layer depth space.
This results in \(\approx215k\)DoFs for each timestep of the serial-in-time method.

The blocks in both the serial- and parallel-in-time methods are solved using hybridisation, which is a technique for exactly reducing the block matrix to a smaller finite-element space supported only on the facets between cells.
The smaller system is similar to a pressure Helmholtz equation, and is solved using LU decomposition with the \((\vec{u},h)\) solution obtained exactly by back-substitution.
For more detail on hybridisation and its implementation in Firedrake, see \citet{gibson2020slate}.
Using this technique, the block ``iteration'' counts are \(k_{s}=1\) and \(k_{p}=1\), as in the advection example. 
Although the blocks considered here are tractable with direct methods, full atmospheric models are not, so hybridisation or similar techniques are often used, and for this reason are interesting to consider for parallel-in-time applications.

The all-at-once system is solved using flexible GMRES \citep[FGMRES,][]{saad1993flexible} to a relative tolerance of \(10^{-11}\) for two \(\alpha\) values \(\alpha=10^{-4}\) and \(\alpha=10^{-6}\), converging in \(m_{p}=3\) and \(m_{p}=2\) iterations respectively independently of \(N_{t}\).\footnote{
   FGMRES is used here not for the flexibility (the preconditioner is constant) but because it requires one less preconditioner application than GMRES or Richardson iterations. For small iteration counts this outweighs the additional memory requirements and gives a significant speedup gain.}
We achieved greater speedup for \(\alpha=10^{-6}\) in the previous example, but we also show \(\alpha=10^{-4}\) here because it is closer both to the recommended values \citep{gander2019convergence}, and to those used for the nonlinear problems we consider next, where the predominant error stems from the reference value \(\hat{u}\) in the preconditioner so decreasing \(\alpha\) further does not lead to greater speedup.
The serial-in-time problem is parallelised in space with \(P_{s}=2\) processors, resulting in \(\approx108\)kDoFs per processor.
In the parallel-in-time method \(P_{p}=2P_{s}N_{t}\) processors are used to maintain the same number of floating points per processor for the complex-valued block solves.

\begin{figure}[t]
    \centering
    \epswitch{
    \includegraphics[width=0.99\linewidth]{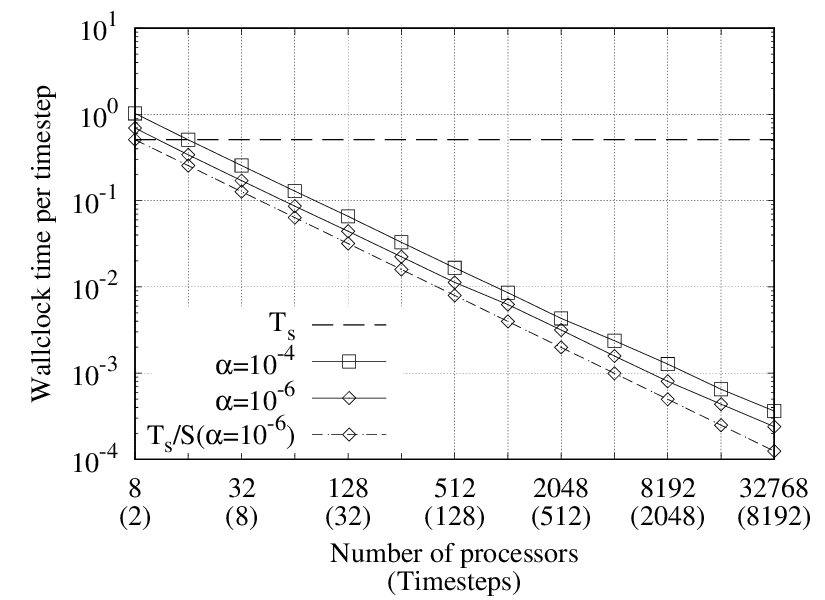}
    }{
    \includegraphics[width=0.99\linewidth]{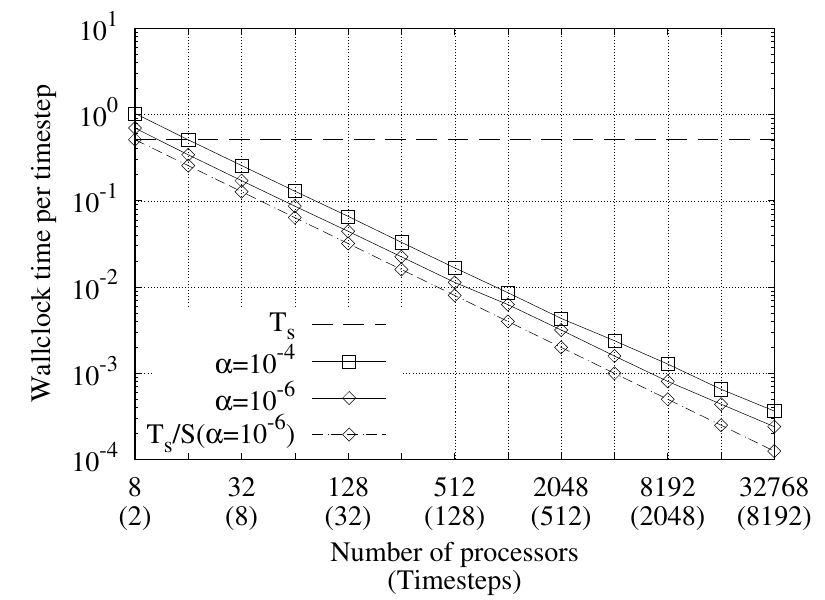}
    }
    \caption{Strong scaling in time of the linear shallow water equations on the rotating sphere with 20k elements, \(dt=900s\), and varying window length \(N_{t}\) for \(\alpha=10^{-4}\) and \(\alpha=10^{-6}\). \(P_{s}=2\) for the serial-in-time method and \(P_{p}=2N_{t}P_{s}\) for the parallel-in-time method. Measured wallclock times per timestep \(T_{p}\), as well as the predictions of the performance model \(T_{s}/S\) (\ref{eq:speedup}) for \(\alpha=10^{-6}\) using the measured \(T_{s}\) (assuming \(T_{c}=0\)). The wallclock time for the serial-in-time method is plotted as the horizontal line. Maximum speedups of 1,402 and 2,126 with \(\alpha=10^{-4}\) and \(\alpha=10^{-6}\) respectively, both at \(P_{p}=32768\) and \(N_{t}=8192\).}
    \label{fig:lswe-speedup}
\end{figure}

The strong scaling of wallclock time taken per timestep is shown in Fig. \ref{fig:lswe-speedup} up to \(N_{t}=8,192\) and \(P_{p}=32,768\), with the performance model prediction for \(\alpha=10^{-6}\).
The crossover to the parallel method being faster occurs at \(N_{t}=4\) for both cases, earlier than for the scalar advection example due to using \(2P_{s}\) processors for each block solve.
Again we see excellent scaling for both cases, achieving 63\%-74\% of the ideal speedup from \(N_{t}=2\) to \(N_{t}=1024\), where the actual speedups are 214 and 323 using 4096 processors.
For \(N_{t}>1024\) the scaling deteriorates slightly, but still reaches maximum speedups over the serial-in-time method of 1,402 and 2,126 with \(N_{t}=8192\) for a total of \(\approx1.75\) billion DoFs using 32,768 processors (512 nodes), or 51\% of the ideal speedup for each value of \(\alpha\).
This speedup is very competitive with previous results in the literature for time-parallel speedup of the linear shallow water equations \citep{schreiber2018beyond}.

ParaDiag can clearly scale excellently and achieve very competitive speedups for linear equations, even hyperbolic ones that have previously proved challenging for parallel-in-time methods.
With reference to the performance model (\ref{eq:speedup}), this scaling relies on three factors.
Firstly, the number of all-at-once iterations \(m_{p}\) is fixed independently of \(N_{t}\).
This holds for constant coefficient equations such as we have considered so far because no averaging is necessary in the preconditioner so the \(\alpha/(1-\alpha)\) convergence rate (\ref{eq:linear-convergence-bound}) is achieved independently of \(N_{t}\).
Secondly, the number of iterations for the complex-valued block solves \(k_{p}\) is fixed independently of \(N_{t}\).
So far we have achieved this either by directly solving the blocks, or using hybridisation to reduce the block system down to a size where a direct solver can be applied.
Lastly, the collective communications required for the space-time transpose must scale well.
In the examples so far this has been the limiting factor on the scaling, although only once \(N_{t}\) is \(\mathcal{O}(10^{3})\).

\subsection{Nonlinear shallow water equations}\label{sec:examples:swe}

\begin{figure}
\begin{center}
\includegraphics[width=0.7\linewidth]{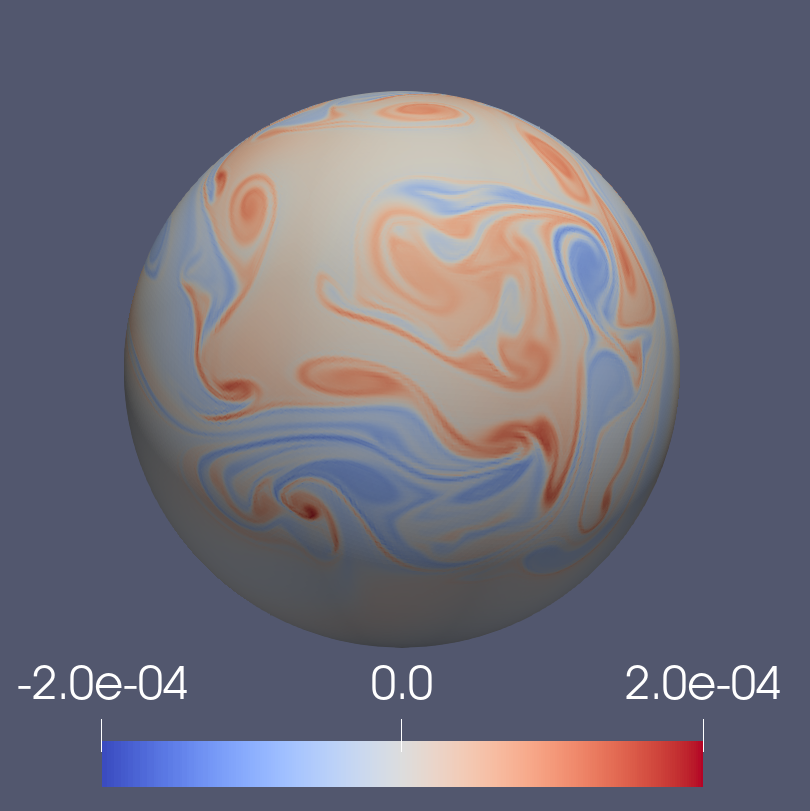}
\end{center}
\caption{\label{fig:swe}Colour plot of the potential vorticity in the rotating shallow water equation testcase.}
\end{figure}

The nonlinear shallow water equations on the rotating sphere retain the nonlinear advection terms which are removed in the linear shallow water equations.
Writing the velocity advection in vector invariant form, the nonlinear equations are
\begin{equation}\label{eq:swe}
\begin{aligned}
    \partial_{t}\vec{u} + (\nabla\times\vec{u} + f)\vec{u}^{\perp} + \frac{1}{2}\nabla|\vec{u}|^{2} + g\nabla h & = 0, \\
    \partial_{t}h + \nabla\cdot(\vec{u}h) & = 0.
\end{aligned}
\end{equation}
The same compatible finite element spaces are used as for the linear shallow water equations, using 
the discretisation presented in \citet{gibson2019compatible}.

\begin{table}[t]
    \caption{Details for 4 refinement levels of the icosahedral sphere mesh used for the unstable jet test case. DoFs includes both the velocity and depth DoFs. \(\Delta x\) is the element edge length and \(\Delta t\) is chosen for a CFL of approximately \(\sigma\approx0.4\). \(P_{s}\) is the number of processors used in the serial-in-time method.\label{tab:sphere-mesh}}
    \centering
    \begin{tabular}{ccccc}
      \toprule
      Mesh &  DoFs (\(\times10^3\)) & \(\Delta x\) (km) & \(\Delta t\) (s) & \(P_{s}\) \\
      \midrule
         4    &   53  & 480  & 900 &  1  \\
         5    &  215  & 240  & 450 &  4  \\
         6    &  860  & 120  & 225 & 16  \\
      7    & 3440  &  60  & 112 & 64  \\
      \bottomrule
    \end{tabular}
\end{table}

We use the classic test case from \citet{galewsky2004initial} of an unstable jet in the Earth's northern hemisphere which breaks down into a highly nonlinear vortical flow over a period of several days.
The strong nonlinearities provide a test of the error introduced by the averaging procedure in the circulant preconditioner.
The flow is evolved from the initial conditions to a final time of 10 days at a range of resolutions and window lengths \(N_{t}\).
The potential vorticity field at 10 days is shown in Fig. \ref{fig:swe}.
All window lengths tested were much smaller than the total number of timesteps required so multiple windows are solved sequentially until the final time is reached, using the last timestep of one window as the initial condition for the next window.

The trapezium rule is used in time and the timestep is scaled with the mesh resolution to keep the advective CFL close to constant around \(\sigma\approx0.4\) across all tests.
The blocks in both the serial- and parallel-in-time methods are solved using the FGMRES Krylov method preconditioned with a geometric multigrid.
The relaxation on each grid level is 4 GMRES iterations preconditioned with a vertex Vanka patch smoother.
This is an additive Schwarz method where each subdomain is associated with a mesh vertex and contains all DoFs in the closure of the star of the vertex.
The coarsest level is solved directly using MUMPS.
For the serial-in-time case this method gives resolution independent convergence rates for fixed CFL.
The patch smoothers are implemented in PETSc with \texttt{PCPATCH} and exposed in Firedrake through the \texttt{firedrake.PatchPC} interface \citep{farrell2021pcpatch}.

\begin{figure}[t]
    \centering
    \epswitch{
    \includegraphics[width=0.99\linewidth]{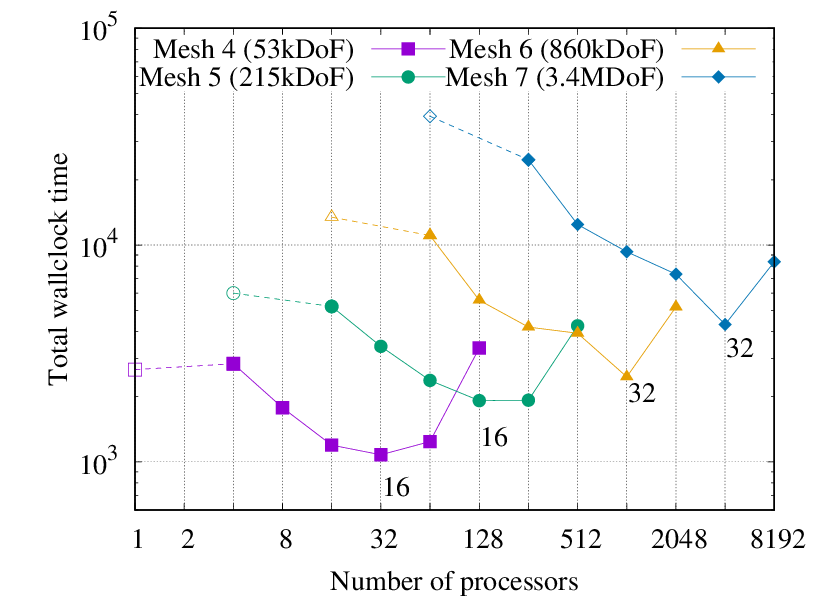}
    }{
    \includegraphics[width=0.99\linewidth]{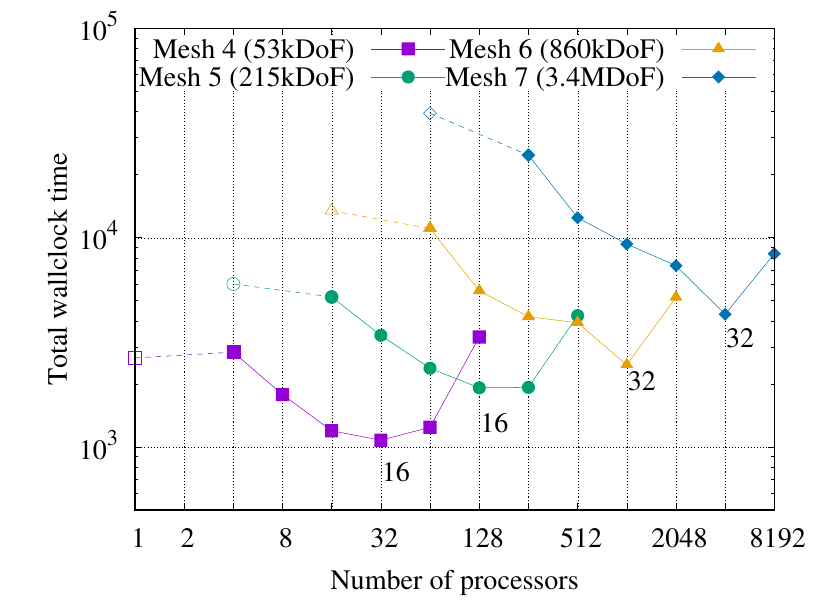}
    }
    \caption{Strong scaling in time for the unstable jet shallow water test case at varying mesh resolutions. Measured wallclock time to simulate 10 days. On each curve: the leftmost data point with an open symbol shows the serial-in-time (but possibly parallel-in-space) result; the second data point onwards show parallel-in-time results at increasing window lengths \(N_{t}\) from 2-64. The number of processors for the parallel-in-time results is \(P_{p}=2N_{t}P_{s}\). Maximum speedups are 2.47, 3.14, 5.45, and 9.12 at refinement levels 4, 5, 6, 7 respectively. Labels indicate \(N_{t}\) at the maximum speedup on each mesh.}
    \label{fig:swe-strong-scaling}
\end{figure}

Each timestep in the serial-in-time method and each window in the parallel-in-time method is solved using an inexact Newton method with the tolerance for the linear solves at each Newton step determined using the first Eisenstat-Walker method \citep{Eisenstat_Walker_1996}.
In the serial-in-time method the convergence criteria is a nonlinear residual below an absolute tolerance of \(10^{4}\), equivalent to around a \(10^{-6}\) relative reduction in the residual.\footnote{The absolute residual values of the residuals are so high because they are calculated by integrating over the surface of the Earth, which is \(\mathcal{O}(10^{14}m^{2})\).}
In the parallel-in-time method a nonlinear residual below \(\sqrt{N_{t}}10^{4}\) is required, i.e. an average residual of \(10^{4}\) at each timestep.
In the linear examples we converged to very tight tolerances, which are useful to confirm the convergence of the method but are impractical because, in reality, discretisation errors will be orders of magnitude larger.
In practice, the convergence tolerance need only be tight enough for the method to be stable, and to reach the expected spatial and temporal convergence rates.
The linear solve at each Newton step of the parallel-in-time method is solved using FGMRES using the circulant preconditioner and a value of \(\alpha=10^{-4}\), and the complex-valued blocks are solved to a relative tolerance of \(\tau=10^{-3}\).
The circulant preconditioner will not achieve the expected convergence unless the complex-valued blocks are solved to a relative tolerance \(\tau\) of \textit{at least} as small as \(\alpha/N_{t}\) (and often smaller) \citep{vcaklovic2023parallel}.
For the linear examples we used a direct solver on the blocks, equivalent to a relative tolerance close to machine precision.
For operational models direct solvers are unfeasible, both because of the large size of the system and because, for nonlinear systems, a refactorisation would be necessary at each Newton iteration.
Here, we use an iterative solver for the blocks so reducing \(\alpha\) requires more expensive block solves.
Empirically we have found that setting \(\alpha\) and \(\tau\) in the range \(10^{-5}-10^{-3}\) gives a good balance between inner/outer Krylov iterations, but exploring and analysing this choice further is left for later work.

Four refinement levels of the icosahedral sphere mesh are used ranging from 50-3,440kDoFs.
See Table \ref{tab:sphere-mesh} for details on each mesh.
The DoFs/processor is fixed at \(\approx\)50kDoFs in all cases, requiring \(P_{p}=2N_{t}P_{s}\) in the parallel-in-time case.

\begin{figure}[t]
    \centering
    \epswitch{
    \includegraphics[width=\linewidth]{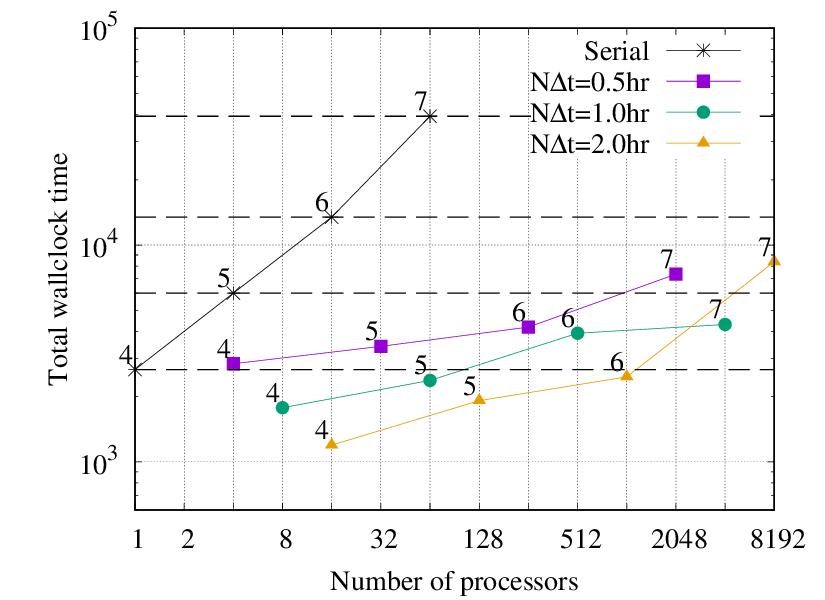}
    }{
    \includegraphics[width=\linewidth]{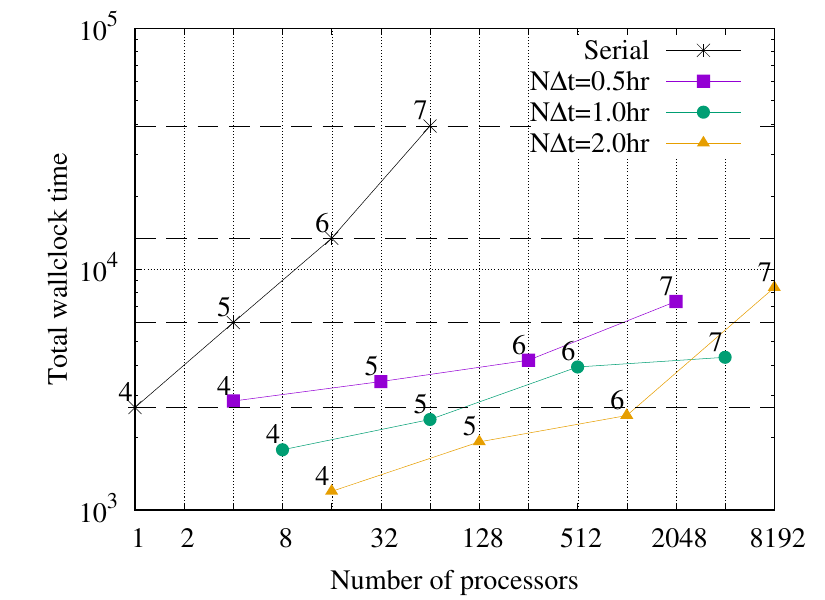}
    }
    \caption{Weak scaling in time and space for the unstable jet shallow water test case. Measured wallclock time to simulate 10 days at various resolutions. Lines are plotted for increasing resolutions at constant window length \(T=N_{t}\Delta t\), i.e. for each mesh refinement \(\Delta t\) and \(\Delta x\) are halved and \(N_{t}\) doubled. Data points are labelled with the mesh resolution. \(N_{t}=2,\;4,\;8\) for \(T=0.5,\;1,\;2\) respectively on mesh 4. The number of processors for the parallel-in-time results is \(P_{p}=2N_{t}P_{s}\).}
    \label{fig:swe-weak-scaling}
\end{figure}

We first show strong scaling of the wallclock time for 10 simulation days at each mesh resolution for increasing window sizes in Fig. \ref{fig:swe-strong-scaling}.
As for the linear problems there is little or no speedup for very small \(N_{t}\), followed by a close to linear increase in the speedup at \(N_{t}\) grows.
However, around \(N_{t}=16-32\) the speedup reaches maximums of between 2.47 and 9.12 for the lowest and highest resolution respectively over the equivalent serial-in-time method.
Although these speedups are significantly less than those achieved for linear problems, they are competitive with current state of the art parallel-in-time methods for the shallow water equations, without any modification to the basic ParaDiag algorithm.

\begin{figure}[t]
    \centering
    \epswitch{
    \includegraphics[width=\linewidth]{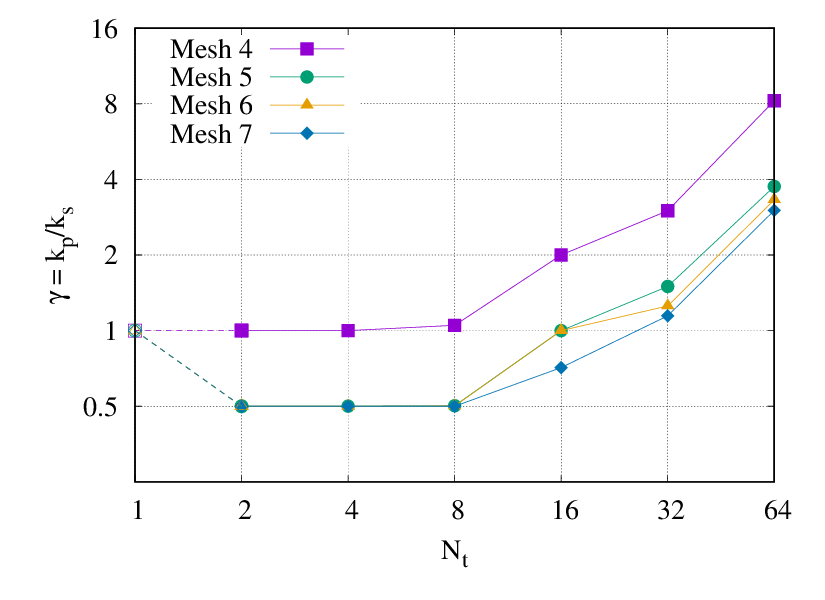}
    }{
    \includegraphics[width=\linewidth]{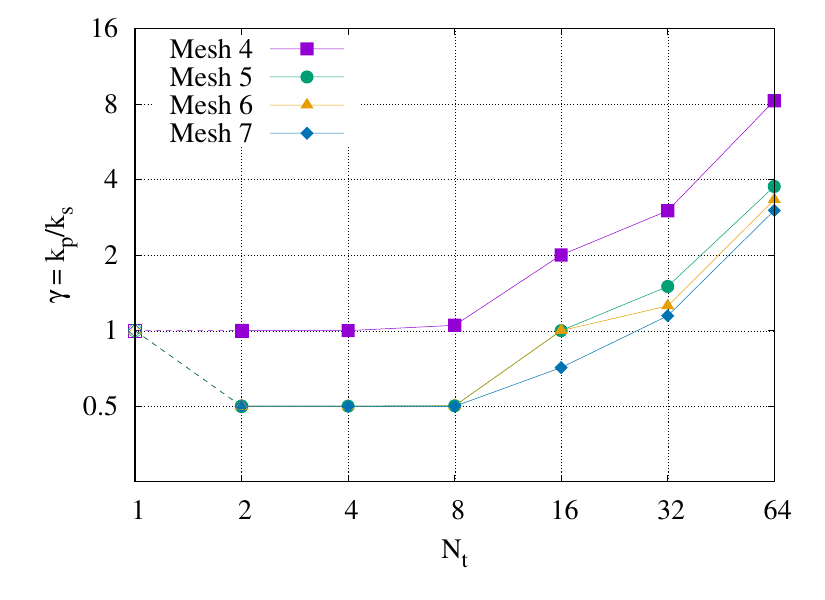}
    }
    \caption{Scaling of the block iteration counts for the nonlinear shallow water equations \(\gamma=k_{p}/k_{s}\), with maximum iteration counts for the complex-valued blocks in the circulant preconditioner \(k_{p}\) and block iteration counts for the real-valued blocks in the serial-in-time method \(k_{s}\). The ratio \(\gamma\) is plotted versus the window length \(N_{t}\) for each mesh resolution. Data point at \(N_{t}=1\) is \(\gamma_{s}=k_{s}/k_{s}=1\).}
    \label{fig:swe-gamma}
\end{figure}

\begin{figure}[t]
    \centering
    \epswitch{
    \includegraphics[width=\linewidth]{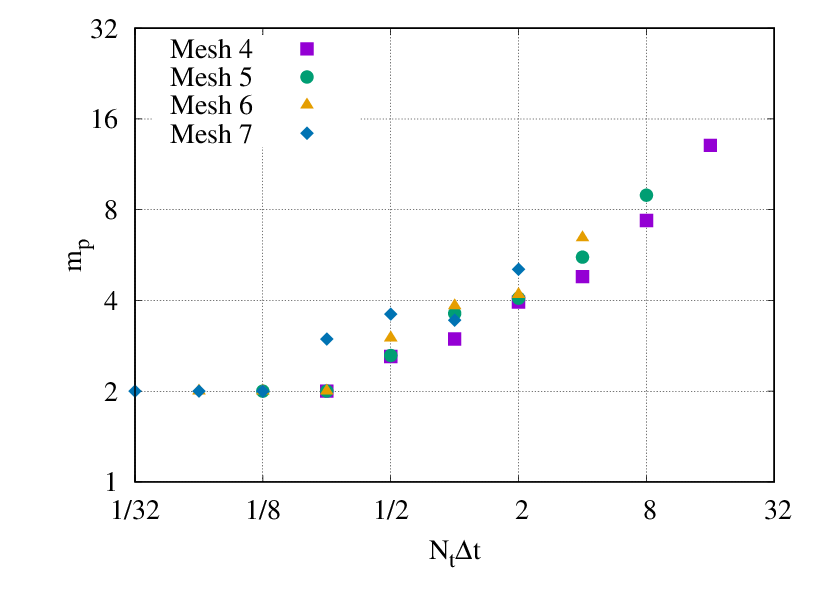}
    }{
    \includegraphics[width=\linewidth]{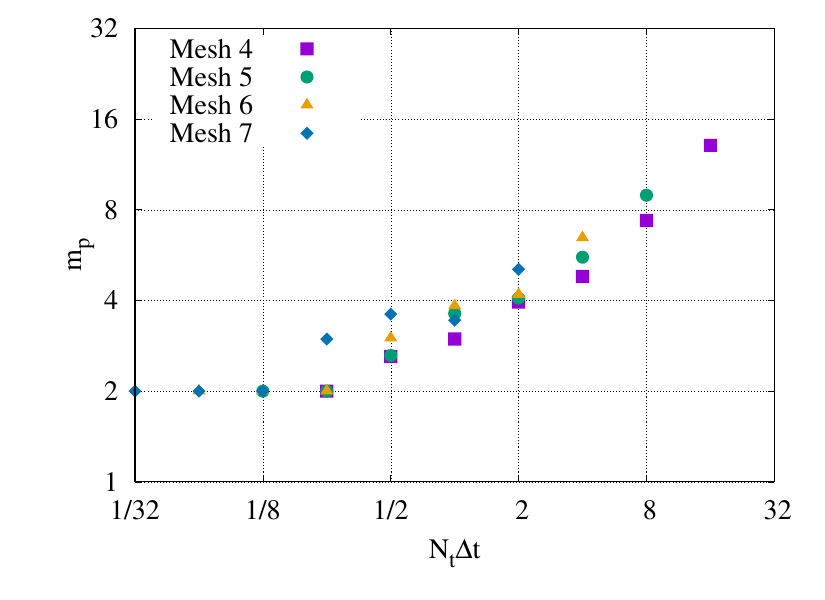}
    }
    \caption{Scaling of the all-at-once iteration counts for the nonlinear shallow water equations. We show iteration counts \(m_p\) for all mesh resolutions and window sizes versus the window duration \(N_{t}\Delta t\) in hours, where \(\Delta t\) depends on the mesh resolution due to the constant CFL.}
    \label{fig:swe-omega}
\end{figure}

For all mesh resolutions there is a sudden and drastic reduction in the speedup at \(N_{t}=64\).
To understand this behaviour we return to the performance model in Sect. \ref{sec:methods:performance-model}.
For all four mesh resolutions, the serial-in-time method required two Newton steps (\(m_{s}=2\)), with one multigrid cycle per step for resolution 4 and two cycles per step for all other resolutions (i.e. \(k_{s}=1\) or \(2\)).
The block iteration parameter \(\gamma=k_{p}/k_{s}\) is shown in Fig. \ref{fig:swe-gamma}, and the number of linear all-at-once iterations required per window  \(m_{p}\) is shown in Fig. \ref{fig:swe-omega}.
Both iteration counts increase with \(N_{t}\), eventually resulting in the performance degradation seen in Fig. \ref{fig:swe-strong-scaling}.
We now examine each parameter in more detail.

In Fig. \ref{fig:swe-gamma}, the increase in the iterations is due to the clustering of the circulant eigenvalues close to the imaginary axis shown in Fig. \ref{fig:circulant-eigenvalues}.
For small window lengths the blocks for meshes 5, 6, and 7 only require a single iteration to reach the required tolerance versus two for the serial-in-time block, which contributes significantly to the improved speedups for higher resolutions.
This lower iteration count is achieved because, although the block coefficients are less favourable in the complex-valued case, we do not require a very tight tolerance.
The major error in the preconditioner is the averaging error, discussed below, which means that as long as the error from the inexact block solves is below the averaging error there is no further gain in the convergence of the all-at-once system from solving the blocks to a tighter tolerance.

In Fig. \ref{fig:swe-omega} the number of linear all-at-once iterations across all Newton steps per window \(m_p\) is plotted against the total window duration \(T = N_{t}\Delta t\).
We see that the data for the different resolutions collapse onto each other and, after an initial plateau at \(m_{p}=2\), the number of all-at-once iterations scales close to linearly with \(T\).
This is consistent with a convergence rate scaling of \(\kappa N_{t}\Delta t\) for the averaged preconditioner \citep{caklovic_paradiag_2023}, assuming that the Lipschitz constant \(\kappa\) is a function of the underlying continuous problem and so is essentially constant with resolution.
This scaling is the primary reason why the maximum speedup increases with mesh resolution in Fig. \ref{fig:swe-strong-scaling}: in order to balance spatial and temporal errors, \(\Delta t\) is decreased proportionally to \(\Delta x\) to maintain a constant CFL, so for the same \(N_{t}\) the higher resolution meshes have better convergence rates of the all-at-once system.\footnote{We note that, even at the highest resolution, our results are coarser than the temporal resolution \(\Delta t=30s\) in the original article \citep{galewsky2004initial}. The improved speedup seen is a product of the convergence behaviour of the averaged preconditioner rather than of the ``over-resolution'' phenomena identified in \citet{gotschel2020twelve}.}

Next we investigate the effectiveness of weak scaling in both time and space.
Previously we were interested in decreasing the wallclock time at a particular resolution.
Now we are interested in solving to a fixed final time and increasing the spatial and temporal resolution simultaneously without increasing the wallclock time.
This represents the situation where we would like to provide a more accurate forecast of a particular length at a particular delivery frequency.

We investigate the weak scaling for the unstable jet case by reinterpreting the data in Fig. \ref{fig:swe-strong-scaling}, selecting a particular window duration \(T=N_{t}\Delta t\) and keeping this fixed as the resolution is increased i.e. each time the mesh is refined \(\Delta x\) and \(\Delta t\) are halved and \(N_{t}\) is doubled.
The results for window durations of \(0.5\), \(1\), and \(2\) hours are shown in Fig. \ref{fig:swe-weak-scaling}, alongside the wallclock times for the serial-in-time method weak scaled only in space.
As expected from the performance model (\ref{eq:time-taken-serial}) the wallclock time for the serial-in-time method increases linearly with resolution despite almost perfect weak scaling in space because the total number of timesteps increases.
On the other hand, the parallel-in-time method scales much more favourably, with the exception of the highest resolution (mesh 7) data point on the longest time window \(T=N_{t}\Delta t=2hr\), which corresponds to the \(N_{t}=64\), \(P_{p}=8192\), data point on the mesh 7 line in Fig. \ref{fig:swe-strong-scaling}.
Figure \ref{fig:swe-weak-scaling} shows that by exploiting time-parallelism the resolution can be increased fourfold in both time and space without increasing the total time to solution.
This can be seen from the mesh 6 solution with \(T=2hr\) requiring the same wallclock time as the serial-in-time mesh 4 solution, or the mesh 7 solution with \(T=1hr\) (and almost with \(T=0.5hr\)) requiring less wallclock time than the serial-in-time mesh 5 solution.
The all-at-once iteration scaling in Fig. \ref{fig:swe-omega} is a key component in achieving this weak scaling behaviour because it shows that the convergence of the all-at-once system with the time averaged preconditioner is independent of both spatial and temporal resolution for fixed \(T\).

\subsection{Compressible Euler equations}\label{sec:examples:euler}

\begin{figure}
\begin{center}
\includegraphics[width=\linewidth]{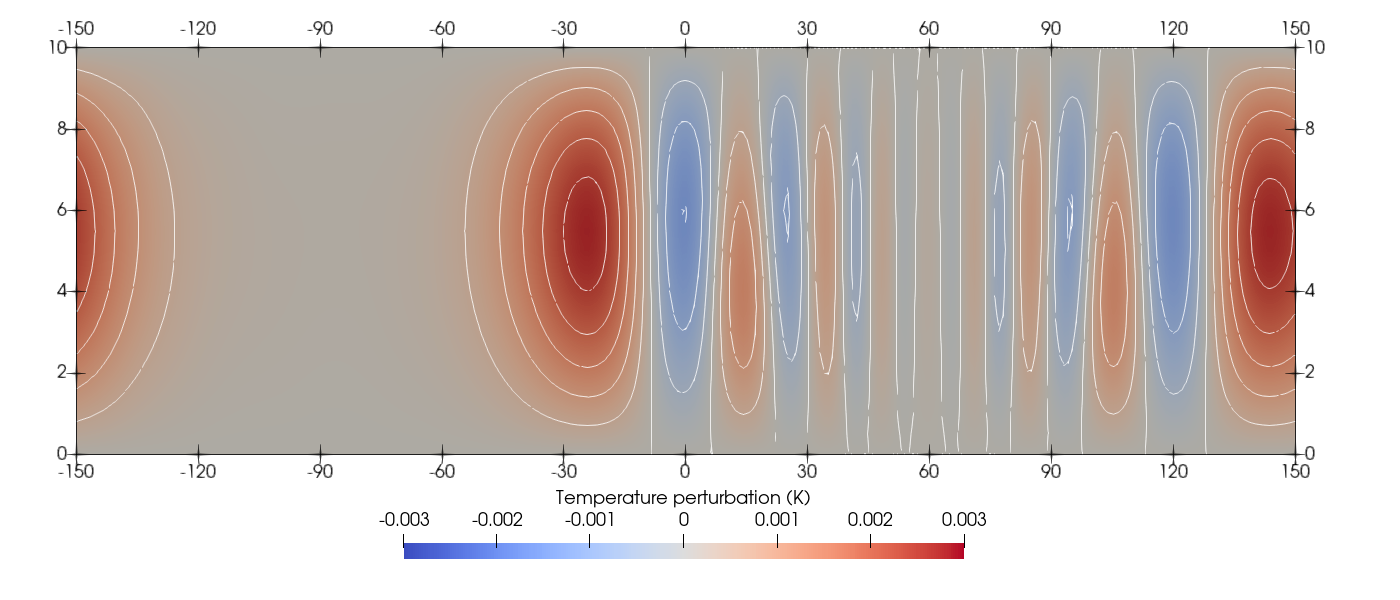}
\end{center}
\caption{\label{fig:euler-temperature}Colour plot of the temperature variation around the back ground state for the compressible Euler internal gravity wave example after 3000 seconds. Contours are drawn every \(5\times10^{-4}K\) and axes labels are in km. Vertical axis scaled by \(\times10\) for clarity.}
\end{figure}

The final example we present is the compressible Euler equations restricted to a 2D vertical slice of the atmosphere.
This model is an important step in the hierarchy of atmospheric models because, unlike the shallow water equations which only support surface gravity waves, the vertical slice model supports internal gravity waves and acoustic waves \citep{Melvin_Dubal_Wood_Staniforth_Zerroukat_2010}.
The compressible Euler equations with prognostic variables of velocity \(\vec{u}\), density \(\rho\), and potential temperature \(\theta\), with the equation of state for the Exner pressure \(\Pi\) are:

\begin{equation}\label{eq:euler}
\begin{aligned}
    \partial_{t}\vec{u} + \vec{u}\cdot\nabla\vec{u} + f\vec{u}^{\perp} + c_{p}\theta\nabla\Pi + g\hat{k} & = 0 \\
    \partial_{t}\theta + \vec{u}\cdot\nabla\theta &= 0 \\
    \partial_{t}\rho + \nabla\cdot(\vec{u}\rho) &= 0 \\
    \Pi^{(1-\kappa)/\kappa} &= R\rho\theta/p_{0}
\end{aligned}
\end{equation}
where \(c_{p}\) is the specific heat at constant pressure, \(R\) is the gas constant, \(\kappa=R/c_{p}\), \(p_{0}\) is the reference temperature, and \(f\) and \(\hat{k}\) are again the Coriolis parameter and unit vector in the vertical direction respectively.
We use the compatible finite element formulation proposed in \citet{Cotter_Shipton_2023} on a Cartesian mesh, where the finite element spaces are defined as tensor products of 1D elements in the horizontal (\(x\)) and vertical (\(z\)) directions.
This tensor product structure is enabled by Firedrake's \texttt{ExtrudedMesh} functionality, where a 1D horizontal base mesh along the ground is first defined, then extruded in the vertical direction \citep{McRae_Bercea_Mitchell_Ham_Cotter_2016,Bercea_McRae_Ham_Mitchell_Rathgeber_Nardi_Luporini_Kelly_2016}.
More details on the properties of this discretisation can be found in \citet{Cotter_Shipton_2023}.
%%% TODO: add this back in %%% The weak form of (\ref{eq:euler}) is shown in \gmd{the supplementary materials}{Appendix \ref{app:vertical-slice}} and more details on the properties of this discretisation can be found in \citet{Cotter_Shipton_2023}.
Results will be shown for a classic test case from \citet{Skamarock_Klemp_1994} of a gravity wave propagating through a uniform background velocity in a nonhydrostatic regime with \(f=0\) and the flow restricted to in-plane motion only.
This test case has been previously used to test the vertical slice discretisations in \citet{Cotter_Shipton_2023} and \citet{Melvin_Dubal_Wood_Staniforth_Zerroukat_2010}.
The domain has a width and height of 300km\(\times\)10km and periodic boundary conditions in the horizontal direction.
The flow is initialised with background temperature and density profiles \(\rho_{ref}\) and \(\theta_{ref}\) in hydrostatic balance, and a uniform horizontal velocity \(\vec{u}_{ref}\).
In the centre of the domain, a small amplitude temperature perturbation with lengthscale 5km is added to the background state, which creates an internal gravity wave that propogates to the right and left.
The temperature perturbation around the background state is shown in Fig. \ref{fig:euler-temperature} after 3000 seconds.

The domain is discretised with a resolution of 1km in both the vertical and horizontal directions, giving 49.2kDoFs per timestep, and the timestep size is \(\Delta t = 12s\) giving an advective Courant number of \(\sigma_{u}\approx0.24\).
Both the serial-in-time and parallel-in-time methods are run for 1024 timesteps, with window sizes of \(N_{t}=2-128\) with the parallel-in-time method (this gives a simulation time longer than the standard end time of 3000s, but results in a minimum of 8 windows with the parallel-in-time method).
For the serial-in-time method \(P_{s}=2\) and for the parallel-in-time method \(P_{p}=2N_{t}P_{s}\).

\citet{Cotter_Shipton_2023} used the trapezium rule for serial-in-time integration, with the real-valued block solved using an additive Schwarz method with columnar patches.
The patches are constructed as the vertical extrusion of a star patch at each vertex of the base mesh i.e. collecting all DoFs associated with a vertical column of facets and the interior DoFs of the neighbouring elements and solved using a direct solver.
\citet{Cotter_Shipton_2023} found that this method has iteration counts that are independent of resolution for fixed Courant number, but grow linearly with the Courant number at fixed resolution.
When increasing resolution in practice, space and time would usually be refined simultaneously, so constant iteration counts for fixed Courant number is a favourable result.
This preconditioner is implemented using Firedrake's \texttt{ASMStarPC} class, which in turns wraps PETSc's PCASM preconditioner, and makes use of the fact that \texttt{ExtrudedMesh} retains the topological information of the base mesh to identify the columns.

The first set of results for this example uses the star patch preconditioner for both the real- and complex- valued blocks.
The patch preconditioner is constructed from the constant-in-time reference state with non-zero velocity (\(\mathbf{u}_{ref}, \rho_{ref}, \theta_{ref})\), which avoids rebuilding the patches through the simulation.
Compared to constructing the patches from the current state at each timestep, constructing the patches from this reference state gives almost identical iteration counts because the gravity wave variations around the background state are relatively small.
For other vertical slice test cases the patches could be constructed around the current state.
In the serial-in-time method each timestep is solved using Newton-FGMRES with the first Eisenstat-Walker adaptive convergence criteria to an absolute tolerance of \(10^{-4}\), giving a relative residual drop of approximately \(10^{-6}\).
With the patch preconditioner, each timestep converged with an average of 27.5 linear iteration across 2 quasi-Newton iterations with \(\Delta t=12s\).
The parallel-in-time method is solved using Newton-FGMRES iterations to an absolute tolerance of \(\sqrt{N_{t}}10^{-4}\) so that the average residual of each timestep matches the serial-in-time method.
The circulant preconditioner uses \(\alpha=10^{-5}\), and each block in Step 2 is solved to a relative tolerance of \(\tau=10^{-5}\).

\begin{table}[t]
    \caption{Maximum wallclock speedups achieved by the parallel-in-time method with three different block solution methods (star patch, or the composite preconditioner with either fixed \(\tau\) or fixed \(k_{p}\)) over the serial-in-time method with two different block solution methods (star patch or the composite preconditioner). The serial-in-time method with the star patch was \(\approx1.3\times\) faster than with the composite preconditioner.\label{tab:euler-speedups}}
    \centering
    \begin{tabular}{ccc}
      \toprule
                                     & \multicolumn{2}{c}{Real-valued block method}  \\
      \cmidrule(lr){2-3}
      \multicolumn{1}{c|}{Complex-valued block method} & Star patch & Composite
      \\
      \midrule
         \multicolumn{1}{c|}{Star patch}                        &                      2.57 &  3.34 \\
         \multicolumn{1}{c|}{Composite, \(\tau=10^{-5}\)}       &                      6.14 &  7.99 \\
         \multicolumn{1}{c|}{Composite, \(k_{p}=\sqrt{N_{t}}\)} &                      8.64 & 11.24 \\
      \bottomrule
    \end{tabular}
\end{table}

\begin{figure}[t]
    \centering
    \epswitch{
    \includegraphics[width=0.99\linewidth]{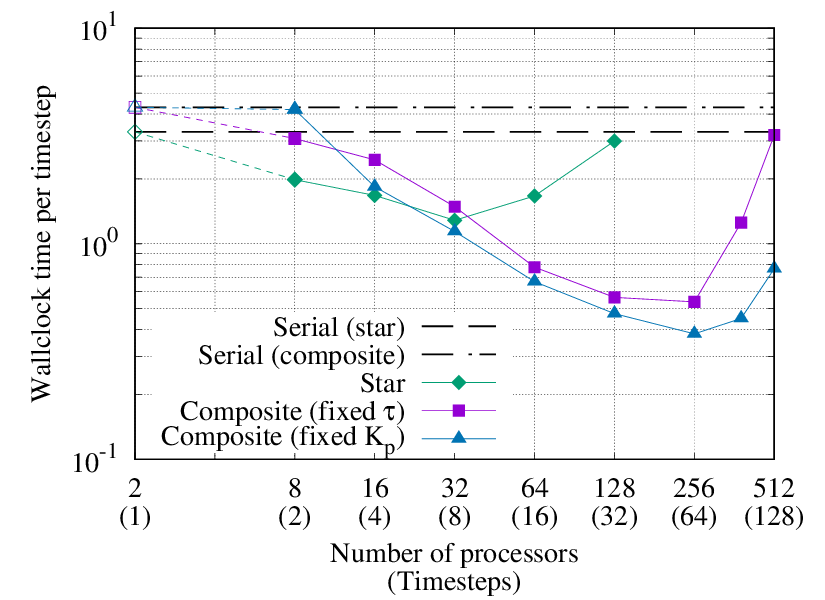}
    }{
    \includegraphics[width=0.99\linewidth]{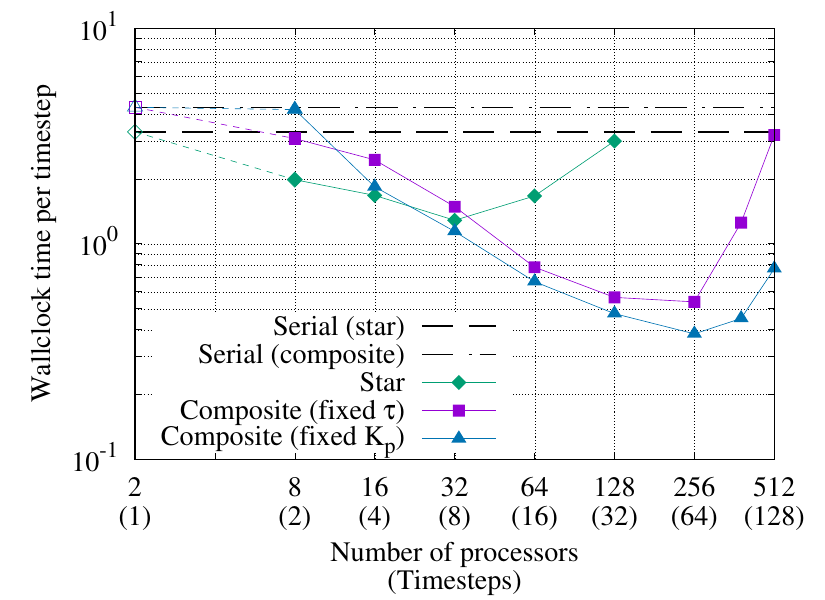}
    }
    \caption{Strong scaling in time of the compressible Euler equations gravity wave test case for varying window length \(N_{t}\). \(P_{s}=2\) for the serial-in-time methods and \(P_{p}=2N_{t}P_{s}\) for the parallel-in-time methods. Measured wallclock times per timestep \(T_{p}\) shown for three block solution methods: the star patch, and the composite preconditioner with fixed \(\tau=10^{-5}\) or fixed \(k_{p}=\sqrt{N_{t}}\). The wallclock times per timestep \(T_{s}\) for the serial-in-time methods with the star patch or composite preconditioner are shown at \(N_{t}=1\) in open symbols and plotted with horizontal lines. Maximum speedups for the parallel-in-time methods over each serial-in-time method shown in Table \ref{tab:euler-speedups}.}
    \label{fig:euler-speedup}
\end{figure}

\begin{figure}[t]
    \centering
    \epswitch{
    \includegraphics[width=0.99\linewidth]{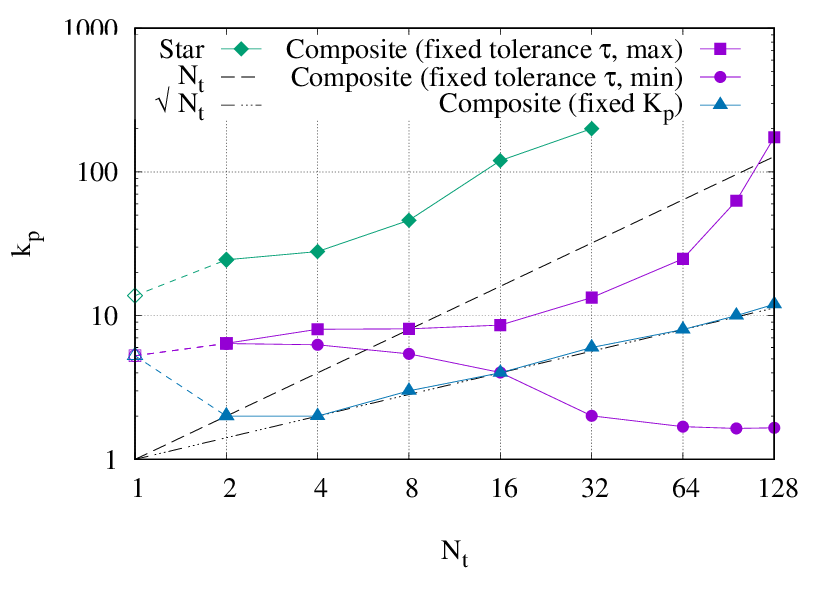}
    }{
    \includegraphics[width=0.99\linewidth]{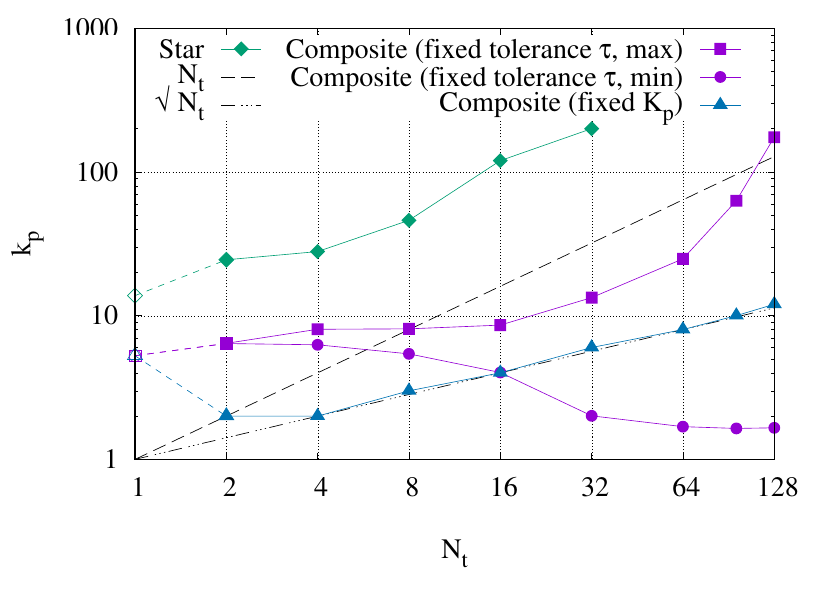}
    }
    \caption{Scaling of the block iteration counts \(k_{p}\) for the compressible Euler equations gravity wave test case at varying window length \(N_{t}\). Three different preconditioning methods shown: star patch; composite preconditioner with fixed \(\tau=10^{-5}\) (for which both maximum and minimum block iteration counts are shown); composite preconditioner with \(k_{p}\) held fixed at \(k_{p}=\sqrt{N_{t}}\) (rounded up). The iteration count for the serial-in-time methods \(k_{s}\) are shown in open symbols at \(N_{t}=1\) (identical for all composite preconditioner methods).}
    \label{fig:euler-gamma}
\end{figure}

\begin{figure}[t]
    \centering
    \epswitch{
    \includegraphics[width=0.99\linewidth]{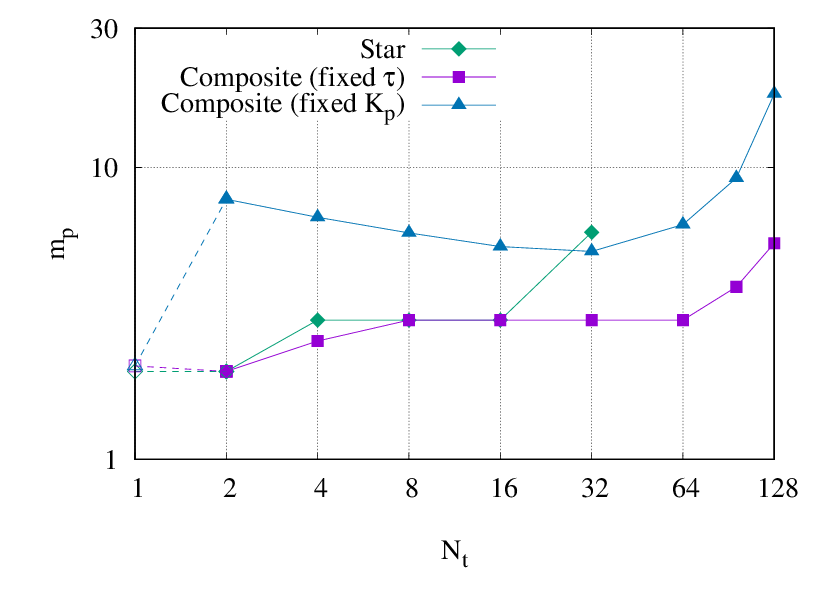}
    }{
    \includegraphics[width=0.99\linewidth]{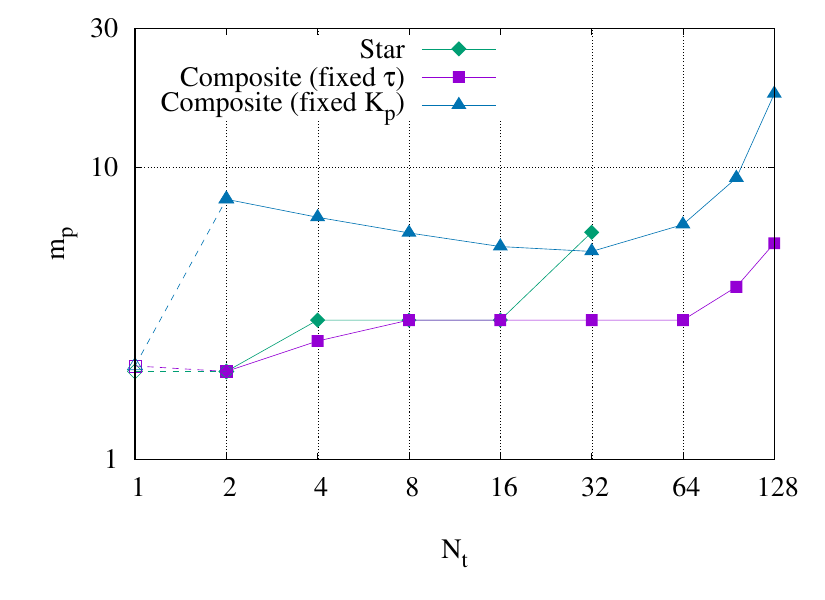}
    }
    \caption{Scaling of the all-at-once outer Krylov iteration counts \(m_{p}\) for the compressible Euler equations gravity wave test case with varying window length \(N_{t}\). Iteration counts shown when the block convergence tolerance is held fixed at \(10^{-5}\) or when the block iterations are held fixed at \(k_{p}=\sqrt{N_{t}}\). The iteration count for the serial-in-time method \(m_{s}\) is shown in open symbols at \(N_{t}=1\).}
    \label{fig:euler-omega}
\end{figure}

The wallclock time achieved with this method is shown in Fig. \ref{fig:euler-speedup} for window lengths \(N_{t}=2-32\), compared to the serial-in-time method (the ``Star patch'' and ``Serial (star)'' curves respectively).
The scaling is underwhelming compared to the previous examples, reaching a maximum speedup of 2.57 at \(N_{t}=8\) (shown in the top left entry in Table \ref{tab:euler-speedups}).
The all-at-once iterations \(m_{p}\) and block iterations \(k_{p}\) are shown in Figs. \ref{fig:euler-omega} and \ref{fig:euler-gamma} respectively.
\(m_{p}\) is quite well-behaved, remaining at 3 for \(N_{t}\leq16\) before jumping to 6 at \(N_{t}=32\).
We would expect \(m_{p}\) to increase more slowly than the previous example because, although the Euler equations are nonlinear, the gravity wave variations around the background state are relatively small.
On the other hand, \(k_{p}\) increases immediately compared to the serial-in-time iteration count \(k_{s}\), and grows quickly with \(N_{t}\).
For this test we had set the maximum number of block iterations to 200 (already an impractically large number).
At \(N_{t}=32\), 18 of the 32 blocks reached the maximum iterations, meaning that the expected tolerance \(\tau\) was not reached, leading to a deterioration in the effectiveness of the circulant preconditioner, and ultimately to the doubling in \(m_{p}\) from \(N_{t}=16\) to \(N_{t}=32\).
We attribute the growth in \(k_{p}\) to the fact that the coefficients \(\psi_{j}\) \eqref{eq:block-coeff-ratio} for some blocks approach the imaginary axis as \(N_{t}\) increases, and that the patch method gives iteration counts that grow linearly with the Courant number.
Although the patch preconditioner is suitable for the serial-in-time method, clearly a different approach is necessary for an effective ParaDiag method.

The patch preconditioner iteration counts grow linearly with the Courant number because the domain of dependence of each application of the preconditioner remains constant, while the physical domain of dependence grows linearly with \(\Delta t\).
Therefore, we propose to compose the patch preconditioner multiplicatively with a second preconditioner.
This second preconditioner is constructed from (\ref{eq:euler}) linearised around a reference state at rest \((\mathbf{u}, \rho, \theta) = (\mathbf{0}, \rho_{ref}, \theta_{ref})\).
An LU factorisation of this matrix is calculated at the start of the simulation and reused throughout the timeseries.%
\footnote{Similar linearisations have previously been used for the inner implicit iterations of semi-implicit time integrators for vertical slice models in \citet{Melvin_Dubal_Wood_Staniforth_Zerroukat_2010} (using a constant in time reference state), and 3D atmospheric models in \citet{Melvin_Benacchio_Shipway_Wood_Thuburn_Cotter_2019} and \citet{Melvin_Shipway_Wood_Benacchio_Bendall_Boutle_Brown_Johnson_Kent_Pring_et_al_2024} (using the previous timestep as the reference state). In these methods the linearised system is solved by reducing down to a pressure Helmholtz equation by Schur factorisation, similarly to the hybridised model in Sect. \ref{sec:examples:lswe}. Applying a similar approach to the current model is left for future work.}
We expect that this composition may improve performance because the two preconditioners are complimentary, in the sense that the linearisation around a state of rest targets the wave terms and provides \textit{global} coupling but lacks the advection terms, whereas the patch smoother handles advection well \textit{locally} but lacks globalisation.
The full block solution strategy is then FGMRES, preconditioned with the composition of: 2 GMRES iterations preconditioned with the linearisation around a reference state at rest, followed by 2 GMRES iterations preconditioned with the column patch preconditioner (we experimented with combinations of 1 or 2 iterations for each inner GMRES method and found (2,2) to have the best performance).

This composition is implemented using PETSc's PCComposite preconditioner.
The linearisation around the reference state of rest is implemented by passing the reference state in the \texttt{appctx} Python dictionary to the \texttt{AuxiliaryRealBlockPC} and \texttt{AuxiliaryComplexBlockPC} preconditioners in the serial- and parallel-in-time methods respectively.
Implementing this solution strategy required only two small modifications to the code compared to using just the patch smoother.
The first was to update the \texttt{'circulant\_block'} PETSc options in the Python options dictionary to use \texttt{'pc\_type': 'composite'} and to add the sub options for each of the two preconditioners described above.
The second was to create a new Firedrake \texttt{Function} to hold the reference state at rest, and pass this into the \texttt{appctx} for the auxiliary preconditioners.

The wallclock time required per timestep for the serial-in-time method, and the parallel-in-time method with increasing window sizes \(N_{t}=2-128\) using the composite block preconditioner, are shown in Fig. \ref{fig:euler-speedup} (the ``Serial (composite)'' and ``Fixed block tolerance'' curves respectively).
The serial-in-time method requires an average of 10.5 block iterations over two Newton iterations per timestep, which is actually slightly slower than the patch preconditioner for this Courant number.
The parallel-in-time method with the composite preconditioner is also slightly slower than with the patch preconditioner for \(N_{t}\leq8\), but continues to scale at larger window sizes.
At \(N_{t}=32\) the speedup over the serial-in-time method with the composite preconditioner is 7.99, and the speedup over the serial-in-time method with the patch preconditioner is 6.14 (second row of Table \ref{tab:euler-speedups}).
The speedup is almost identical at \(N_{t}=64\), before deteriorating at larger \(N_{t}\).
Nonetheless, this demonstrates that the composite preconditioner gives a more performant method, with a \(2.39\times\) improvement on the maximum speedup.

The average number of block inner Krylov iterations for the most difficult block \(k_{p}\), and of all-at-once outer Krylov iterations \(m_{p}\), shown in Figs. \ref{fig:euler-gamma} and \ref{fig:euler-omega} respectively for the composite preconditioner (the ``Fixed block tolerance (max)'' curves), both increase with \(N_{t}\), but much more gradually than with the patch preconditioner.
The block iterations \(k_{p}\) are almost constant until \(N_{t}=16\) before increasing for larger window lengths.
The all-at-once iterations \(m_{p}\) remain almost constant until a longer window length \(N_{t}\leq64\) but a linear increase is seen for \(N_{t}>64\).
This is not due to under converged blocks in the circulant preconditioner as was the case with the patch preconditioner - no blocks reached 200 iterations - but simply due to the longer window lengths.
These combined effects cause the loss of wallclock scaling around \(N_{t}=32\) seen in Fig. \ref{fig:euler-speedup}.

The number of iterations required for \(\tau=10^{-5}\) for the easiest complex-valued block is also shown in Fig. \ref{fig:euler-gamma} (the ``Fixed block tolerance (min)'' curve), and actually decreases with \(N_{t}\).
This block always corresponds to the furthest right point of each line in Fig. \ref{fig:circulant-eigenvalues}, i.e. the block whose coefficients resemble a very small real-valued timestep, which is known to be a favourable condition for iterative solvers.
The trends of the maximum and minimum block iterations suggests a further modification to the solution strategy.
Instead of fixing the required residual drop for the block solves, we could instead fix the number of iterations \(k_{p}\) to be constant for all blocks, somewhere between the maximum and minimum iterations required for a fixed residual drop.
The hope is that ``oversolving`` the easiest blocks will mitigate ``undersolving'' the most difficult blocks enough that the performance of the all-at-once circulant preconditioner does not deteriorate too significantly.
As an example, we fix \(k_{p}=\sqrt{N_{t}}\) (rounding up when \(\sqrt{N_{t}}\) is not an integer).
This particular scaling is somewhat arbitrary, but importantly it scales sublinearly with \(N_{t}\).

The ``Fixed block iterations'' plot in Fig. \ref{fig:euler-omega} shows that the number of all-at-once iterations required is increased by this strategy, by a factor of between 1.7-2.1 for \(N_{t}=8-64\) but rising to a factor of 3.7 at \(N_{t}=128\).
In comparison \(k_{p}\) is reduced by a factor of 3.1 at \(N_{t}=64\).
In Fig. \ref{fig:euler-speedup} we can see that fixing the block iterations in this way improves the speedup for all window lengths, and extends the strong scaling from \(N_{t}=32\) to \(N_{t}=64\) compared to fixing the block tolerance.
At \(N_{t}=64\) a speedup of 11.24 over the serial-in-time method with the composite preconditioner is achieved, and a speedup up over the serial-in-time method with the patch preconditioner of 8.64, constituting an improvement of \(1.41\times\) over using a fixed block tolerance, and \(3.37\times\) greater speedup than the parallel-in-time method with the patch preconditioner (last row in Table \ref{tab:euler-speedups}).
If a more sophisticated solver than LU was used for the linearisation around a state of rest, e.g. a Schur complement reduction, we would expect the speedup of the composite method to improve further.
We clarify that we are not advocating for fixing \(k_{p}\) specifically to \(\sqrt{N_{t}}\) as a general strategy, rather demonstrating that in some cases it is possible to improve overall performance by ``trading'' computational work between the inner block iterations and the outer all-at-once iterations.
More work is needed in the future to investigate leveraging inexact block solves and differing convergence levels between the blocks to improve performance, guided by the estimates of \citet{vcaklovic2023parallel} in their Lemma 3 for the error produced by an inexact application of the circulant preconditioner due to finite machine precision and inexact block solves.
They use this error estimate to derive a method for adapting \(\alpha\) for fixed \(\tau\) and discuss the trade-offs between block tolerance and the outer convergence.
However, varying \(\alpha\) has less impact for nonlinear problems where the dominant error is the time-averaged reference state, and would require recalculating any factorisations in the block solution method at every outer Krylov iteration, so specifying \(\tau\) may be a more effective strategy.

This example shows not only that ParaDiag can produce speedups over serial-in-time methods for the compressible Euler equations, but also that sophisticated solvers can be created and tested using only standard PETSc, Firedrake, and asQ components.
As stated in Sect. \ref{sec:asq} as part of the aims of asQ, this flexibility is essential for developing ParaDiag into a practical method for use in real applications.

\conclusions[Summary and Outlook]\label{sec:conclusions}  %% \conclusions[modified heading if necessary]

\subsection{Summary}

We have introduced a new software library, asQ, for developing the ParaDiag-II family of parallel-in-time methods.
ParaDiag-II provides a parallelisable method for solving multiple timesteps of classical timesteppers, with the implicit \(\theta\)-method used in asQ.
By building on the Firedrake finite element library, Unified Form Language (UFL), and PETSc (the Portable, Extensible Toolkit for Scientific Computation), the design of asQ is intended to improve the productivity of method developers by minimising the effort required to implement a new model or method, thereby maximising the time spent testing said models or methods.

asQ automates the \textit{construction} of the all-at-once system and the associated Jacobian and preconditioners using Firedrake.
Implementing this construction is often a very time consuming step in the development cycle, which here is achieved solely from a user provided high level mathematical description of the finite element model using UFL.
While the \textit{construction} is automated, the user is given full control over the \textit{solution} of the all-at-once system.
This control is predominantly achieved with PETSc options, provided via Python dictionaries or command line arguments.
The options interface enables a wide array of solution methods provided by PETSc, Firedrake, and asQ to be selected and switched between with zero, or minimal, changes to the code.
asQ also supports the user providing additional states or finite element forms to construct more complex solution methods which cannot be specified solely through the options dictionaries.

Space-time parallelism is implemented in MPI, and the vast majority of asQ's API is collective over both time and space, with the user only required to specify how many processors to use for parallelism in the time dimension.
As such, users need very little expertise with parallel computing to be able to run in parallel, and to profile the performance of their methods.

We have used asQ to demonstrate ParaDiag methods on a set of test cases relevant to atmospheric modelling.
The demonstrations have shown, firstly, that asQ achieves the expected convergence rates for ParaDiag and that the parallel implementation scales to over 10,000 processors.
Secondly, we have shown that ParaDiag is capable of achieving excellent speedup over the equivalent serial-in-time method for linear, constant coefficient equations.
For nonlinear problems, some speedup can be achieved, but at a much more modest level due to two factors.
\begin{enumerate}
    \item[1)] The constant-in-time reference state approximation in the preconditioner leads to increasing errors as the window duration grows.
    \item[2)] The complex coefficients of the block linear systems in the preconditioner become less favourable as the number of timesteps in the window grows.
\end{enumerate}
The success of ParaDiag methods for nonlinear and nonconstant-coefficient problems will be determined by how effectively these two issues can be overcome.

\subsection{Outlook}
There are a variety of interesting directions for future research, a few of which we describe here, first for further development of the asQ library, and secondly for development of ParaDiag methods.

\subsubsection{asQ library}

Currently asQ supports only the implicit \(\theta\)-method time integrator, but implementing other time integrators would be a very useful extension.
For example, implicit Runge-Kutta methods have been recently been implemented with Firedrake in the Irksome library \citep{Farrell_Kirby_Marchena-Menendez_2021,Kirby_MacLachlan_2024}.
Irksome generates the UFL and preconditioners for the coupled stage system, which could in turn be used to create the Runge-Kutta all-at-once system and corresponding circulant preconditioner in asQ.

Other diagonalisations could also be implemented, such as the direct diagonalisation of ParaDiag-I \citep{maday2008parallelization}.
The roundoff error for direct diagonalisation grows faster with \(N_{t}\) than in the \(\alpha\)-circulant approach, however, direct diagonalisation may be more competitive for nonlinear problems where the averaging error currently limits ParaDiag to fairly small \(N_{t}\).

There are also open questions on ideal implementation strategies.
For example, some of the complex-valued blocks in the circulant preconditioner are much easier to solve than others.
The parallel efficiency of this step could be improved by load balancing i.e. solving multiple easier blocks on the same communicator.

\subsubsection{ParaDiag methods}

Future improvements of ParaDiag, especially for nonlinear problems, must tackle the two issues identified above: preconditioning the complex-valued blocks; and overcoming the averaging error for long time windows.
We are currently investigating developments in both directions using the nonlinear shallow water equations as a testbed.

Effective preconditioning for the complex-valued blocks will depend on the PDE being solved. 
For atmospheric models, the blocks include both wave propagation and advection terms.
Given the difficulty in solving both high wavenumber Helmholtz equations and high Courant number advection equations, preconditioners that either tackle both components explicitly or reduce the system to a size tractable by direct methods are likely necessary.
This could entail composite preconditioners, similar to that used for the compressible Euler equations above, or preconditioners based on Schur complement factorisations.

Diagonalisation in time requires the preconditioner to be a sum of Kronecker products, which in turn requires a constant-in-time spatial Jacobian.
As such, overcoming the averaging error will require alternate solution strategies, for example splitting the timeseries into smaller chunks which can each be preconditioned with a circulant preconditioner, thereby decreasing the averaging error in each chunk.
At the linear level, this could resemble a block Jacobi preconditioner for the all-at-once system (see the \texttt{SliceJacobiPC} in Appendix \ref{app:other-asq}).
At the nonlinear level, one of the simplest options is to parallelise over the quasi-Newton iterations of several all-at-once systems, which is essentially a pipelined nonlinear Gauss-Seidel method, for which we have some promising preliminary results.

ParaDiag methods is currently a very active area in the parallel-in-time field, with many possible research directions.
The library presented here provides a sandbox for implementing and testing these methods, enabling faster development of new approaches as well as scalable performance demonstrations.

%% The following commands are for the statements about the availability of data sets and/or software code corresponding to the manuscript.
%% It is strongly recommended to make use of these sections in case data sets and/or software code have been part of your research the article is based on.

\codedataavailability{%
All code used in this manuscript is free and open source.
The asQ library is available at https://github.com/firedrakeproject/asQ and is released under the MIT license.
The data in this manuscript was generated using the Python scripts in (\cite{hope_collins_2024_asq_scripts}, also available in the asQ repository), the Singularity container in \citep{hope_collins_2024_asq_container}, and the versions of Firedrake, PETSc, and their dependencies in \citep{firedrake_zenodo_2024_asq}.
}

%%% TODO: add these back in for the journal submission:
%%% Finite element forms - appendix

\appendix
\section{Finite element forms}\label{app:finite-element-forms}

In this appendix the finite element weak forms used for each example in Sect. \ref{sec:examples} are given.
First, we define some common nomenclature.
The spatial coordinate is \(\vec{x}\in\Omega\subset\mathbb{R}^{d}\), where \(\Omega\) is the domain of interest with boundary \(\partial\Omega\), and \(d\) is the spatial dimension.
\(L^{2}\left(\Omega\right)\) is the space of scalar functions on \(\Omega\) with finite \(L^{2}\) norm, and \(H_{\text{div}}\left(\Omega\right)\subset\left(L^{2}\left(\Omega\right)\right)^{d}\) is the space of vector valued functions with divergence also in \(L^{2}\left(\Omega\right)\).
The integration measures used are: \(d\vec{x}\) over element interiors; \(d\vec{s}\) over element facets \(\Gamma\) in the interior of \(\Omega\); \(d\vec{S}\) over element facets on the boundary \(\partial\Omega\).
\(|\vec{v}|\) is the (local) norm of \(\vec{v}\).
If \(\varphi^{-}\) and \(\varphi^{+}\) are the values of \(\varphi\) on either side of a facet, then \([\![\varphi]\!]=\varphi^{+}-\varphi^{-}\) is the jump over the facet, and \(\{\varphi\}=\left(\varphi^{-}+\varphi^{+}\right)/2\) is the average over the facet.
\(\tilde{\varphi}\) is the value of \(\varphi\) on the upwind side of the facet according to a velocity \(\vec{u}\).

\subsection{Scalar advection}\label{app:scalar-advection}

The scalar advection \eqref{eq:scalar-advection} with advecting velocity \(\vec{u}\) is solved by finding \(q \in V \subset L^{2}\left(\Omega\right)\) such that
\begin{align}\label{eq:app:scalar-advection}
    \int_{\Omega} \left(q\partial_{t}\phi - q\nabla\cdot(\phi\vec{u}) \right)\text{d}\vec{x} 
    + \int_{\Gamma} \tilde{q}[\![\phi\vec{u}\cdot\vec{n}]\!]\text{d}\vec{s}
    = 0\quad \forall \phi \in V
\end{align}
In Sect. \ref{sec:examples:advection} we choose \(V=\text{DG}_{1}\), the space of piecewise linear polynomials.

\subsection{Linear shallow water}\label{app:lswe}

The linear shallow water equations on the rotating sphere \eqref{eq:lswe} are solved by finding \(\left(\vec{u}, h\right)\in V_{\vec{u}}\times V_{h}\subset H_{\text{div}}\left(\Omega\right) \times L^{2}\left(\Omega\right)\), such that
\begin{align}\label{eq:app:lswe}
    \int_{\Omega}\left(\vec{v}\cdot\partial_{t}\vec{u} + \vec{v}\cdot(f\vec{u}^{\perp}) - gh\nabla\cdot\vec{v}\right) \text{d}\vec{x} = 0 \quad &\forall \vec{v} \in V_{u}, \\
    \int_{\Omega}\left(\phi\partial_{t}h + H\phi\nabla\cdot\vec{u}\right)\text{d}\vec{x} = 0 \quad &\forall \phi \in V_{h}
\end{align}
where \(f=2\omega z/R\) is the Coriolis parameter, \(\omega\) and \(R\) are the rotation rate and radius of the sphere respectively, \(z\) is the vertical coordinate from the centre of the sphere, \(g\) is gravity, \(H\) is the constant mean depth, and \(\vec{u}^{\perp}=\vec{u}\times\hat{\vec{k}}\), where \(\hat{\vec{k}}\) is the unit vector normal to the surface of the sphere.
In Sect. \ref{sec:examples:lswe} we choose \(V_{h}=\text{DG}_{k-1}\), and \(V_{\vec{u}}=\text{BDM}_{k}\) the Brezzi-Douglas-Marini elements \citep{Brezzi_Douglas_Marini_1985}, with \(k=2\).

\subsection{Nonlinear shallow water}\label{app:swe}
The nonlinear shallow water equations on the rotating sphere \eqref{eq:swe} are solved by finding \(\left(\vec{u}, h\right)\in V_{\vec{u}}\times V_{h}\subset H_{\text{div}}\left(\Omega\right) \times L^{2}\left(\Omega\right)\), such that
\begin{align}\label{eq:app:swe}
    \int_{\Omega}\left((\vec{v}\cdot\partial_{t}\vec{u} + f\vec{v}\cdot\vec{u}^{\perp} - g(h+b)\nabla\cdot\vec{v}\right)\text{d}\vec{x} \\ \nonumber
    - \int_{\Omega}\left(\nabla\times(\vec{v}\times\vec{u})\times\vec{u} - \frac{1}{2}|\vec{u}|^{2}(\nabla\cdot v)\right) \text{d}\vec{x} \\ \nonumber
    - \int_{\partial\Omega}\left(\vec{n}\times(\vec{u}\times\vec{v})\cdot\tilde{\vec{u}}\right) \text{d}\vec{s} = 0 \quad &\forall \vec{v} \in V_{u},\\
    \int_{\Omega}\left(\phi\partial_{t}h - h\vec{u}\cdot\nabla\phi\right)\text{d}\vec{x} + \int_{\partial\Omega}[\![\phi]\!]\tilde{h}\tilde{\vec{u}}\cdot\vec{n}\text{d}\vec{s} = 0 \quad &\forall \phi \in V_{h}
\end{align}
where \(f\), \(g\), and \(\vec{u}^{\perp}\) are the same as in \ref{app:lswe}, and \(b\) is the topography.
In Sect. \ref{sec:examples:swe} we choose \(V_{\vec{u},h}\) identical to in \ref{sec:examples:lswe}.

\subsection{Vertical slice compressible flow}\label{app:vertical-slice}

The Euler equations restricted to a vertical slice \eqref{eq:euler} are  solved by finding \((\vec{u}, \theta, \rho)\in V_{\vec{u}}\times V_{\theta}\times V_{\rho}\) such that:
\begin{align}
    \int_{\Omega}\left(\vec{v}\cdot\partial_{t}\vec{u} + f\vec{v}\cdot\vec{u}^{\perp} + \vec{v}\cdot\hat{\vec{k}}g + \mu\vec{v}\cdot\hat{\vec{k}}\vec{u}\cdot\hat{\vec{k}}\right)\text{d}\vec{x} \\ \nonumber
    + \int_{\Omega}\left(\nabla_{h}\times\left(\vec{v}\times\vec{u}\right)\times\vec{u} - \nabla\cdot\vec{v}\frac{1}{2}|\vec{u}|^{2}\right)\text{d}\vec{x} \\ \nonumber
    + \int_{\partial\Omega}[\![\vec{n}\times\left(\vec{u}\times\vec{v}\right)]\!]\cdot\vec{\tilde{u}}\text{d}\vec{s} \\ \nonumber
    - \int_{\Omega}\nabla_{h}\cdot\left(\vec{v}\theta\right)c_{p}\Pi\text{d}\vec{x} + \int_{\Gamma_{v}}[\![\vec{n}\cdot\vec{v}\theta]\!]c_{p}\{\Pi\}\text{d}\vec{S} = 0 \quad &\forall \vec{v} \in V_{u},\\
    \int_{\Omega}\left(q\partial_{t}\theta - q\vec{u}\cdot\nabla_{h}\theta\right)\text{d}\vec{x} \\ \nonumber
    + \int_{\Gamma_{v}}[\![q\vec{u}\cdot\vec{n}]\!]\tilde{\theta}\text{d}\vec{S} - \int_{\Gamma_{v}}[\![q\theta\vec{u}\cdot\vec{n}]\!]\text{d}\vec{S} \\ \nonumber
    + \int_{\Gamma}C_{0}h^{2}|\vec{u}\cdot\vec{n}|[\![\nabla_{h}q]\!]\cdot[\![\nabla_{h}\theta]\!]\text{d}\vec{S} = 0 \quad &\forall q \in V_{\theta},\\
    \int_{\Omega}\left(\phi\partial_{t}\rho - \rho\vec{u}\cdot\nabla_{h}\phi\right)\text{d}\vec{x} + \int_{\Gamma}[\![\phi\vec{u}\cdot\vec{n}]\!]\tilde{\rho}\text{d}\vec{S} = 0 \quad &\forall \phi \in V_{\rho}
\end{align}
This discretisation is described in detail by \citet{Cotter_Shipton_2023}.
\(\hat{\vec{k}}\) is the vertical unit vector, \(g\) is again the gravity, \(\mu\) is a spatially varying viscosity parameter to damp reflections at the top of the domain, \(c_{p}\) is the constant pressure specific heat, \(C_{0}\) is a stabilisation constant, \(h\) is a measure of facet edge length, \(\nabla_{h}\) is the gradient evaluated locally in each cell, and \(\Gamma_{v}\) is the set of vertical facets.
Each function space \(V=V^{h}\otimes V^{v}\) is a tensor product of spaces defined in the horizontal \(V^{h}\) and vertical \(V^{v}\) directions.
The density space is \(V_{\rho}=\text{DG}^{h}_{k-1}\otimes\text{DG}^{v}_{k-1}\) i.e. fully discontinuous.
The temperature space is \(V_{\theta}=\text{DG}^{h}_{k-1}\otimes\text{CG}^{v}_{k}\) i.e. discontinuous in the horizontal, continuous in the vertical.
The velocity space is defined for the horizontal \(u\)  and vertical \(w\)  velocity components separately \(V_{\vec{u}}=V_{u}\oplus V_{w}\), where \(V_{w}=V_{\theta}\), and \(V_{u}=\text{CG}^{h}_{k}\otimes \text{DG}^{v}_{k-1}\) i.e. continuous in the horizontal and discontinuous in the vertical.
In Sect. \ref{sec:examples:euler} we choose \(k=2\).

\section{Other asQ components}\label{app:other-asq}

We give a brief description of some asQ components further to those described in Sect. \ref{sec:asq:components}.
These components were not described in Sect. \ref{sec:asq:components} because they are not required for the core ParaDiag method and were not used for the examples in Sect. \ref{sec:examples}.
However, they are useful for developing further ParaDiag methods so we include them here.

\subsubsection*{\texttt{LinearSolver}}
Given an \texttt{AllAtOnceForm}, the \texttt{LinearSolver} sets up a PETSc KSP for a linear system where the matrix is the \texttt{AllAtOnceJacobian} of the given form.
The linear system can then be solved with a given \texttt{AllAtOnceCofunction} for the right hand side.
This is different from the \texttt{AllAtOnceSolver} in that it does not automatically include the initial conditions (\ref{eq:theta-rhs}) in the right hand side, so the solution from a \texttt{LinearSolver} is not a timeseries (unless the \texttt{AllAtOnceCofunction} has been calculated from the initial conditions).
The \texttt{LinearSolver} has two main uses. First, constructing preconditioners from the all-at-once Jacobian (e.g. the \texttt{SliceJacobiPC} described below), and second, testing properties of linear solution strategies by setting specific right hand sides (e.g. Fourier modes).
If the \texttt{AllAtOnceForm} is nonlinear, then the \texttt{AllAtOnceJacobian} is the linearisation of the \texttt{AllAtOnceForm} around its \texttt{AllAtOnceFunction}.

\subsubsection*{\texttt{JacobiPC}}
The classical point Jacobi preconditioner approximates a matrix by a diagonal matrix.
Block Jacobi preconditioners extend this idea by using a block-diagonal approximation, where each block is (an approximation of) the block of the original matrix that couples a set of \(m\) DoFs.
E.g. the elements of the block \(B\) corresponding to the DoFs \(\{l,l+1,\dots,l+m-1\}\) of a matrix \(A\) are \(B_{i,j}=A_{l+i,l+j},\; i,j\in\{0,1,\dots,m-1\}\).

The \texttt{JacobiPC} class is a block Jacobi preconditioner for the all-at-once Jacobian (\ref{eq:nonlinear-jacobian}), with \(N_{t}\) blocks each corresponding to the DoFs of a single timestep, i.e. \(m=N_{x}\) and the \(j\)-th block is \(M/\Delta t + \theta\nabla_{u}f(u^{j},t^{j})\).
Unlike the circulant preconditioner, each block can be linearised around a different state so can exactly match the diagonal blocks in the Jacobian.
The construction of the blocks and the solver parameters used can be customised to those of the \texttt{CirculantPC}.
However, because the all-at-once Jacobian is block lower triangular, a Krylov method with this preconditioner requires \(N_{t}\) iterations before achieving any substantial drop in the residual \citep{wathen_observations_2022}.

\subsubsection*{\texttt{SliceJacobiPC}}
The \texttt{SliceJacobiPC} class is a second block Jacobi preconditioner for the all-at-once Jacobian (\ref{eq:nonlinear-jacobian}), however here we refer to the blocks of the preconditioner as ``slices'' to avoid confusion with the use of ``blocks'' in the rest of the paper.
The \texttt{SliceJacobiPC} has \(N_{s}\) slices, each constructed from \(k\) (consecutive) timesteps where \(kN_{s}=N_{t}\), i.e. \(m=kN_{x}\) and each block is the all-at-once Jacobian for \(k\) timesteps.
Each slice is then (approximately) inverted with a \texttt{LinearSolver}.
Previously, PETSc's solver composition enabled taking solver options for the serial-in-time method and using them for the inner blocks of \texttt{CirculantPC} and \texttt{JacobiPC}.
In the same way, the solver composition enables taking solver options for an \texttt{AllAtOnceSolver} and using them for each slice of a \texttt{SliceJacobiPC}.
For example we could approximate each slice with a separate \texttt{CirculantPC} - each using the time-averaged state of their own slice rather than of the entire timeseries.

\noappendix       %% use this to mark the end of the appendix section. Otherwise the figures might be numbered incorrectly (e.g. 10 instead of 1).

%% Regarding figures and tables in appendices, the following two options are possible depending on your general handling of figures and tables in the manuscript environment:

%% Option 1: If you sorted all figures and tables into the sections of the text, please also sort the appendix figures and appendix tables into the respective appendix sections.
%% They will be correctly named automatically.

%% Option 2: If you put all figures after the reference list, please insert appendix tables and figures after the normal tables and figures.
%% To rename them correctly to A1, A2, etc., please add the following commands in front of them:

%\appendixfigures  %% needs to be added in front of appendix figures

%\appendixtables   %% needs to be added in front of appendix tables

%% Please add \clearpage between each table and/or figure. Further guidelines on figures and tables can be found below.

\authorcontribution{
Conceptualisation: JHC, AH, WB, LM, CC.
Data curation: JHC.
Formal analysis: JHC, CC.
Funding aquisition: CC.
Investigation: JHC.
Methodology: JHC, CC.
Project administration: CC.
Software: JHC, AH, WB, LM, CC.
Supervision: CC.
Validation: JHC, WB, CC.
Visualisation: JHC.
Writing - original draft preparation: JHC, CC.
Writing - review \& editing: JHC, AH, WB, LM, CC.
} %% this section is mandatory

\competinginterests{The authors declare that they have no conﬂict of interest.} %% this section is mandatory even if you declare that no competing interests are present

%\disclaimer{TEXT} %% optional section

\begin{acknowledgements}
This work was supported by the Engineering and Physical Sciences Research Council (EP/W015439/1 \& EP/R029628/1); the Natural Environment Research Council (NE/R008795/1); UK Research \& Innovation and the UK Met Office through the ExCALIBUR programme (SPF EX20-8).
This work used the ARCHER2 UK National Supercomputing Service (https://www.archer2.ac.uk).
\end{acknowledgements}

%% REFERENCES

%% The reference list is compiled as follows:

% \begin{thebibliography}{}
% 
% \bibitem[AUTHOR(YEAR)]{LABEL1}
% REFERENCE 1
% 
% \bibitem[AUTHOR(YEAR)]{LABEL2}
% REFERENCE 2
% 
% \end{thebibliography}

%% Since the Copernicus LaTeX package includes the BibTeX style file copernicus.bst,
%% authors experienced with BibTeX only have to include the following two lines:
%%
\bibliographystyle{copernicus.bst}
\bibliography{asQ_paper.bib}

\end{document}